\documentclass[10pt,a4paper]{amsart}
\usepackage{amsfonts,layout}

\newcount\minute
\newcount\hour
\def\currenttime{%
	\minute\time
	\hour\minute
	\divide\hour60
	\the\hour:\multiply\hour60\advance\minute-\hour\the\minute}
\def\draftnote{{\it \today \quad  \currenttime \hfill  tex-file :   \jobname}}

\catcode`@=11
\@addtoreset{equation}{section}

\catcode`@=12
\newtheorem{Theorem}{Theorem}[section]
\newtheorem{Definition}{Definition}[section]
\newtheorem{Proposition}{Proposition}[section]
\newtheorem{Lemma}{Lemma}[section]

\newtheorem{Remark}{Remark}[section]

\newtheorem{Hyp.}{Hyp.}[section]

\begin{document}

\title[]{The cost of controlling strongly degenerate parabolic equations}

\author{P. Cannarsa} 
\address{Dipartimento di Matematica, Universit\`a di Roma "Tor Vergata",
Via della Ricerca Scientifica, 00133 Roma, Italy}
\email{cannarsa@mat.uniroma2.it}

\author{P. Martinez} 
\address{Institut de Math\'ematiques de Toulouse; UMR 5219, Universit\'e de Toulouse; CNRS \\ 
UPS IMT F-31062 Toulouse Cedex 9, France} \email{patrick.martinez@math.univ-toulouse.fr}


\author{J. Vancostenoble} 
\address{Institut de Math\'ematiques de Toulouse; UMR 5219, Universit\'e de Toulouse; CNRS \\ 
UPS IMT F-31062 Toulouse Cedex 9, France} \email{Judith.Vancostenoble@math.univ-toulouse.fr}


\thanks{The first author was partly supported by the University of Roma Tor Vergata (Consolidate the Foundations 2015) and Istituto Nazionale di Alta Matematica (GNAMPA 2017 Reseach Projects). The second author was partly supported by Istituto Nazionale di Alta Matematica (GNAMPA 2017 Reseach Projects)}

\begin{abstract}
We consider the typical one-dimensional strongly degenerate parabolic operator $Pu= u_t - (x^\alpha u_x)_x$ with $0<x<\ell$ and $\alpha\in(0,2)$,
controlled either by a boundary control acting at $x=\ell$, or by a locally distributed control.
Our main goal is to study the dependence of the so-called controllability cost needed to drive an initial condition to  rest with respect to the degeneracy parameter $\alpha$. We prove that the control cost blows up with an explicit exponential rate, as $e^{C/((2-\alpha)^2 T)}$, when  $\alpha \to 2^-$ and/or $T\to 0^+$.

Our analysis builds on earlier results and methods (based on functional analysis and complex analysis techniques) developed by several authors such as Fattorini-Russel, Seidman, G\"uichal, Tenenbaum-Tucsnak and Lissy for the classical heat equation. In particular, we use the moment method and related constructions of suitable biorthogonal families, as well as new fine properties of the Bessel functions $J_\nu$ of large order $\nu$ (obtained by ordinary differential equations techniques).

\end{abstract}
\maketitle

\section{Introduction}

\subsection{Presentation of the problem and of the main results} \hfill

The aim of this paper is to study the null controllability cost for the typical 1D degenerate parabolic operator 
\begin{equation}
\label{intro-heat-deg}
Pu=u_t - (x^\alpha u_x)_x\qquad(x\in(0,1),\ t>0)
\end{equation}
under the action of a boundary control $H$:
\begin{equation}
\label{eq-intro-bd}
\begin{cases}
u_t - (x^\alpha u_x)_x =0 , & \qquad x\in(0,1),\ t>0,\\
(x^\alpha u_x) (0,t)=0, & \qquad t>0, \\
u(1,t)= H(t),  & \qquad t>0, \\
u(x,0)=u_0(x), &\qquad x\in(0,1),
\end{cases}
\end{equation}

\noindent and under the action of a locally distributed control $h$:
\begin{equation}
\label{eq-intro}
\begin{cases}
u_t - (x^\alpha u_x)_x =h (x,t) \chi_{[a,b]}(x)  & \qquad x\in(0,1),\ t>0,\\
(x^\alpha u_x) (0,t)=0, & \qquad t>0, \\
u(1,t)=0,  & \qquad t>0, \\
u(x,0)=u_0(x), &\qquad x\in(0,1) .
\end{cases}
\end{equation}
In \cite{sicon2008}, we established the following property:
\begin{quote}
given $\alpha \geq 1$, $T>0$, $0<a<b<1$, then, for any $u_0 \in L^2(0,1)$, problem \eqref{eq-intro}
admits a control $h \in L^2 ((a,b)\times (0,T))$ that drives the solution to $0$ in time $T>0$
{\it if and only if} $\alpha < 2$.
\end{quote}
In the same way,
\begin{quote}
given $\alpha \geq 1$, $T>0$, then, for any $u_0 \in L^2(0,1)$, problem \eqref{eq-intro-bd}
admits a control $H \in L^2 (0,T)$ that drives the solution to $0$ in time $T>0$
{\it if and only if} $\alpha < 2$.
\end{quote}
 The aim of this paper is to understand the behavior of the cost of control as: 
\begin{itemize}
\item  $\alpha \to 2^-$ (since $\alpha=2$ is the threshold for null controllability) ,
\item and/or $T\to 0^+$, an issue related to the so-called 'fast control problem'.
\end{itemize}
It is well-known that the cost of control blows up when $T\to 0^+$ (at least for nondegenerate parabolic equations, as we will recall in the following), 
and it is expected to blows up when $\alpha \to 2^-$. In this work we will prove precise upper and lower bounds for this blow-up:
denoting by  $u^{(h)}$ the solution of \eqref{eq-intro} for any given $h \in L^2 ((a,b)\times (0,T))$,
we  prove that the null controllability costs, defined as
$$ C_{ctr-bd}(\alpha,T) := \sup _{\Vert u_0 \Vert _{L^2(0,1)}=1} \inf \{\Vert H \Vert _{H^1 (0,T)}, u^{(H)}(T)=0 \} , $$
and
$$  C_{ctr-loc}(\alpha,T) := \sup _{\Vert u_0 \Vert _{L^2(0,1)} =1} \inf \{\Vert h \Vert _{L^2((0,1)\times ((0,T))}, u^{(h)}(T)=0 \}, $$
blow up \begin{itemize}
\item as $\alpha \to 2^-$, 
\item and/or as $T\to 0^+$, 
\end{itemize}
at a precise speed: there exist positive constants $C, C'$ such that

$$ 
e^{-\frac{1}{C} \frac{1}{(2-\alpha)^{4/3}}(\ln \frac{1}{2-\alpha} + \ln \frac{1}{T})}  e^{\frac{C}{T (2-\alpha) ^2}}
\leq C_{ctr-bd}(\alpha,T) \leq 
e^{\frac{C'}{T (2-\alpha) ^2 }} ,
$$
and, in a similar way,
$$ 
e^{-\frac{1}{C} \frac{1}{(2-\alpha)^{4/3}}(\ln \frac{1}{2-\alpha} + \ln \frac{1}{T})}  e^{\frac{C}{T (2-\alpha) ^2}}
\leq C_{ctr-loc}(\alpha,T) \leq 
e^{\frac{C'}{T (2-\alpha) ^2 }} .
$$
(See precise statements in Theorems \ref{thm-cost-bd-low}, \ref{thm-cost-bd-up}, \ref{thm-cost-ld-low} and \ref{thm-cost-ld-up}.)


\subsection{Relation to literature} \hfill

This question of the cost of null controllability when some parameter comes into play has been studied for several equations and in several situations:

\begin{itemize}
\item the 'fast control problem', that is, the cost of null controllability with respect to time $T$ as $T\to 0^+$, has been investigated for the heat operator 
\begin{equation}
\label{intro-heat-fast}
Pu= u_t - \Delta u
\end{equation}
(with a boundary or localized control) and the Schr\"odinger equation by several authors, see, in particular, the works by Seidman et al \cite{Seidman, Seid-Avdon}, G\"uichal \cite{Guichal}, Miller \cite{Miller1, Miller2, Miller3, Miller4}, Tenenbaum and Tucsnak \cite{Tucsnak, Tucsnak2}, and the more recent papers by Lissy \cite{Lissy1} (for dispersive equations) and Benabdallah et al \cite{Assia-cost} (for parabolic systems);
\item the 'vanishing viscosity limit', that is the cost of null controllability of a heat operator with the addition of a transport term when the diffusion coefficient goes to zero:
\begin{equation}
\label{intro-heat-van}
P_\varepsilon u=u_t - \varepsilon u_{xx} + M u_x 
\end{equation}
(again with a boundary or localized control)  has been investigated by Coron and Guerrero \cite{Coron}, Guerrero and Lebeau \cite{Guerrero-Lebeau}, Glass \cite{Glass}, Glass and Guerrero \cite{Glass-Guerrero}, and Lissy \cite{Lissy2};
\item the 1D degenerate parabolic equation, controlled by a boundary control
acting at the degeneracy point (and $\alpha \to 1^-$, $1$ being the threshold value of well-posedness in this case, see \cite{cost-weak}).
\end{itemize}


\subsection{Description of the method and connection with the literature} \hfill

For the proof of our results we  follow the classical strategy which consists in:
\begin{itemize} 
\item the spectral analysis of the associated stationary operator (see Proposition \ref{*prop-vp}) in order to determine the eigenvalues and eigenfunctions of our problem by typical ODE techniques,
\item the use of the moment method, that was developed in the seminal papers by Fattorini and Russell \cite{FR1,FR2}, to give, at least formally, a sequence of relations satisfied by the desired control,
\item the construction and the properties of suitable biorthogonal families which are the main tool (at this point we will use two extensions of the results 
of Seidman-Avdonin-Ivanov \cite{Seid-Avdon} and G\"uichal \cite{Guichal}, that we proved in \cite{cost-weak} and \cite{CMV-biortho-general}), 
\item the construction of suitable controls, mainly based on the behavior of the eigenfunctions of the spectral problem in the control region.
\end{itemize}

Hence a starting point is the study of the spectral problem. In the context of degenerate parabolic equations, it is classical that the eigenfunctions of the problem are expressed in terms of Bessel functions of order $\nu_\alpha = \frac{\alpha -1}{2-\alpha}$, and the eigenvalues in terms of the zeros of these Bessel functions, see Kamke \cite{Kamke}. This was a crucial observation in the work by Gueye~\cite{Gueye}, where the null controllability of the degenerate heat equation for $\alpha \in [0,1)$ was addressed for the first  time when the control acts at the degeneracy point, and in \cite{cost-weak} where we obtained optimal bounds for the cost of control for such a problem. 
For strongly degenerate parabolic equations, an additional source of difficulty is that the order of the useful Bessel functions blows up as $\alpha \to 2^-$. To cope with such difficulties several classical results from Watson \cite{Watson} and Qu-Wong \cite{QuWong} will be needed. 

It turns out  that there is a common phenomenon in the classical fast control problem (\eqref{intro-heat-fast} when $T \to 0^+$), the vanishing viscosity problem (\eqref{intro-heat-van} when $\varepsilon \to 0 ^+$), and the 
 null controllability of the degenerate parabolic equation \eqref{intro-heat-deg} when the degeneracy parameter approaches its critical value: {\it the eigenvalues concentrate when parameters go to their critical values}. Such a concentration phenomenon can be observed:
\begin{itemize}
\item for the vanishing viscosity problem, in \cite{Coron};
\item  for degenerate parabolic equations, once the eigenvalues have been computed, in Lemma \ref{lem-concentration};
\item  for the classical heat equation and the fast control problem, once time has been renormalized to a fixed value, see Remark~\ref{re:conc_heat}.
\end{itemize}
This common feature is the key point in understanding the behavior of the control in every context. Indeed, the construction of suitable biorthogonal families is strongly related to gap properties: the gap $\lambda _{n+1}-\lambda _{n} \to 0$ when the degeneracy parameter goes to its critical limit, and the speed at which it goes to zero govern the upper and lower estimates for the associated biorthogonal families (since the norm of such biorthogonal families involve large products of the inverse of such differences), hence for the null controllability cost.

In the context of degenerate parabolic equations, in order to obtain optimal bounds,  
\begin{itemize}
\item we refine classical results providing sharp gap estimates for the zeros of Bessel functions of large order, see Lemma \ref{lem-Sturm},
\item we will combine these gap estimates with  some recent results \cite{cost-weak, CMV-biortho-general}
that complete classical results of Fattorini-Russel \cite{FR2}, obtaining explicit and precise (upper and lower) estimates for biorthogonal families, even in short time, under some gap conditions, namely:
$$ \begin{cases} 
\forall n\geq 1, \quad \gamma _{min} \leq \sqrt{\lambda _{n+1}} - \sqrt{\lambda _{n}} \leq \gamma _{max} \\
\forall n\geq N^*, \quad \gamma _{min}^* \leq \sqrt{\lambda _{n+1}} - \sqrt{\lambda _{n}} \leq \gamma _{max}^*
\end{cases} $$
these gap conditions are a little more general than the asymptotic development of the eigenvalues used in Tenenbaum-Tucsnak \cite{Tucsnak} and Lissy \cite{Lissy1, Lissy2}:
$$ \lambda _n = rn^2 + O(n),$$
but, which is more important, they allow us to obtain precise estimates for the biorthogonal family when some parameter comes into play, as it happens here or in 2D problems (see, e.g., \cite{Karine-Grushin1, FR3}); (the proof of the general results obtained in \cite{cost-weak, CMV-biortho-general} concerning biorthogonal families  is based on some complex and hilbertian analysis techniques developped by Seidman-Avdonin-Ivanov \cite{Seid-Avdon}, G\"uichal \cite{Guichal} and the adjonction of some well-chosen parameter, inspired from Tenenbaum-Tucsnak \cite{Tucsnak});
\item we complete the analysis of the asymptotic behavior of Bessel functions of large order, see Proposition \ref{prop-fctpr}; this issue is related to the so-called 'transition zone' (see Watson \cite{Watson} and Krasikov \cite{Krasikov}) and recent results of Privat-Tr\'elat-Zuazua \cite{Privat2015} p. 957, even though we prove Proposition \ref{prop-fctpr} directly by estimating the norm of the solution of a second-order differential equation depending on a large parameter. 
\end{itemize}


\subsection{Plan of the paper} \hfill

The plan of the paper is the following.
\begin{itemize}
\item In section \ref{s-results-1D}, we state our main results concerning the null controllability costs (Theorems \ref{thm-cost-bd-low}-\ref{thm-cost-ld-up}) and the spectral properties of the problem (see Propositions \ref{*prop-vp} and \ref{prop-fctpr}).
\item In section \ref{*s3}, we recall the main properties of Bessel functions, and prove Propositions \ref{*prop-vp} and \ref{prop-2.2}.
\item In section \ref{sec-moment}, we establish useful identities by the moment method.
\item In section \ref{sec-preuveThm1}, we prove Theorem \ref{thm-cost-bd-low}; the proof is based on a recent result \cite{CMV-biortho-general} based on hilbertian techniques developed by G\"uichal \cite{Guichal}, the concentration of the eigenvalues (Lemma \ref{lem-concentration}) and a precised form of a classical property concerning the gap of the zeros of the Bessel function of large parameter (Lemma \ref{lem-Sturm}).
\item In section \ref{sec-Thm2} we prove Theorem \ref{thm-cost-bd-up}; the proof is based on the construction of some suitable biorthogonal family, based on \cite{cost-weak} and inspired by the construction of Seidman et al \cite{Seid-Avdon}.
\item In section \ref{sec-Thm4}, we prove Theorem \ref{thm-cost-ld-up}, which is a direct consequence of the biorthogonal family constructed in section \ref{sec-Thm2}, assuming Proposition \ref{prop-fctpr}.
\item In section \ref{sec-Thm3}, we prove Theorem \ref{thm-cost-ld-low}; the proof is based on energy methods and uses Theorem \ref{thm-cost-bd-low}.
\item  In section \ref{sec-complfp}, we study the eigenfunctions in the control region, proving Proposition \ref{prop-fctpr}; the proof is based on ODE techniques.
\end{itemize}


\section{The cost of null controllability: main results}
\label{s-results-1D}


\subsection{The controllability problems} \hfill

We study the cost of null controllability of a degenerate parabolic equation, using either a  boundary control acting at the non degeneracy point or a locally distributed control. More precisely, we fix $\ell >0$, $\alpha \geq 1$, $T>0$, and for any $u_0 \in L^2(0,\ell)$, we wish to find a control $H$ such that the solution of 

\begin{equation}
\label{*pbm-controle1}
\begin{cases}
u_t - (x^\alpha u_x)_x =0 , & \qquad x\in(0,\ell),\ t>0,\\
(x^\alpha u_x) (0,t)=0, & \qquad t>0, \\
u(\ell,t)= H(t),  & \qquad t>0, \\
u(x,0)=u_0(x), &\qquad x\in(0,1),
\end{cases}
\end{equation}
also satisfies $u(\cdot, T) =0.$

Similarly, given $0 < a<b< \ell$ and for any $u_0 \in L^2(0,\ell)$, we wish to find a control $h$ such that the solution of 
\begin{equation}
\label{*pbm-controle2}
\begin{cases}
u_t - (x^\alpha u_x)_x =h (x,t) \chi_{[a,b]}(x)  & \qquad x\in(0,\ell),\ t>0,\\
(x^\alpha u_x) (0,t)=0, & \qquad t>0, \\
u(\ell,t)=0,  & \qquad t>0, \\
u(x,0)=u_0(x), &\qquad x\in(0,\ell),
\end{cases}
\end{equation}
also satisfies $u(\cdot, T) =0.$

\noindent In space dimension 1, these two problems are very close. From \cite{sicon2008}, we know that such controls exist if and only if $\alpha \in [1,2)$. 


\subsection{Functional setting and well-posedness}  \hfill

\subsubsection{Functional setting and well-posedness for a locally distributed control}  \hfill

For  $1 \leq \alpha < 2$, we consider the following spaces :
$$
 H^1_{\alpha} (0,\ell):=\{ u \in L^2(0,\ell) \ \mid \ u \text{ locally absolutely continuous in } (0,\ell],
 x^{\alpha/2} u_x \in  L^2(0,\ell)  \} ,
 $$
$$
 H^1_{\alpha,0} (0,\ell):=\{ u \in H^1_{\alpha} (0,\ell) \ \mid \ u(\ell)=0 \} ,$$
and
$$
 H^2_{\alpha} (0,\ell):=\{ u \in H^1_{\alpha} (0,\ell) \ \mid 
 x^\alpha u_x \in H^1(0,\ell) \} 
$$

Then, the operator $A:D(A)\subset L^2(0,\ell)\to L^2(0,\ell)$ will be  defined by
\begin{equation*}
\begin{cases}
\forall u \in D(A), \quad    Au:= (x^\alpha  u_x)_x, & \\
D(A) :=   \{ u \in H^1_{\alpha ,0}(0,\ell)  \ \mid \ x^\alpha  u_x \in  H^1(0,\ell) \}
=  H^2_{\alpha }(0,\ell) \cap  H^1_{\alpha ,0}(0,\ell)
,&\\
 \qquad  \quad =  \{ u \in L^2(0,\ell) \ \mid \ u \text{ locally absolutely continuous in } (0,\ell],&\\
\qquad  \qquad  \qquad
x^\alpha u \in H^1_0(0,\ell), \  x^\alpha  u_x \in  H^1(0,\ell) \text{ and }  (x^\alpha  u_x) (0) =0 \}.&
\end{cases}
\end{equation*}
Notice that, if $u \in D(A)$, then $u$ satisfies the Neumann boundary condition $(x^\alpha  u_x) (0) =0$ at $x=0$ and the  Dirichlet boundary condition $u(\ell)=0$ at $x=\ell$.

Then the following results hold, (see, e.g., \cite{campiti} and \cite{CMV4}).

\begin{Proposition}
\label{prop-A} $A: D(A) \subset L^2(0,\ell)\to L^2(0,\ell)$ is a self-adjoint negative operator
with dense domain.
\end{Proposition}
Hence, $A$ is the infinitesimal generator of an analytic semigroup of contractions $e^{tA}$
on $L^2(0,\ell)$. 
Given a source term $h$ in  $L^2((0,\ell)\times (0,T))$ and an initial condition $v_0 \in  L^2(0,\ell)$, consider the problem
\begin{equation}
\label{eq-v-h}
\begin{cases}
v_t - (x^\alpha v_x)_x = h(x,t), \\
(x^\alpha v_x) (0,t)=0, \\
v(\ell,t) = 0, \\
v(x,0)=v_0(x)  .
\end{cases}
\end{equation}
 The function $v \in \mathcal C ^0 ([0,T]; L^2(0,\ell))  \cap L^2(0,T; H^1_{\alpha,0} (0,\ell))$
given by the variation of constant formula
$$ v(\cdot ,t) = e^{tA} v_0 + \int _0 ^t e^{(t-s)A} h(\cdot, s) \, ds
$$ 
is called the mild solution of \eqref{eq-v-h}. We say that a function 
$$v \in
\mathcal C ^0 ([0,T]; H^1_{\alpha,0} (0,\ell))  \cap   H^1(0,T; L^2(0,\ell)) \cap    L^2 (0,T; D(A)) $$
is a strict solution of \eqref{eq-v-h} if $v$ satisfies $v_t - (x^\alpha v_x)_x = h(x,t)$ almost everywhere
in $(0,\ell )\times (0,T)$, and the initial and boundary conditions  for all $t\in [0,T]$ and all $x\in [0,\ell]$.

\begin{Proposition}
\label{prop-2.2}
If $v_0 \in H^1 _{\alpha,0}(0,\ell)$, then the mild solution of \eqref{eq-v-h} is the unique strict solution of \eqref{eq-v-h}.
\end{Proposition}
The proof of Proposition \ref{prop-2.2} follows from classical results, see subsection \ref{sub-prop-2.2}.


\subsubsection{Functional setting and well-posedness for a boundary control}  \hfill

To define the solution of the boundary value problem \eqref{*pbm-controle1}, we transform it into a problem with homogeneous boundary conditions and a source term (depending on the control $h$): 
formally, if $u$ is a solution of \eqref{*pbm-controle1}, then the function $v$ defined by
\begin{equation}
\label{*v:=u,G}
 v(x,t) = u(x,t) - \frac{x^{2-\alpha}}{\ell^{2-\alpha}} H(t)
\end{equation}
satisfies the auxiliary problem
\begin{equation}
\label{*eq-v-G'}
\begin{cases}
v_t - (x^\alpha v_x)_x = - \frac{x^{2-\alpha}}{\ell^{2-\alpha}} H'(t) +\frac{2-\alpha}{\ell^{2-\alpha}}H(t), \\
(x^\alpha v_x) (0,t)=0, \\
v(\ell,t) = 0, \\
v(x,0)=u_0(x) - \frac{x^{2-\alpha}}{\ell^{2-\alpha}} H(0).
\end{cases}
\end{equation}
Reciprocally, given $H \in H^1 (0,T)$, consider the solution $v$ of 
\begin{equation}
\label{*eq-v-g}
\begin{cases}
v_t - (x^\alpha v_x)_x = - \frac{x^{2-\alpha}}{\ell^{2-\alpha}} H'(t) + \frac{2-\alpha}{\ell^{2-\alpha}}H(t) , \\
(x^\alpha v_x) (0,t)=0, \\
v(\ell,t) = 0, \\
v(x,0)=v_0(x) .
\end{cases}
\end{equation}
Then the function $u$ defined by
\begin{equation}
\label{*u:=v,g} 
u(x,t) = v(x,t) + \frac{x^{2-\alpha}}{\ell^{2-\alpha}} H(t)  
\end{equation}
satisfies
\begin{equation}
\label{*eq-u-g}
\begin{cases}
u_t - (x^\alpha u_x)_x = 0, \\
(x^\alpha u_x) (0,t) = 0 , \\
u(\ell,t) =  H(t) , \\
u(x,0)=v_0(x) + \frac{x^{2-\alpha}}{\ell^{2-\alpha}} H(0) .
\end{cases}
\end{equation}
This motivates the following definition of what is the solution of the boundary value problem \eqref{*pbm-controle1}:

\begin{Definition}
a) We say that $u \in C ([0,T]; L^2(0,\ell)) \cap L^2(0,T; H^1 _\alpha (0,\ell))$ is the mild solution of 
\eqref{*pbm-controle1} if $v$ defined by \eqref{*v:=u,G} is the mild solution of \eqref{*eq-v-G'}.

b) We say that $$u \in C ([0,T]; H^1 _{\alpha} (0,\ell)) \cap H^1(0,T; L^2 (0,\ell))\cap L^2(0,T; H^2 _{\alpha} (0,\ell))$$ is the strict solution of \eqref{*pbm-controle1} if $v$ defined by \eqref{*v:=u,G} is the strict solution of \eqref{*eq-v-G'}.
\end{Definition}
Then we immediately obtain
\begin{Proposition}
\label{prop-u-G}
a) Given $u_0 \in L^2 (0,\ell)$, $H\in H^1(0,T)$, problem \eqref{*pbm-controle1} admits a unique mild solution.

b) Given $u_0 \in H^1 _{\alpha,0} (0,\ell)$, $H \in H^1(0,T)$ such that $u_0(\ell)=H(0)$, problem \eqref{*pbm-controle1} admits a unique strict solution. 
\end{Proposition}
\noindent The proof of Proposition \ref{prop-u-G} follows immediately, noting that
$$ \tilde H (x,t) := \frac{x^{2-\alpha}}{\ell^{2-\alpha}} H(t)  $$ satisfies
$$ \tilde H \in C ([0,T]; H^1 _{\alpha} (0,\ell)) \cap H^1(0,T; L^2 (0,\ell))\cap L^2(0,T; H^2_{\alpha} (0,\ell)) .$$


\subsection{Null controllability results for the boundary control}  \hfill

Consider
\begin{equation}
\label{*def-cost-bd}
C_{ctr-bd}(\alpha,T,\ell) := \sup _{\Vert u_0 \Vert _{L^2(0,\ell)}=1} \inf \{\Vert H \Vert _{H^1 (0,T)}, u^{(H)}(T)=0 \} , 
\end{equation}
where $u^{(H)}$ is the solution of problem \eqref{*pbm-controle1}. Then we prove the following estimates:


\subsubsection{Lower bound of the null controllability cost}\hfill

\begin{Theorem} 
\label{thm-cost-bd-low}
There exists a constant $C_u>0$ independent of $\alpha \in [1,2)$, of $\ell >0$ and of $T>0$ such that

\begin{equation}
\label{*borne-bd-low}
C_{ctr-bd}(\alpha,T,\ell)
\geq 
C_u \frac{\ell ^{2-\alpha}}{\sqrt{(2-\alpha) T \ell} }
e^{- \pi ^2 \frac{T}{\ell ^{2-\alpha}}}
 e^{C_u \frac{\ell ^{2-\alpha}}{T (2-\alpha) ^2}}
e^{-\frac{1}{C_u} (\frac{1}{(2-\alpha)^{4/3}} +  \frac{\ell ^{1-\alpha /2}}{2-\alpha} )(\ln \frac{\ell ^{1-\alpha /2}}{2-\alpha} + \ln \frac{1}{T})} .
\end{equation}
\end{Theorem}

\begin{Remark}
{\rm This proves that the cost blows up when $T\to 0^+$,
or $\alpha \to 2^-$, or $\ell \to +\infty$, and at least exponentially fast. When $\ell$ is fixed and $T\leq T_0$, this simplifies into
$$ C_{ctr-bd}(\alpha,T,\ell)
\geq 
C_u 
 e^{\frac{C_u}{T (2-\alpha) ^2}}
e^{-\frac{1}{C_u} \frac{1}{(2-\alpha)^{4/3}}(\ln \frac{1}{2-\alpha} + \ln \frac{1}{T})} .$$
}
\end{Remark}


\subsubsection{Upper bound of the null controllability cost} \hfill

\begin{Theorem} 
\label{thm-cost-bd-up}
There exists a constant $C_u>0$ independent of $\alpha \in [1,2)$, of $\ell >0$ and of $T>0$ such that

\begin{equation}
\label{*borne-bd-up}
C_{ctr-bd}(\alpha,T,\ell) \leq 
\frac{C_u }{\sqrt{(2-\alpha ) T \ell}} 
e^{- \frac{1}{C_u} \frac{T }{\ell^{2-\alpha}}}
e^{C_u\frac{\ell^{2-\alpha}}{T (2-\alpha) ^2 }} .
\end{equation}
\end{Theorem}

\begin{Remark} {\rm This proves that the cost blows up exactly exponentially fast
as $T\to 0^+$,
or $\alpha \to 2^-$, or $\ell \to +\infty$. When $\ell$ is fixed and $T\leq T_0$, this simplifies into
$$ C_{ctr-bd}(\alpha,T,\ell) \leq 
\frac{C_u }{\sqrt{(2-\alpha ) T}} 
e^{\frac{C_u}{T (2-\alpha) ^2 }} .
$$
}
\end{Remark}


\subsection{Null controllability results for the locally distributed control}  \hfill

Consider
\begin{equation}
\label{*def-cost-loc}
C_{ctr-loc}(\alpha,T,\ell) := 
\sup _{\Vert u_0 \Vert _{L^2(0,\ell)}=1} \inf \{\Vert h \Vert _{L^2((a,b)\times ((0,T))}, u^{(h)}(T)=0 \} ,
\end{equation}
where $u^{(h)}$ is the solution of problem \eqref{*pbm-controle2}. For this problem, we are mainly interested in the dependence with respect to the degeneracy ($\alpha \to 2^-$) and to time (fast controls, when $T \to 0$), see Remarks \ref{rq-low-loc} and \ref{rq-up-loc}.
And we prove the following estimates:


\subsubsection{Lower bound of the null controllability cost}\hfill

\begin{Theorem} 
\label{thm-cost-ld-low}
Given $\ell >0$, and $0<a<b<\ell$, there exists a constant $\tilde C = \tilde C(a,b,\ell)>0$ independent of $\alpha \in [1,2)$ and of $T>0$ such that
\begin{equation}
\label{*borne-loc-low}
C_{ctr-loc}(\alpha,T,\ell) 
\geq 
\tilde C
 e^{\frac{\tilde C}{T (2-\alpha) ^2}}
e^{-\frac{1}{\tilde C} \frac{1}{(2-\alpha)^{4/3}} (\ln \frac{1}{2-\alpha} + \ln \frac{1}{T})
- \frac{1}{\tilde C} T} - 1 .
\end{equation}

\end{Theorem}

\begin{Remark}
\label{rq-low-loc}
{\rm In the proof of Theorem \ref{thm-cost-ld-low} we obtain an explicit expression of $\tilde C(a,b,\ell)$. And of course Theorem \ref{thm-cost-ld-low} proves that the cost blows up (exponentially fast) when $T\to 0^+$,
or $\alpha \to 2^-$: when $T\to 0$
and/or $\alpha \to 2^-$, this simplifies into
$$
C_{ctr-loc}(\alpha,T,\ell) 
\geq \tilde C
 e^{\frac{\tilde C}{T (2-\alpha) ^2}}
e^{-\frac{1}{\tilde C} \frac{1}{(2-\alpha)^{4/3}} (\ln \frac{1}{2-\alpha} + \ln \frac{1}{T})} .
$$ 
}
\end{Remark}


\subsubsection{Upper bound of the null controllability cost} \hfill

\begin{Theorem} 
\label{thm-cost-ld-up}
There exists a constant $C_{u} = C_u>0$ independent of $\alpha \in [1,2)$, of $T>0$ and of $0<a<b<\ell$, and 
$\gamma _0 ^* = \gamma _0 ^* (a,b,\ell) >0$ such that

\begin{equation}
\label{*borne-loc-up}
C_{ctr-loc}(\alpha,T, \ell) \leq
\frac{C_u}{\gamma  ^* _0 } 
e^{C_u \frac{\ell^{2-\alpha}}{T (2-\alpha)^2 }} 
e^{- \frac{1}{C_u} \frac{T}{\ell ^{2-\alpha} }} 
\max \{\frac{1}{\sqrt{T} (2-\alpha)} , \frac{1}{\ell ^{1-\alpha/2}} \} .
\end{equation}
\end{Theorem}

\begin{Remark} 
\label{rq-up-loc}
{\rm This proves that the cost blows up exactly exponentially fast
as $T\to 0^+$,
or $\alpha \to 2^-$, or $\ell \to +\infty$. When $\ell$ is fixed and $T\leq T_0$, this simplifies into
$$ 
C_{ctr-loc}(\alpha,T, \ell) \leq
C_u
e^{C_u \frac{1}{T (2-\alpha)^2 }} .
$$
}
\end{Remark}


\subsection{The eigenvalue problem}  \hfill

The knowledge of the eigenvalues and associated eigenfunctions of the degenerate diffusion operator $ u \mapsto - (x^\alpha u')'$, i.e. the solutions $(\lambda, \Phi )$ of
\begin{equation}
\label{*pbm-vp}
\begin{cases}
 - (x^\alpha \Phi '(x))' =\lambda \Phi (x) & \qquad x\in (0,\ell),\\
(x^\alpha \Phi '(x)) (0)=0 ,\\ 
\Phi (\ell)=0. & 
\end{cases}
\end{equation}
 will be essential for our purposes.
 
\subsubsection{Eigenvalues and eigenfunctions} \hfill
 
It is well-known that Bessel functions play an important role in this problem, see, e.g., Kamke \cite{Kamke}. For $\alpha \in [1,2)$, let
$$ \nu _\alpha := \frac{ \alpha -1 }{2-\alpha}, 
\qquad \kappa _\alpha:= \frac{2-\alpha}{2}.$$
Given $\nu \geq 0$, we denote by $J_\nu$ the Bessel function of first kind and of order $\nu$ (see section \ref{sect-Bessel-rappels}) and denote $j_{\nu,1}< j_{\nu,2} < \dots < j_{\nu,n} <\dots$ the sequence of positive zeros of $J_\nu$. 
Then we have the following: 

\begin{Proposition}
\label{*prop-vp}
The eigenvalues $\lambda$ for problem \eqref{*pbm-vp} are given by 
\begin{equation}
\label{*vp}
\forall n \geq 1, \qquad \lambda_{\alpha, n} = \ell ^{\alpha -2} \kappa _\alpha ^2 j_{\nu _\alpha ,n}^2 
\end{equation}
and the corresponding normalized (in  $L^2(0,\ell)$) eigenfunctions takes the form
\begin{equation}
\label{*fp}
 \Phi_{\alpha, n}(x)=  \frac{\sqrt{2 \kappa _\alpha }}{\ell^{\kappa _\alpha} \vert J'_{\nu_\alpha} (j_{\nu_\alpha,n} ) \vert} 
x^{(1-\alpha)/2} J_{\nu _\alpha} (j_{\nu_\alpha,n} (\frac{x}{\ell}) ^{\kappa_\alpha}), \qquad x \in (0,\ell).
\end{equation}
Moreover the family $(\Phi_{\alpha, n})_{n\geq 1}$ forms an orthonormal basis of $L^2(0,\ell)$.
\end{Proposition}

\begin{Remark}\label{*rq1} 
{\rm 
Gueye \cite{Gueye} proved Proposition \ref{*prop-vp} in the case $\alpha \in [0,1)$ and when $\ell =1$. The case $\alpha \in [1,2)$ and $\ell \neq 1$ is very similar.
}
\end{Remark}


\subsubsection{The eigenfunctions in the control region} \hfill

We will prove the following property:

\begin{Proposition}
\label{prop-fctpr}
Given $0<a<b<\ell$, there exists $\gamma  ^* _0 = \gamma  ^* _0 (a,b,\ell) >0$ such that
\begin{equation}
\label{mmmp}
\forall \alpha \in [1,2), \forall m\geq 1, \quad \int _a ^b \Phi _{\alpha ,m} (x) ^2 \, dx \geq \gamma  ^* _0  (2-\alpha) .
\end{equation}
\end{Proposition}
It is classical in the nondegenerate case (Lagnese \cite{Lagnese}) that 
$$ \inf _m \int _a ^b \Phi _{\alpha ,m} ^2  >0 ;$$
but, in our purpose of estimating the cost of null controllability, it is necessary to have a lower bound of $\int _a ^b \Phi _{\alpha ,m} ^2$
with respect to the degeneracy parameter $\alpha$ when $\alpha \to 2^-$, and the dependence is given in Proposition \ref{prop-fctpr}.
This does not come easily, since $\Phi _{\alpha ,m} $ is solution of a second-order differential equation depending on a large parameter. We will overcome this difficulty with ODE techniques.


\section{Proof of Propositions \ref{*prop-vp} and \ref{prop-2.2}} 
\label{*s3}

In this section, we study the spectral problem \eqref{*pbm-vp} and the properties of the eigenvalues and eigenfunctions, and as a first application we deduce the well-posedness result stated in Proposition \ref{prop-2.2}.

Let us study the spectral problem. First, one can observe that if $\lambda$ is an eigenvalue, then $\lambda >0$: indeed, multiplying \eqref{*pbm-vp} by $\Phi $ and integrating by parts, then
$$ \lambda \int _0 ^\ell \Phi ^2 = \int _0 ^\ell x^\alpha \Phi _x ^2 ,$$
which implies first $\lambda \geq 0$, and next that $\Phi =0$ if $\lambda =0$.


\subsection{The link with the Bessel's equation} \hfill

There is a change a variables that allows one to transform the eigenvalue problem \eqref{*pbm-vp} into a
differential Bessel's equation (see in particular Kamke \cite[section  2.162, equation (Ia), p. 440]{Kamke},
and Gueye \cite{Gueye}): 
assume that $\Phi$ is a solution of \eqref{*pbm-vp} associated to the eigenvalue $\lambda$; then one easily checks that
the function $\Psi$ defined by
\begin{equation}
\label{*eq-lien}
\Phi  (x) =: x^{\frac{1-\alpha}{2}} \Psi \Bigl(\frac{2}{2-\alpha} \sqrt{\lambda} x^{\frac{2-\alpha}{2}} \Bigr)
\end{equation}
is solution of the following boundary problem:
\begin{equation}
\label{*pb-bessel}
\begin{cases}
y^2 \Psi ''(y) + y \Psi  '(y) + (y^2 - (\frac{\alpha -1}{2-\alpha}) ^2) \Psi (y) = 0, \quad y\in (0, \frac{2 }{2-\alpha}\sqrt{\lambda} \ell ^{\kappa _\alpha} ), \\
(2-\alpha) y^{\frac{1}{2-\alpha}} \Psi  '(y) - (\alpha -1) y^{\frac{\alpha -1}{2-\alpha}} \Psi (y) \to 0 \text{ as } y \to 0, \\
\Psi \Bigl(\frac{2}{2-\alpha} \sqrt{\lambda} \ell^{\frac{2-\alpha}{2}} \Bigr) = 0 .
\end{cases}
\end{equation}


\subsection{Bessel's equation and Bessel's functions of order $\nu$} 
\label{sect-Bessel-rappels}

For reader convenience, we recall here the definitions concerning Bessel's equation and functions together with some useful properties of these functions and of their zeros. {\it Throughout this section, we assume that $\nu \in \mathbb R_+$}.

 
\subsubsection{Bessel's equation and Bessel's functions of order $\nu$} \hfill

The Bessel's functions of order $\nu$ are the solutions of the following differential equation (see \cite[section 3.1, eq. (1), p. 38]{Watson} or \cite[eq (5.1.1), p. 98]{Lebedev}): 
\begin{equation}
\label{*eq-bessel-ordre-nu}
y^2 \Psi ''(y) + y \Psi'(y) +(y^2-\nu^2) \Psi(y)=0, \qquad y\in (0,+\infty).
\end{equation}
The above equation is called {\it Bessel's equation for functions of order $\nu$}. 
Of course the fundamental theory of ordinary differential equations says that the solutions of \eqref{*eq-bessel-ordre-nu} generate a vector space $S_\nu$ of dimension 2. In the following we recall what can be chosen as a basis of $S _\nu$.


\subsubsection{Fundamental solutions of Bessel's equation when $\nu \notin \Bbb N$} \hfill

Assume that $\nu \notin \Bbb N$. When looking for solutions of  \eqref{*eq-bessel-ordre-nu} of the form of series of ascending powers of $y$, one can  construct two series that are solutions:
$$ \sum_{m \geq 0}  \frac{(-1)^m}{m! \ \Gamma (\nu+m+1) } \left( \frac{y}{2}\right) ^{\nu+2m}
\ \text{ and } \ 
\sum_{m \geq 0}  \frac{(-1)^m}{m! \ \Gamma (-\nu+m+1) } \left( \frac{y}{2}\right) ^{-\nu+2m},$$
where $\Gamma$ is the Gamma function
(see \cite[section 3.1, p. 40]{Watson}).
 The first of the two series converges for all values of $y$ and defines the so-called Bessel function of order $\nu$ and {\it of the first kind} which is denoted by $J_\nu$:
\begin{equation}
\label{*def-Jnu}
 J_\nu (y)
:= \sum_{m = 0} ^\infty  \frac{(-1)^m}{m! \ \Gamma (m+\nu+1) } \left( \frac{y}{2}\right) ^{2m+\nu}
=\sum_{m = 0} ^\infty  c_{\nu,m} ^+ y^{2m+\nu} , \qquad y \geq 0,
\end{equation}
 (see \cite[section 3.1, (8), p. 40]{Watson} or \cite[eq. (5.3.2), p. 102]{Lebedev}).
The second series converges for all positive values of $y$ and is evidently $J_{-\nu}$:
\begin{equation}
\label{*def-J-nu}
 J_{-\nu} (y)
:= \sum_{m = 0} ^\infty  \frac{(-1)^m}{m! \ \Gamma (m-\nu+1) } \left( \frac{y}{2}\right) ^{2m-\nu}
=\sum_{m = 0} ^\infty  c_{\nu,m} ^- y^{2m-\nu} , \qquad y > 0 .
\end{equation}

When $\nu \not \in  \mathbb N$, the two functions $J_\nu$ and $J_{-\nu}$ are linearly independent and therefore the pair $(J_\nu,J_{-\nu})$ forms a fundamental system of solutions of  \eqref{*eq-bessel-ordre-nu},
(see  \cite[section 3.12, eq. (2), p. 43]{Watson}).


\subsubsection{Fundamental solutions of Bessel's equation when $\nu =n \in \Bbb N$} \hfill
\label{*sub-bessel-entier}

Assume that $\nu =n \in \Bbb N$. When looking for solutions of  \eqref{*eq-bessel-ordre-nu} of the form of series of ascending powers of $y$, one sees that $J_n$ and $J_{-n}$ are still solutions of \eqref{*eq-bessel-ordre-nu}, 
where $J_n$ is still given by \eqref{*def-Jnu} and $J_{-n}$ is given by \eqref{*def-J-nu}; when $\nu =n \in \Bbb N$, $J_{-n}$ can be written
\begin{equation}
\label{*def-J-n}
J_{-n} (y) = \sum_{m \geq n}  \frac{(-1)^m}{m! \ \Gamma (m-n+1) } \left( \frac{y}{2}\right) ^{-n+2m} .
\end{equation}
However now $J_{-n} (y) = (-1)^n J_{n} (y)$, hence $J_n$ and $J_{-n}$ are linearly dependent, (see \cite[section 3.12, p. 43]{Watson} or \cite[eq. (5.4.10), p. 105]{Lebedev}). The determination of a fundamental system of solutions in this case requires further investigation. In this purpose, one introduces the Bessel's functions of order $\nu$ and {\it of the second kind}:
among the several definitions of Bessel's functions of second order, we recall here the definition by Weber. The Bessel's functions of order $\nu$ and {\it of second kind} are denoted by $Y_\nu$ and defined by (see \cite[section 3.54, eq. (1)-(2), p. 64]{Watson} or \cite[eq. (5.4.5)-(5.4.6), p. 104]{Lebedev}):
$$\begin{cases}
\forall \nu \not \in \mathbb N, &\qquad 
Y_\nu(y) := 
\displaystyle{ \frac{J_\nu(y) \cos(\nu \pi)-J_{-\nu} (y)}{\sin(\nu \pi)}},
\\
\forall n \in \mathbb N, &\qquad 
 Y_n(y) := \lim_{\nu\to n} Y_\nu (y).
\end{cases}
$$
For any $\nu \in \mathbb R_+$, the two functions $J_\nu$ and $Y_\nu$ always  are linearly independent, see  \cite[section 3.63, eq. (1), p. 76]{Watson}.
In particular, in the case $\nu=n\in\mathbb N$, the pair $(J_n,Y_n)$ forms a fundamental system of solutions of the Bessel's equation for functions of order $n$.  

In the case $\nu = n \in \mathbb N$,  it will be useful to expand $Y_n$ under the form of a series of ascending powers. This can be done using Hankel's formula, 
see \cite[section 3.52, eq. (3), p. 62]{Watson} 
or \cite[eq. (5.5.3), p. 107]{Lebedev}:
\begin{multline}
\label{*expand-Yn}
\forall n \in \mathbb N^\star, \qquad Y_n(y)=
\frac{2}{\pi} J_n(y) \log \left( \frac{y}{2}\right) 
- \frac{1}{\pi} \sum_{m=0}^{n-1} \frac{(n-m-1)!}{m!} \left(\frac{y}{2}\right) ^{2m-n}
\\
- \frac{1}{\pi} \sum_{m=0}^{+\infty} \frac{(-1)^m}{m!(n+m)!}\left(\frac{y}{2}\right) ^{n+2m} 
\left[ \frac{\Gamma '(m+1)}{\Gamma (m+1)}  + \frac{\Gamma ' (m+n+1)}{\Gamma (m+n+1)}  \right],
\end{multline}
where $\frac{\Gamma '}{\Gamma}$ is the logarithmic derivative of the Gamma function, and satisfies
$\frac{\Gamma '(1)}{\Gamma (1)} = - \gamma$ (here $\gamma$ denotes Euler's constant) and 
$$\frac{\Gamma '(m+1)}{\Gamma(m+1)} = 1+\frac{1}{2} + \ldots \frac{1}{m} -\gamma \ \text{ for all } m \in \mathbb N.$$
In the case $n=0$, the first sum in \eqref{*expand-Yn} should be set equal to zero.


\subsubsection{Zeros of Bessel functions of order $\nu$ of the first kind} \hfill

The function $J_\nu$ has an infinite number of real zeros  which are simple with the possible exception of $x=0$ (\cite[section 15.21, p. 478-479 applied to $C_\nu=J_\nu$]{Watson} or \cite[section 5.13, Theorem 2, p. 127]{Lebedev}). We denote by $(j_{\nu,n})_{n\geq 1}$ the strictly increasing sequence of the positive zeros of $J_{\nu}$:
$$ 0< j_{\nu,1} < j_{\nu,2} < \dots < \ j_{\nu,n} < \dots$$
and we recall that 
$$  j_{\nu,n} \to +\infty \text{ as } n \to +\infty.$$
We will also often use the following bounds on the zeros, proved in Lorch and Muldoon \cite{Lorch}:
\begin{equation}
\label{*eq-Lorch}
\forall \nu \geq \frac{1}{2}, \forall n\geq 1, \quad 
\pi (n + \frac{\nu}{4}-\frac{1}{8}) \leq j_{\nu, n} \leq\pi (n + \frac{\nu}{2}-\frac{1}{4}) .
\end{equation}
Note also that (\cite{Lorch}):
\begin{equation}
\label{*eq-Lorch2}
\forall \nu \in [0, \frac{1}{2}], \forall n\geq 1, \quad 
\pi (n + \frac{\nu}{2}-\frac{1}{4}) \leq j_{\nu, n} \leq \pi (n + \frac{\nu}{4}-\frac{1}{8}) .
\end{equation}

We will also use the following asymptotic development of the first zero $j_{\nu,1}$ 
of $J_\nu$ with respect to $\nu$ (\cite[section 15.81, p. 516]{Watson}) when $\nu \to + \infty$:
$$ j_{\nu,1} = \nu + 1,855757 \nu ^{1/3} + O(1) ,$$
and a similar develoment for $j_{\nu,2}$, extracted from \cite{QuWong} where it is proved that
\begin{equation}
\label{eq-Qu-Wong} 
\nu - \frac{a_k}{2^{1/3}} \nu ^{1/3} < j_{\nu,k}  < \nu - \frac{a_k}{2^{1/3}} \nu ^{1/3}
+ \frac{3}{20} a_k ^2 \frac{2^{1/3}}{\nu ^{1/3}} ,
\end{equation}
which is valid for all $\nu >0$, all $k\geq 1$, and where $a_k$ is the $k$-th negative zero of the Airy function.


\subsection{Proof of Proposition \ref{*prop-vp}} \hfill
\label{*s4}

As noted before, $\lambda=0$ is not an admissible eigenvalue for problem
\eqref{*pbm-vp}, hence $\lambda >0$. So, using \eqref{*eq-lien}, we can transform problem \eqref{*pbm-vp} into problem \eqref{*pb-bessel}. In the following, because of the difference in the construction of a fundamental system of solutions of \eqref{*eq-bessel-ordre-nu}, we treat the following cases separately: $\nu _\alpha \not \in \mathbb N$, $\nu _\alpha =n \in\mathbb N^*$
and $\nu_\alpha = 0$.


\subsubsection{Case $\nu _\alpha \not \in \mathbb N$.} 

Let us assume that $\nu _\alpha \not \in \mathbb N$. Then we have
$$ \Phi = C_+ \Phi _+ + C_- \Phi _-$$
where $\Phi _+$ and $\Phi _-$ are defined by
\begin{equation*}
\label{*base-phi}
\Phi _+ (x) := x^{\frac{1-\alpha}{2}} J_{\nu _\alpha} (\frac{2}{2-\alpha} \sqrt{\lambda} x^{\frac{2-\alpha}{2}}),
\quad
\Phi _- (x) := x^{\frac{1-\alpha}{2}} J_{-\nu _\alpha} (\frac{2}{2-\alpha} \sqrt{\lambda} x^{\frac{2-\alpha}{2}}).
\end{equation*}
Using the series expansion of $J_{\nu _\alpha}$ and $J_{-\nu _\alpha}$, one obtains
\begin{equation}
\label{*serie-phi-cas2}
\Phi _+ (x) = \sum _{m=0} ^\infty \tilde{c} _{\nu _\alpha ,m} ^+  x ^{(2-\alpha) m},
\quad 
\Phi _- (x) = \sum _{m=0} ^\infty \tilde{c} _{\nu _\alpha ,m} ^- x ^{1-\alpha + (2-\alpha) m} ,
\end{equation}
where the coefficients $\tilde{c} _{\nu _\alpha ,m} ^+$ and $\tilde{c} _{\nu _\alpha ,m} ^-$ are defined by
\begin{equation}
\label{*coeffs-phi}
 \tilde{c} _{\nu _\alpha ,m} ^+ := c_{\nu _\alpha ,m} ^+ \Bigl( \frac{2}{2-\alpha} \sqrt{\lambda} \Bigr) ^{2m+\nu _\alpha},
\quad 
\tilde{c} _{\nu _\alpha ,m} ^- := c_{\nu _\alpha ,m} ^- \Bigl( \frac{2}{2-\alpha} \sqrt{\lambda} \Bigr) ^{2m-\nu _\alpha} .
\end{equation}
We deduce that
$$ \Phi _+ (x) \sim _{0} \tilde{c} _{\nu _\alpha ,0} ^+,
\quad x^{\alpha /2} \Phi _+ ' (x) \sim _{0}  (2-\alpha) \tilde{c} _{\nu _\alpha ,1} ^+  x ^{1-\alpha /2} ,$$
$$ \Phi _- (x) \sim _{0} \tilde{c} _{\nu _\alpha ,0} ^- x ^{1-\alpha},
\quad x^{\alpha /2} \Phi _- ' (x) \sim _{0}  (1-\alpha) \tilde{c} _{\nu _\alpha ,0} ^-  x ^{-\alpha /2} ,$$
hence $\Phi _+\in H^1 _\alpha (0,\ell)$, while $\Phi _- \notin H^1 _\alpha (0,\ell)$.
Therefore, $\Phi = C_+ \Phi _+ + C_- \Phi _- \in H^1 _\alpha (0,\ell)$
implies that $C_- =0$ and $\Phi = C_+ \Phi _+$. Moreover, $x^\alpha \Phi _+' (x) \to 0$ as $x\to 0$, hence
the boundary condition in $0$ is automatically satisfied. Finally, the boundary condition $\Phi (\ell)=0$ implies that there is some $C_+$ and some $m \in \Bbb N$, $m\geq 1$ such that
$$ \lambda = \kappa _\alpha ^2  j_{\nu _\alpha, m}^2 \ell ^{\alpha -2}
\quad \text{ and } \quad \Phi (x) = C_+ x^{\frac{1-\alpha}{2}} J_{\nu _\alpha} ( j_{\nu _\alpha, m} (\frac{x}{\ell})^{\kappa_\alpha}) .$$
In the same way, any $\Phi (x) := Cx^{\frac{1-\alpha}{2}} J_{\nu _\alpha} ( j_{\nu _\alpha, m} (\frac{x}{\ell})^{\kappa_\alpha}) $ is solution of \eqref{*pbm-vp}, 
and the family $(\Phi _n (x) := x^{\frac{1-\alpha}{2}} J_{\nu _\alpha} ( j_{\nu _\alpha, n} (\frac{x}{\ell})^{\kappa _\alpha }) )_n$ forms an orthogonal family of $L^2 (0,\ell)$, which is complete since the family is composed by the eigenfunctions of the operator $T_\alpha$:
$$ T_\alpha: L^2(0,\ell) \to L^2 (0,\ell), \quad f \mapsto T_\alpha (f) :=u_f $$
where $u_f \in D(A) $ is the solution of the problem $-Au_f = f $, which is self-adjoint and compact (Appendix in \cite{fatiha}).
Finally, it remains to norm this orthogonal family:
\begin{multline*}
\int _0 ^\ell x^{1-\alpha} J_{\nu _\alpha} ^2 ( j_{\nu _\alpha, n} (\frac{x}{\ell})^{\kappa _\alpha }) \, dx
= \ell ^{2-\alpha} \int _0 ^1 y^{1-\alpha} J_{\nu _\alpha} ^2 ( j_{\nu _\alpha, n} y^{\kappa _\alpha }) \, dy
\\
= \frac{\ell ^{2-\alpha} }{\kappa _\alpha} \int _ 0 ^1 z  J_{\nu _\alpha} ^2 ( j_{\nu _\alpha, n} z) \, dz
= \frac{\ell ^{2-\alpha} }{\kappa _\alpha}  \frac{1}{2}
 [J_{\nu _\alpha +1} (j_{\nu _\alpha,n} ) ]^2
= \frac{\ell ^{2-\alpha} [J' _{\nu _\alpha} (j_{\nu _\alpha,n} ) ]^2 }{2\kappa _\alpha} ,
\end{multline*}
which gives us that the family given by \eqref{*fp} forms an orthonormal basis of $L^2(0,\ell)$. 
This ends the proof of Proposition \ref{*prop-vp} when $\alpha \in [1,2)$ is such that $\nu _\alpha \notin \Bbb N$. 


\subsubsection{Case $\nu  _\alpha =n _\alpha  \in \mathbb N ^*$.}

 Let us assume that $\nu _\alpha =n _\alpha \in \mathbb N ^*$. In this case, we have recalled in subsection \ref{*sub-bessel-entier} that a fundamental system of the differential equation \eqref{*eq-bessel-ordre-nu} is given by $J_{ n_\alpha}$ and $Y_{n _\alpha}$.
This gives us that $\Phi$ is a linear combination of $\Phi _+$ and $\Phi _{+,-}$, where
\begin{equation}
\label{*def-phi+-}
\Phi _{+,-} (x) := x^{\frac{1-\alpha}{2}} Y_{n_\alpha}  (\frac{2}{2-\alpha} \sqrt{\lambda} x^{\frac{2-\alpha}{2}}).
\end{equation}
As we have done above, we now study if
$\Phi _{+,-} \in H^1 _\alpha (0,\ell)$. First we need its decomposition in series: it follows from \eqref{*expand-Yn} that

\begin{equation}
\label{*serie-+-}
\Phi _{+,-} (x)
= \frac{2}{\pi} \Phi _{+} (x) \log \left( \frac{1}{2-\alpha} \sqrt{\lambda} x^{\frac{2-\alpha}{2}}\right)  
+  \sum_{m=0}^{n_\alpha-1} \hat a_m x ^{(1-\alpha) +(2-\alpha)m}
+  \sum_{m=0}^{+\infty} \hat b_m x ^{(2-\alpha)m},
\end{equation}
where
$$\hat a_m:= 
- \frac{1}{\pi}\frac{(n_\alpha -m-1)!}{m!}\left(\frac{\sqrt{\lambda}}{2\kappa _\alpha }\right) ^{2m-n_\alpha}$$
and 
$$\hat b_m:=- \frac{1}{\pi} \frac{(-1)^m}{m!(n_\alpha+m)!}\left(\frac{\sqrt{\lambda}}{2\kappa_\alpha}\right) ^{2m+ n_\alpha} \left[ \frac{\Gamma ' (m+1)}{\Gamma (m+1)} + \frac{\Gamma ' (m+n_\alpha+1)}{\Gamma (m+n_\alpha+1)}  \right].$$
We study the three functions that appear in the formula of $\Phi _{+,-}$.
First $$ \Phi _{+,-,1} (x) := \frac{2}{\pi} \Phi _{+} (x) \log \left( \frac{1}{2-\alpha} \sqrt{\lambda} x^{\frac{2-\alpha}{2}}\right)$$
satisfies
$$ \Phi _{+,-,1} (x) \sim _0 \frac{2\kappa_\alpha}{\pi} \tilde{c} _{n _\alpha ,0} ^+ \log x, \quad x^{\alpha /2} \Phi _{+,-,1}' (x) \sim _0 \frac{2\kappa_\alpha}{\pi} \tilde{c} _{n _\alpha ,0} ^+ x^{-1+\alpha /2},$$
hence $\Phi _{+,-,1} \in H^1 _\alpha (0,\ell)$ since $\alpha >1$.
Next $$ \Phi _{+,-,2} (x) := \sum_{m=0}^{n_\alpha-1} \hat a_m x ^{(1-\alpha) +(2-\alpha)m} $$
satisfies
$$ \Phi _{+,-,2} (x) \sim _0  \hat{a}_0 x^{1-\alpha}, \quad x^{\alpha /2} \Phi _{+,-,2}' (x) \sim _0 (1-\alpha)  \hat{a}_0 x^{-\alpha /2},$$
hence $\Phi _{+,-,2} \notin H^1 _\alpha (0,\ell)$, since $ \hat{a}_0 \neq 0$.
Finally, $$ \Phi _{+,-,3} (x) := \sum_{m=0}^{+\infty} \hat b_m x ^{(2-\alpha)m} $$
satisfies
$$ \Phi _{+,-,3} (x) \sim _0 \hat b _0, \quad x^{\alpha /2} \Phi _{+,-,3}' (x) \sim _0 (2-\alpha) \hat b_1 x^{1-\alpha /2},$$
hence $\Phi _{+,-,3} \in H^1 _\alpha (0,\ell)$.
Thus $\Phi _{+,-} = \Phi _{+,-,1} + \Phi _{+,-,2} + \Phi _{+,-,3} \notin H^1 _\alpha (0,\ell)$, and if
$\Phi = C_+ \Phi _+ + C_{+,-} \Phi _{+,-} \in H^1 _\alpha (0,\ell)$ then necessarily $C_{+,-}=0$, and 
$\Phi = C_+ \Phi _+$. Then we are in the same position as in the previous case and the conclusion is the same.


\subsubsection{Case $\nu  _\alpha =0$ (hence $\alpha =1$).}  

In this case, the first sum in the decomposition of $Y_0$ is equal to zero, hence we have $ \Phi _{+,-} = \Phi _{+,-,1} + \Phi _{+,-,3}$. Moreover,
$$ \Phi _{+,-,1} (x) \sim _0 \frac{2\kappa_1}{\pi} \tilde{c} _{0,0} ^+ \ln x, \quad x^{\alpha /2} \Phi _{+,-,1}' (x) \sim _0 \frac{2\kappa_1}{\pi} \tilde{c} _{0,0} ^+ x^{-1/2},$$
hence $\Phi _{+,-,1} \notin H^1 _\alpha (0,\ell)$. On the contrary $\Phi _{+,-,3} \in H^1 _\alpha (0,\ell)$,
which implies that, once again, $\Phi _{+,-} = \Phi _{+,-,1} + \Phi _{+,-,3} \notin H^1 _\alpha (0,\ell)$, and the conclusion is the same. \qed 


\subsection{Proof of Proposition \ref{prop-2.2}} \hfill
\label{sub-prop-2.2}

Since $\{\Phi _{\alpha,n} , n\geq 1 \}$ is an orthonormal basis of $L^2 (0, \ell)$,
it suffices to observe that 
\begin{equation}
\label{*P}
H^1 _{\alpha,0} (0, \ell)
= \{u \in L^2 (0,\ell), \sum _{n=1} ^\infty \lambda _{\alpha,n} (u,\Phi _{\alpha,n})_{L^2(0,\ell)}  ^2 < \infty \}
(= D((-A)^{1/2})).
\end{equation}
Indeed, since $A$ generates an analytic semigroup of negative type on $X=L^2(0, \ell)$, the conclusion follows from the variation of constant formula
$$ u(t) = e^{tA}u_0 + \int _0 ^t e^{(t-s)A} h(\cdot, s) \, ds $$
and well-known maximal regularity results which ensure  that both maps
$$ t \mapsto e^{tA}u_0 \quad \text{ and } \quad t \mapsto \int _0 ^t e^{(t-s)A} h(\cdot, s) \, ds $$
belong to $H^1 (0,T; X) \cap L^2 (0,T; D(A)) \cap \mathcal C ^0 ([0,T]; D((-A)^{1/2}))$
whenever $u_0 \in D((-A)^{1/2})$ and $h \in L^2(0,T;X)$ (see, e.g., \cite{bensoussan1992}).
Finally, in order to check \eqref{*P}, it suffices to observe that for any $u\in D(A)$, given by 
$$ u= \sum _{n=1} ^\infty  (u,\Phi _{\alpha,n})_{L^2(0,\ell)} \Phi _{\alpha,n} ,$$
we have that 
$$ \int _0 ^\ell a(x) u_x ^2 \, dx = - (Au,u) _{L^2(0,\ell)}
= \sum _{n=1} ^\infty \lambda _{\alpha,n} (u,\Phi _{\alpha,n})_{L^2(0,\ell)}  ^2 .\qed $$


\section{Preliminaries: the moment method}
\label{sec-moment}

We follow the strategy initiated by Fattorini and Russell \cite{FR1, FR2}.
The precise estimates given in Theorems \ref{thm-cost-bd-low}-\ref{thm-cost-ld-up} are based on 
identities given by the moment method. We separate the boundary case from the locally distributed case.


\subsection{The boundary control problem \eqref{*pbm-controle1}}

\subsubsection{The moment problem satisfied by a control $H \in L^2(0,T)$} \hfill

In this part, we analyze the problem with formal computations. First, we expand the initial condition $u_0 \in L^2(0,\ell)$:
there exists $(\mu_{\alpha ,n}^0)_{n\geq 1} \in \ell ^2(\mathbb N^\star )$ such that
$$ u_0 (x) =\sum _{n\geq 1} \mu_{\alpha ,n} ^0 \Phi_{\alpha ,n} (x).
$$
Next we expand the solution $u$ of \eqref{*pbm-controle1}: 
$$ u(x,t)= \sum_{n\geq 1} \beta_{\alpha ,n} (t) \Phi_{\alpha ,n} (x), \qquad x\in (0,\ell), \ t \geq 0$$
with
$$\sum_{n\geq 1 } \beta_{\alpha ,n} (t)^2 <+\infty .$$
Therefore the controllability condition $u(\cdot,T) =0$ becomes
$$\forall n \geq 1, \qquad \beta_{\alpha ,n} (T)= 0.$$

On the other hand, we observe that 
$w_{\alpha ,n} (x,t):= \Phi_{\alpha ,n} (x) e^{\lambda_{\alpha ,n} (t-T)}$ is solution of the adjoint problem:
\begin{equation}
\label{*pbm-controle3-adjoint}
\begin{cases}
(w_{\alpha ,n})_t +(x^\alpha (w_{\alpha ,n})_x)_x =0 & \qquad x\in(0,\ell),\ t>0,\\
(x^\alpha (w_{\alpha ,n})_x) (0,t)=0, \\
w_{\alpha ,n} (\ell,t)=0 & \qquad t>0.
\end{cases}
\end{equation}
Multiplying \eqref{*pbm-controle1} by $w_{\alpha ,n}$ and \eqref{*pbm-controle3-adjoint} by $u$, we obtain 
\begin{multline*}
0 = \int_0^T\int_0^\ell w_{\alpha ,n} (u_t -(x^\alpha u_x)_x) + u ((w_{\alpha ,n})_t+(x^\alpha (w_{\alpha ,n})_x)_x) \\
= \int_0^\ell [ w_{\alpha ,n} u ]_0^T dx - \int_0^T [w_{\alpha ,n} x^\alpha u_x ]_0^\ell dt  + \int_0^T [u x^\alpha (w_{\alpha ,n})_x ]_0^\ell dt\\
= \int_0^\ell u(x,T) \Phi_{\alpha ,n}(x)   dx  - \int_0^\ell u(x,0) \Phi_{\alpha ,n}(x)  e^{-\lambda_{\alpha ,n} T } dx
+ \int_0^T u(\ell,t)  (x^\alpha (w_{\alpha ,n})_x)  (\ell,t)  dt\\
= \beta_{\alpha ,n}(T) -  e^{-\lambda_{\alpha ,n} T }  \mu_{\alpha ,n}^0 + \int_0^T  H(t) e^{\lambda_{\alpha ,n} (t-T) } (x^\alpha \Phi_{\alpha ,n} ')(x=\ell) dt .
\end{multline*}
It follows that, if the control $H$ drives the solution to $0$ at time $T$, then
$$   r_{\alpha ,n} \int_0^T  H(t) e^{- \lambda_{\alpha ,n} (T-t) }  dt = e^{-\lambda_{\alpha ,n} T }  \mu_{\alpha ,n} ^0 ,$$
where we have set
\begin{equation}
\label{*infl-fp}
r_{\alpha ,n} = (x^\alpha \Phi_{\alpha ,n} ')(x=\ell).
\end{equation}

Hence, the controllability condition $u(\cdot,T) =0$ implies that
\begin{equation}
\label{*moment-bd}
\forall n \geq 1, \quad  r_{\alpha ,n} \int_0^T   H(t) e^{ \lambda_{\alpha ,n} t }  dt 
= \mu_{\alpha ,n} ^0. 
\end{equation}


\subsubsection{The moment problem satisfied by a control $H \in H^1 (0,T)$} \hfill

Moreover, since we want a solution of the moment problem that belongs to $H^1(0,T)$, it will be more interesting to see what its derivative has to satisfy. Integrating by parts, we have
$$\int _0 ^T H(t) e^{\lambda _{\alpha,n}t} \, dt
= [ \frac{1}{\lambda _{\alpha,n}} H(t)  e^{\lambda _{\alpha,n}t} ] _0 ^T
- \int _0 ^T \frac{1}{\lambda _{\alpha,n}} H'(t)  e^{\lambda _{\alpha,n}t} \, dt .$$
Hence the derivative $H'$ has to satisfy
\begin{equation}
\label{*moment-bd-G'}
 -  \frac{r_{\alpha ,n}}{\lambda _{\alpha,n}} \int _0 ^T H'(t)  e^{\lambda _{\alpha,n}t} \, dt
= \mu_{\alpha ,n} ^0 
- \frac{r_{\alpha ,n}}{\lambda _{\alpha,n}}  \Bigl[ H(T)  e^{\lambda _{\alpha,n}T} - H(0) \Bigr] .
\end{equation}
We will provide a solution of this problem that satisfies $H(0)=0=H(T)$.


\subsubsection{A formal solution to the moment problem, using a biorthogonal family} \hfill
\label{*sect-formal}

Assume that there is a family $(\sigma_{\alpha, m}^+)_{m \geq 1}$ of functions $\sigma_{\alpha, m}^+ \in L^2(0,T)$, which is biorthogonal to the family $(e^{\lambda_{\alpha ,n} t })_{n\geq 1}$, which means that:

\begin{equation}
\label{biortho+1}
 \forall m, n \geq 1, \quad \int _0 ^T \sigma_{\alpha, m}^+ (t) e^{\lambda_{\alpha ,n} t } \, dt = \delta _{mn}=
\begin{cases} 1 \text{ if } m=n , \\
0 \text{ if } m \neq n .
\end{cases}
\end{equation}
Then, at least formally, the function
$$ H(t) := \sum _{m=1} ^\infty \frac{\mu_{\alpha ,m} ^0}{r_{\alpha ,m}} \sigma_{\alpha, m}^+ (t) $$
satisfies the moment problem \eqref{*moment-bd}. To enter into our functional setting, we would need to verify that this gives a function belonging to $H^1 (0,T)$, then at least to $L^2 (0,T)$. For this, we will need suitable bounds on $\Vert \sigma_{\alpha, m}^+ \Vert _{L^2(0,T)}$, first with respect to $m$ (to ensure the convergence of the series that defines $H$), then with respect to $\alpha$, to measure the null controllability cost.

Since our functional setting demands the control to belong to $H^1(0,T)$, we are going to repeat the same arguments, but with the moment problem \eqref{*moment-bd-G'}: set
$$ \lambda_{\alpha ,0} := 0 ,$$
and assume that we are able to construct a family 
$(\sigma_{\alpha, m} ^+)_{m \geq 0}$ 
of functions 
$\sigma_{\alpha, m} ^+ \in L^2(0,T)$, 
which is biorthogonal to the family 
$(e^{\lambda_{\alpha ,n} t })_{n\geq 0}$, 
which means that:

\begin{equation}
\label{biortho+0}
 \forall m, n \geq 0, \quad \int _0 ^T \sigma_{\alpha, m}^+ (t) e^{\lambda_{\alpha ,n} t } \, dt = \delta _{mn}=
\begin{cases} 1 \text{ if } m=n , \\
0 \text{ if } m \neq n .
\end{cases}
\end{equation}
Then consider 
$$ K(t) := - \sum _{m=1} ^\infty \frac{\lambda _{\alpha,m} \mu_{\alpha ,m} ^0}{r_{\alpha ,m}} \sigma_{\alpha, m}^+ (t), \quad \text{ and } \quad H(t) := \int _0 ^t K(\tau) \, d\tau .$$
Then at least formally $K$ solves the following moment problem
\begin{equation*}
\forall n\geq 1, \quad  -  \frac{r_{\alpha ,n}}{\lambda _{\alpha,n}} \int _0 ^T K(t)  e^{\lambda _{\alpha,n}t} \, dt
= \mu_{\alpha ,n} ^0  ;
\end{equation*}
moreover, if $K\in L^2(0,T)$, then $H\in H^1(0,T)$, clearly $H'=K$, and $H(0)=0$, and moreover $H(T)=0$ thanks 
to the additional property that the family $(\sigma_{\alpha, m}^+)_{m \geq 1}$ is orthogonal to
$e^{\lambda_{\alpha ,0} t } =1$; hence $H$ will be in $H^1(0,T)$ and will satisfy the moment problem \eqref{*moment-bd-G'}.
 It remains to check that all this makes sense, in particular that $K\in L^2(0,T)$. Clearly, we will need suitable $L^2$ bounds on the biorthogonal sequence $(\sigma_{\alpha,m}^+)_{m\geq 1}$, that will come from the study of the eigenvalues $\lambda _{\alpha,n}$, and from the behavior of the real sequence $(r_{\alpha ,m}^2)_m$.
 

\subsection{The locally distributed control problem \eqref{*pbm-controle2}}
\label{sec-mmloc}

\subsubsection{The moment problem satisfied by a control $h \in L^2((a,b)\times (0,T))$} \hfill

First we expand the initial condition $u_0 \in L^2(0,\ell)$: 
there exists $(\mu_{\alpha ,n} ^0)_{n\geq 1} \in \ell ^2(\mathbb N^\star)$
such that
$$ u_0 (x) =\sum _{n\geq 1} \mu_{\alpha ,n} ^0 \Phi_{\alpha ,n} (x), 
\qquad x \in (0,\ell).$$
Next we expand the solution $u$ of \eqref{*pbm-controle2}: 
$$u(x,t)= \sum_{n\geq 1} \beta_{\alpha ,n} (t) \Phi_{\alpha ,n} (x), \qquad x\in (0,\ell), \ t \in (0,T),
\quad \text{ with } \sum_{n\geq 1} \beta_{\alpha ,n} (t)^2 <+\infty.$$
Once again multiplying \eqref{*pbm-controle2} by $w_{\alpha ,n} (x,t):= \Phi_{\alpha ,n} (x) e^{\lambda_{\alpha ,n} (t-T)}$, which is solution of the adjoint problem \eqref{*pbm-controle3-adjoint}, one gets:
\begin{multline*}
\int _0 ^T \int _0 ^\ell h(x,t) \chi_{[a,b]}(x) w_{\alpha ,n} (x,t) \, dx\, dt
= \int _0 ^T \int _0 ^\ell  w_{\alpha ,n} (x,t) (u_t - (x^\alpha u_x)_x)
\\
= \int _0 ^\ell  [w_{\alpha ,n} u ]_0 ^T - \int _0 ^T \int _0 ^\ell (w_{\alpha ,n})_t u
- \int _0 ^T [  w_{\alpha ,n} (x^\alpha u_x)]_0 ^\ell + \int _0 ^T \int _0 ^\ell  (w_{\alpha ,n})_x x^\alpha u_x
\\
= \int _0 ^\ell  \Phi_{\alpha ,n} u(T) - e^{-\lambda_{\alpha ,n} T}\int _0 ^\ell  \Phi_{\alpha ,n} u_0
- \int _0 ^T [  w_{\alpha ,n} (x^\alpha u_x)]_0 ^\ell  
\\
+ \int _0 ^T [ x^\alpha (w_{\alpha ,n})_x u ]_0 ^\ell - \int _0 ^T \int _0 ^\ell \Bigl( (w_{\alpha ,n})_t + ( x^\alpha (w_{\alpha ,n})_x )_x \Bigr) u
\\
=\int _0 ^\ell  \Phi_{\alpha ,n} u(T) - e^{-\lambda_{\alpha ,n} T}\int _0 ^\ell  \Phi_{\alpha ,n} u_0 .
\end{multline*}
Hence, if $h$ drives the solution $u$ to $0$ in time $T$, we obtain the following moment problem:
\begin{equation}
\label{*moment2}
\forall n \geq 1, \qquad 
\int_0^T \int_0^\ell h(x,t) \chi_{[a,b]}(x) \Phi_{\alpha ,n} (x) e^{\lambda_{\alpha ,n} t} dx  dt
=-\mu_{\alpha ,n} ^0  .
\end{equation}


\subsubsection{A formal solution to the moment problem, using a biorthogonal family} \hfill
\label{*s5}

Assume for a moment that there exists a family $(\sigma_{\alpha, m}^+)_{m\geq 1}$ in $L^2(0,T)$ that satisfies \eqref{biortho+1}.
Then let us define
\begin{equation}
\label{*contr} 
h(x,t):= \sum_{m\geq 1} -\mu_{\alpha ,m} ^0 \sigma_{\alpha ,m}^+ (t) \frac{\Phi_{\alpha ,m}(x)}{\int_a^b \Phi_{\alpha ,m}^2 } .
\end{equation}
Let us  prove that, formally, $h$ is solution of the moment problem \eqref{*moment2}:
\begin{multline*}
\int_0^T \int_0^\ell h(x,t) \chi_{[a,b]}(x) \Phi_{\alpha ,n}(x) e^{\lambda_{\alpha ,n} t} dxdt
\\
=
\int_a^b \int_0^T \left(  \sum_{m\geq 1} -\mu_{\alpha ,m} ^0  \sigma_{\alpha ,m}^+(t) \frac{\Phi_{\alpha ,m}(x)}{\int_a^b \Phi_{\alpha ,m}^2 }  \right) \Phi_{\alpha ,n}(x) e^{\lambda_{\alpha ,n} t } dtdx
\\
=\int_a^b \sum_{m\geq 1}-\mu_{\alpha ,m} ^0  \frac{\Phi_{\alpha ,m}(x)\Phi_{\alpha, n}(x) }{\int_a^b \Phi_{\alpha ,m}^2 }
\left( \int_0^T \sigma_{\alpha ,m}^+(t) e^{\lambda_{\alpha ,n} t } dt \right) dx
\\
= \sum_{m\geq 1} -\mu_{\alpha ,m} ^0  \delta_{mn} \frac{\int_a^b \Phi_{\alpha ,m}(x)\Phi_{\alpha ,n}(x)dx  }{\int_a^b \Phi_{\alpha ,m}^2 }  
= -\mu_{\alpha ,n} ^0 .
\end{multline*}
Hence, formally, $h$ defined by \eqref{*contr}  solves the moment problem. It remains to check that all this makes sense, in particular that $h \in L^2((0,\ell)\times (0,T))$. Clearly, we will need suitable $L^2$ bounds on the biorthogonal sequence $(\sigma_{\alpha,m}^+)_{m\geq 1}$, that will come from the study of the eigenvalues $\lambda _{\alpha,n}$, and from the behavior of the real sequence $(\int_a^b \Phi_{\alpha ,m}^2)_m$
(given in Proposition \ref{prop-fctpr}).


\section{Proof of Theorem \ref{thm-cost-bd-low}}
\label{sec-preuveThm1}

In this section, we are going to work on the moment problem \eqref{*moment-bd} given by the moment method to obtain the desired lower bound on the null controllability cost for \eqref{*pbm-controle1}. The proof will use in particular ideas of G\"uichal \cite{Guichal}. 


\subsection{The contribution of the eigenfunctions to the blow-up of the null controllability cost} \hfill

Assume that $H \in H^1(0,T)$ drives the solution $u$ of \eqref{*pbm-controle1} to $0$ in time $T$. Then $H$ satisfies \eqref{*moment-bd}. Let us compute the coefficient that appears:
\begin{multline}
\label{*infl-fp2}
\vert r_{\alpha ,n} \vert 
= \vert  (x^\alpha \Phi_{\alpha ,n} ')(x=\ell) \vert 
= \ell ^\alpha \frac{\sqrt{2 \kappa _\alpha }}{\ell^{\kappa _\alpha} \vert J'_{\nu_\alpha} (j_{\nu_\alpha,n} ) \vert} \ell^{(1-\alpha)/2} j_{\nu_\alpha,n} \frac{\kappa_\alpha}{\ell}  \vert J'_{\nu_\alpha} (j_{\nu_\alpha,n} ) \vert 
\\
= 
\ell^{\alpha - \kappa_\alpha + (1-\alpha)/2 -1} \sqrt{2} {\kappa _\alpha} ^{3/2} j_{\nu_\alpha,n} 
= 
\sqrt{2} \, \ell ^{(2\alpha-3)/2} \, \kappa _\alpha ^{3/2} \,  j_{\nu_\alpha,n}
\\
= 
\sqrt{2} \, \ell ^{(\alpha-1)/2} \, \sqrt{\kappa _\alpha} \,  \sqrt{\lambda_{\alpha,n}}  .
\end{multline}
 This implies that the null controllability cost blows up, at least at a rational rate: indeed, we deduce from \eqref{*moment-bd} that
\begin{equation*}
\forall n \geq 1, \quad  \Vert H \Vert _{L^2(0,T)} \, \Vert e^{ \lambda_{\alpha ,n} t }  \Vert _{L^2(0,T)} \geq \frac{\vert \mu_{\alpha ,n} ^0 \vert}{\vert r_{\alpha ,n} \vert}, 
\end{equation*}
hence 
$$ \forall n \geq 1, \quad  \Vert H \Vert _{L^2(0,T)} 
\geq \frac{\vert \mu_{\alpha ,n} ^0 \vert}{\sqrt{2} \, \ell ^{(\alpha-1)/2} \, \sqrt{\kappa _\alpha} \,  \sqrt{\lambda_{\alpha,n}}} \sqrt{\frac{2\lambda_{\alpha,n}}{e^{2\lambda_{\alpha,n}T}-1}}.$$
Fix $u_0=\Phi _{\alpha,1}$. Then any control that drives $\Phi _{\alpha,1}$ to $0$ in time $T$ satisfies
$$ \Vert H \Vert _{L^2(0,T)} 
\geq \frac{1}{\ell ^{(\alpha-1)/2} \, \sqrt{\kappa _\alpha} \, \sqrt{e^{2\lambda_{\alpha,1}T}-1}}  .
$$
This implies a first bound from below for the null controllability cost:
$$ C_{ctr-bd} \geq \frac{1}{\ell ^{(\alpha-1)/2} \, \sqrt{\kappa _\alpha} \, \sqrt{e^{2\lambda_{\alpha,1}T}-1}}  .$$
In particular, just looking the behavior with respect to $\alpha \in [1,2)$, we see
 that there exists $C_{T,\ell}$ independent of $\alpha \in [1,2)$ such that
$$ C_{ctr-bd} \geq \frac{C_{T, \ell}}{\sqrt{2-\alpha}} .$$
This gives a first estimate of blow-up (that we will improve in the following).


\subsection{
A connection between null controllability and the existence of biorthogonal sequences
} \hfill

We notice the following fact: fix $m\geq 1$ and consider the initial condition $u_0 = \Phi _{\alpha ,m}$; let $H_{\alpha,m}$ be a control that drives the solution of \eqref{*pbm-controle1} to $0$ in time $T$; then the sequence $(r_{\alpha ,m} H_{\alpha,m})_{m\geq 1}$ is biorthogonal to $(e^{\lambda _{\alpha,n} t})_{n\geq 1}$ in $L^2(0,T)$. Indeed,
$H_{\alpha,m}$ satisfies \eqref{*moment-bd}:
$$\forall n \geq 1, \quad r_{\alpha ,n} \int _0 ^T H_{\alpha,m} (t) e^{\lambda _{\alpha,n} t} \, dt = \mu _{\alpha,n}^0 = \delta _{mn} ,$$
hence
$$ \forall n \geq 1, \quad \int _0 ^T \Bigl( r_{\alpha ,m} H_{\alpha,m} (t)\Bigr) e^{\lambda _{\alpha,n} t} \, dt = 
r_{\alpha ,m} \frac{\delta_{mn}}{r_{\alpha ,n}}
=
\begin{cases}
1 & \text{ if } m=n, \\
0 & \text{ if } m \neq n ,
\end{cases}$$
hence 
$$ \forall m,n \geq 1, \quad \int _0 ^T \Bigl( r_{\alpha ,m} H_{\alpha,m} (t)\Bigr) e^{\lambda _{\alpha,n} t} \, dt = \delta_{mn} ,$$
which means that the sequence $(r_{\alpha ,m} H_{\alpha,m})_{m\geq 1}$ is biorthogonal to $(e^{\lambda _{\alpha,n} t})_{n\geq 1}$ in $L^2(0,T)$.

In the literature, there exists several bounds from below for the biorthogonal families, we refer in particular to G\"uichal \cite{Guichal} and Hansen \cite{Hansen}. In the following we will use two extensions of the one of G\"uichal \cite{Guichal}, obtained in \cite{cost-weak} and in \cite{CMV-biortho-general}. 


\subsection{The concentration of the eigenvalues} \hfill

The following observation is fundamental in the understanding of the blow-up of the null controllability cost:

\begin{Lemma}
\label{lem-concentration}
The eigenvalues concentrate when $\alpha \to 2^-$:
$$ \forall n\geq 1, \quad \lambda_{\alpha, n+1} - \lambda _{\alpha,n} \to 0 \quad \text{ as } \alpha \to 2^- .$$
\end{Lemma}
Before proving Lemma \ref{lem-concentration}, let us explain why this property if clearly important in the understanding of the blow-up of the null controllability cost: as noted before, if null controllability holds, and if $H_{\alpha,m}$ is a control that drives the solution of \eqref{*pbm-controle1} with $u_0 = \Phi _{\alpha ,m}$ to $0$ in time $T$, then $(r_{\alpha ,m} H_{\alpha,m})_{m\geq 1}$ is biorthogonal to $(e^{\lambda _{\alpha,n} t})_{n\geq 1}$ in $L^2(0,T)$; now, if additionnally some eigenvalues concentrate, for example $\lambda_{\alpha, 2} - \lambda _{\alpha,1} \to 0$ as $\alpha \to 2^-$, then $r_{\alpha ,1} H_{\alpha,1}$ will have to satisfy
$$ \int _0 ^T \Bigl( r_{\alpha ,1} H_{\alpha,1} (t)\Bigr) e^{\lambda _{\alpha,1} t} \, dt = 1, \quad \text{ and } \quad \int _0 ^T \Bigl( r_{\alpha ,1} H_{\alpha,1} (t)\Bigr) e^{\lambda _{\alpha,2} t} \, dt = 0 ,$$
hence
$$ \int _0 ^T \Bigl( r_{\alpha ,1} H_{\alpha,1} (t)\Bigr) ( e^{\lambda _{\alpha,1} t} - e^{\lambda _{\alpha,2} t}) \, dt = 1 ;$$
but this will only be possible if $\Vert r_{\alpha ,1} H_{\alpha,1} \Vert$ is sufficiently large,
since $\Vert e^{\lambda _{\alpha,1} t} - e^{\lambda _{\alpha,2} t} \Vert _{L^2(0,T)}$ will be small. We will come back on this later.

\noindent {\it Proof of Lemma \ref{lem-concentration}.} 
We note that
$$\lambda_{\alpha, n+1} - \lambda _{\alpha,n} = \kappa _\alpha ^2 (j_{\nu_\alpha, n+1}^2 - j _{\nu_\alpha,n}^2 )
= \kappa _\alpha ^2 (j_{\nu_\alpha, n+1} - j _{\nu_\alpha,n} ) (j_{\nu_\alpha, n+1} + j _{\nu_\alpha,n} ). $$
It is classical (\cite{Kom-Lor} p. 135) that 
\begin{itemize}
\item if $\nu\in [0, \frac{1}{2}]$, the sequence $(j_{\nu, n+1} - j _{\nu,n})_n$ is nondecreasing and converges to $\pi$,
\item if $\nu\geq \frac{1}{2}$, the sequence $(j_{\nu, n+1} - j _{\nu,n})_n$ is nonincreasing and converges to $\pi$.
\end{itemize}
Then, when $\nu _\alpha \geq \frac{1}{2}$ (i.e. when $\alpha \in [\frac{4}{3},2)$), 
the sequence $(j_{\nu_\alpha, n+1} - j _{\nu_\alpha,n})_n$ is nonincreasing, hence
$$\lambda_{\alpha, n+1} - \lambda _{\alpha,n}
\leq \kappa _\alpha ^2 (j_{\nu_\alpha, 2} - j _{\nu_\alpha,1} ) (j_{\nu_\alpha, n+1} + j _{\nu_\alpha,n} ) ,$$
and using \eqref{*eq-Lorch},
$$\lambda_{\alpha, n+1} - \lambda _{\alpha,n}
\leq \kappa _\alpha ^2 (j_{\nu_\alpha, 2} - j _{\nu_\alpha,1} ) 
\Bigl( \pi (n+1 + \frac{\nu_\alpha}{2}-\frac{1}{4}) 
+ \pi (n + \frac{\nu_\alpha}{2}-\frac{1}{4}) \Bigr) .$$
Using \eqref{eq-Qu-Wong}, we obtain
\begin{multline*}
j_{\nu_\alpha, 2} - j _{\nu_\alpha,1} \leq \Bigl( \nu_\alpha - \frac{a_2}{2^{1/3}} \nu_\alpha ^{1/3}
+ \frac{3}{20} a_2 ^2 \frac{2^{1/3}}{\nu_\alpha ^{1/3}} \Bigr) - \Bigl( \nu_\alpha  - \frac{a_1}{2^{1/3}} \nu _\alpha ^{1/3} \Bigr) 
\\
= \frac{a_1 - a_2}{2^{1/3}} \nu_\alpha ^{1/3}
+ \frac{3}{20} a_2 ^2 \frac{2^{1/3}}{\nu_\alpha ^{1/3}} .
\end{multline*}
Hence there is some $C$ independent of $\alpha \in [\frac{4}{3},2)$ such that
\begin{equation}
\label{*gap12} 
j_{\nu_\alpha, 2} - j _{\nu_\alpha,1} \leq C  \nu_\alpha ^{1/3} ,
\end{equation}
and 
$$ \lambda_{\alpha, n+1} - \lambda _{\alpha,n} 
\leq C  \nu_\alpha ^{1/3} \kappa _\alpha ^2 (n + \nu_\alpha)
\leq C ( \kappa_\alpha ^{2/3} + \kappa_\alpha ^{5/3} n). \qed $$

\begin{Remark}
\em \label{re:conc_heat}
A similar concentration phenomenon can be pointed out in the fast control problem for the classical heat equation
\begin{equation}\label{eq:conc_heat}
\begin{cases}
u_t -  u_{xx} =h (x,t) \chi_{[a,b]}(x)  & \qquad x\in(0,1),\ 0<t<T,\\
u(0,t)=0=u(1,t), & \qquad 0<t<T, \\
u(x,0)=u_0(x), &\qquad x\in(0,1) ,\\
u(x,T)=0, &\qquad x\in(0,1) .
\end{cases}
\end{equation}
Indeed, as is well-known, the eigenvalues of the stationary operator associated with \eqref{eq:conc_heat} are $\lambda_n=\pi^2n^2$ for all $n>0$. On the other hand, if we are interested in studying the behaviour of the above system for controls yielding $u(\cdot,T)=0$ as $T\to 0^+$, then it might be useful to normalize the time, hence to look at the normalized solution
$$v(x,\tau)=u(x,\tau T) .$$ 
This function $v$ is solution of the problem
\begin{equation*}
\begin{cases}
v_\tau - T v_{xx} =T h(x,\tau T) \chi_{[a,b]}(x)  & \qquad x\in(0,1),\ 0<\tau<1,\\
v(0,\tau)=0=v(1,\tau), & \qquad 0<\tau<1, \\
v(x,0)=u_0(x), &\qquad x\in(0,1) ,\\
v(x,1)=0, &\qquad x\in(0,1) .
\end{cases}
\end{equation*}
Clearly, the eigenvalues of the stationary operator associated with this last problem are given by the sequence $\{T\pi^2n^2\}_{n\ge 1}$, which concentrates as  $T\to 0^+$.

\end{Remark}


\subsection{An additionnal property of the eigenvalues} \hfill

As we recalled, it is classical (\cite{Kom-Lor} p. 135) that 
\begin{itemize}
\item if $\nu\in [0, \frac{1}{2}]$, the sequence $(j_{\nu, n+1} - j _{\nu,n})_n$ is nondecreasing and converges to $\pi$,
\item if $\nu\geq \frac{1}{2}$, the sequence $(j_{\nu, n+1} - j _{\nu,n})_n$ is nonincreasing and converges to $\pi$.
\end{itemize}
Hence there exists a rank $N_\nu$ such that
$$ i \geq N_\nu \implies j_{\nu , i+1} - j_{\nu , i} \leq 2\pi .$$
However, the asymptotic development \eqref{eq-Qu-Wong} tells us that
\begin{equation}
\label{*p16}
 j_{\nu ,2} - j_{\nu , 1} \sim _{ \nu \to \infty} \frac{a_1-a_2}{2^{1/3}} \nu ^{1/3}.
 \end{equation}
Hence this rank $N_\nu$ probably satisfies $N_\nu \to +\infty$ as $\nu \to \infty$.
In the following, we estimate this $N_\nu$ (using the classical theory of Sturm concerning second order differential equations); we will need this estimate later.

\begin{Lemma}
\label{lem-Sturm}
Given $\nu \geq \frac{1}{2}$, then 
\begin{equation}
\label{*gap-loin}
\forall n > \nu, \quad j_{\nu, n+1} - j_{\nu,n} \leq 2 \pi .
\end{equation}
\end{Lemma}

\noindent {\it Proof of Lemma \ref{lem-Sturm}.} We follow and use the proofs of section 7.3 in \cite{Kom-Lor}: first we note that 
$$ y_\nu (x) := \sqrt{x} J_\nu (x) $$
satisfies the second-order differential equation
$$ y_\nu '' (x) + h_\nu (x) y_\nu (x) = 0 ,$$
with $$ h_\nu (x) = 1 - \frac{\nu ^2 - \frac{1}{4}}{x^2} .$$
Of course, $y_\nu$ and $J_\nu$ have the same positive zeros. We are going to use the following classical property of Sturm type (see Proposition 7.6 in \cite{Kom-Lor}): assume that 
\begin{itemize}
\item $f,g: [a,b] \to \Bbb R$ are continuous and satisfy
$$ \forall x \in [a,b], \quad f(x)<g(x) ,$$
\item $u,v$ are functions of class $C^2$ satisfying 
$$  \forall x \in [a,b], \quad u'' + fu=0, \quad v'' + gv =0 ,$$
\item $a,b$ are two consecutive zeros of $u$, 
\end{itemize}
then $v$ has at least one zero in $(a,b)$. 

We recall that, in a classical way (\cite{Kom-Lor}), this implies that $J_\nu$ has an infinite number of positive zeros: indeed: 
$$\forall x> \nu, \quad h_\nu (x) > \frac{1}{4\nu^2} ,$$
hence choosing 
\begin{gather*}
k\geq 1, \quad a:= 2k \nu \pi, \quad b:=2(k+1) \nu \pi, \\
f(x):= \frac{1}{4\nu^2}, \quad u(x):= \sin \frac{x}{2\nu} , \\
g(x):= h_\nu (x), \quad v(x):= y_\nu (x), 
\end{gather*}
 we can apply the Sturm property, and we derive that $y_\nu$ (hence $J_\nu$) has at least one zero on $ (2k \nu \pi, 2(k+1) \nu \pi)$. From \eqref{*eq-Lorch} we also have
$$ \forall k > \nu, \quad j_{\nu,k} > \pi (\nu + \frac{1}{4} (\nu - \frac{1}{2}) ) =:\gamma _\nu,$$
and then we can apply the Sturm property with
\begin{gather*}
 k> \nu, \quad a:= j_{\nu, k}, \quad b:=j_{\nu, k} + \frac{\pi}{\sqrt{h_\nu (\gamma _\nu)}}, \\
f(x):= h_\nu (\gamma _\nu), \quad u(x):= \sin \Bigl( \sqrt{h_\nu (\gamma _\nu)} (x-j_{\nu,k})\Bigr)  , \\
g(x):= h_\nu (x), \quad v(x):= y_\nu (x), 
\end{gather*}
and we deduce that $y_\nu$ has at least one zero inside $( j_{\nu, k}, j_{\nu, k} + \frac{\pi}{\sqrt{h_\nu (\gamma _\nu)}})$, hence
$$ j_{\nu,k+1} < j_{\nu, k} + \frac{\pi}{\sqrt{h_\nu (\gamma _\nu)}} .$$
Hence 
$$ \forall k >\nu, \quad  j_{\nu,k+1} - j_{\nu, k} < \frac{\pi}{\sqrt{h_\nu (\gamma _\nu)}}
= \frac{\pi}{\sqrt{1 - \frac{\nu ^2 - \frac{1}{4}}{\gamma _\nu ^2}}}.$$
It can be easily checked that 
$$ \forall \nu \geq \frac{1}{2}, \quad \frac{\pi}{\sqrt{1 - \frac{\nu ^2 - \frac{1}{4}}{\gamma _\nu ^2}}} \leq 2\pi :$$
indeed, if $\nu \geq \frac{1}{2}$, then
$$
\Bigl( 1 - \frac{\nu ^2 - \frac{1}{4}}{\gamma _\nu ^2} \Bigr) - \frac{1}{4}
= \frac{3}{4}- \frac{\nu ^2 - \frac{1}{4}}{\gamma _\nu ^2}
= \frac{3 \gamma _\nu ^2 - 4 (\nu ^2 - \frac{1}{4})}{4 \gamma _\nu ^2},$$
and
$$ 3 \gamma _\nu ^2 - 4 (\nu ^2 - \frac{1}{4})
= 3\Bigl( \pi (\nu + \frac{1}{4} (\nu - \frac{1}{2}) )\Bigr)  ^2 - 4(\nu ^2 - \frac{1}{4}) $$
and the discriminant of this quantity is negative, hence the quantity remains positive.
This implies \eqref{*gap-loin}. \qed


\subsection{A lower bound of the norm of any sequence biorthogonal to $(e^{\lambda _{\alpha,n}t})_n$ when $\nu _\alpha \in [0, \frac{1}{2}]$} \hfill

If $\nu _\alpha \in [0, \frac{1}{2}]$, the gap $(j_{\nu_\alpha, n+1} - j_{\nu_\alpha, n})_n$ is nondecreasing and converges to $\pi$, hence
$$ \forall n \geq 1, \quad j_{\nu_\alpha, n+1} - j_{\nu_\alpha, n} \leq \pi ,$$
hence
$$ \forall n \geq 1, \quad \sqrt{\lambda _{\alpha,n+1}} - \sqrt{\lambda _{\alpha,n}}
= \ell ^{\frac{\alpha}{2}-1} \kappa _\alpha (j_{\nu_\alpha, n+1} - j_{\nu_\alpha, n})
\leq \ell ^{\frac{\alpha}{2}-1} \kappa _\alpha \pi ,$$
hence
$$ \forall n \geq 1, \quad \sqrt{\lambda _{\alpha,n+1}} - \sqrt{\lambda _{\alpha,n}}
\leq \gamma _{max} \quad \text{ with } \gamma _{max } = \ell ^{\frac{\alpha}{2}-1} \kappa _\alpha \pi .$$

Let us apply the following extension of G\"uichal \cite{Guichal}:

\begin{Theorem} (Theorem 2.5 in \cite{cost-weak})
\label{thm-guichal-gen}
Assume that 
$$ \forall n\geq 0, \quad \lambda_n \geq 0, $$
and that there is some $0 < \gamma _{\text{min}} \leq \gamma _{\text{max}}$ such that
\begin{equation}
\label{gap-max}
\forall n \geq 0, \quad \gamma _{\text{min}} \leq \sqrt{\lambda _{n+1}} - \sqrt{\lambda _{n}}  \leq \gamma _{\text{max}} .
\end{equation}
Then there exists 
$c_u>0$ independent of $T$, and $m$ such that: any family $(\sigma _{m} ^+)_{m\geq 0}$ which is biorthogonal to the family $(e^{\lambda _{n}t})_{n\geq 0}$ in $L^2(0,T)$ 
satisfies:
\begin{equation}
\label{csq-gap-max}
\Vert \sigma _{m} ^+ \Vert _{L^2(0,T)} ^2
\geq e^{-2\lambda _{ m} T} e^{\frac{1}{2\gamma _{\text{max}} ^2 T}} b(T,\gamma_{max},m) ,
\end{equation}
with
\begin{equation}
\label{b-T-m-max}
b(T,\gamma_{max},m) = \frac{c_u ^2}{ C(m, \gamma _{\text{max}}, \lambda _0)^2 \, T } (\frac{1}{2\gamma _{\text{max}}^2 T})^{2m} \frac{1}{(4\gamma _{max} ^2T+1)^2} .
\end{equation}
and
\begin{equation}
\label{valeur-C}
C(m, \gamma _{max}, \lambda _0)
= m! \, 2 ^{m+ [\frac{2 \sqrt{\lambda _0}}{\gamma _{max}}] +1}
\, (m+ [\frac{2 \sqrt{\lambda _0}}{\gamma _{max}}] +1) .
\end{equation}
\end{Theorem}

Using Theorem \ref{thm-guichal-gen} with 
$\gamma _{max } = \ell ^{\frac{\alpha}{2}-1} \kappa _\alpha \pi$, 
one obtains that any family $(\sigma _m ^+)_{m\geq 1}$ which is biorthogonal to the family $(e^{\lambda _{\alpha ,n}t})_{n\geq1}$ in $L^2(0,T)$ satisfies a lower bound with the classical dominant exponential factor
of the type $e^{C/T}$:

$$ \Vert \sigma _{m} ^+ \Vert _{L^2(0,T)} ^2
\geq e^{-2\lambda _{\alpha, m} T} e^{C_u \frac{\ell^{2-\alpha}}{T \kappa _\alpha ^2}}
\frac{1}{T} \Bigl(\frac{\ell^{2-\alpha}}{T \kappa _\alpha ^2} \Bigr)^{2m}
\Bigl(\frac{\frac{\ell^{2-\alpha}}{T \kappa _\alpha ^2}}{1+\frac{\ell^{2-\alpha}}{T \kappa _\alpha ^2}} \Bigr)^2 \frac{C_u ^{2m}}{m!^2} 2 ^{-2\frac{\ell^{1-\alpha/2}}{ \kappa _\alpha}} 
\frac{1}{(m+1)^2 +  \frac{\lambda _{\alpha,1} \ell ^{1-\frac{\alpha}{2}} }{\kappa _\alpha ^2 \pi ^2}}.$$
This will immediately give an exponential blow-up of the cost as $T\to 0^+$, as explained in subsection \ref{subsec-expbl}, but the interesting behavior is when $\alpha \to 2^-$, and we study it in the following.


\subsection{A lower bound of the norm of any sequence biorthogonal to $(e^{\lambda _{\alpha,n}t})_n$ when $\nu _\alpha \geq \frac{1}{2}$} \hfill

This is the interesting case, where $\alpha \to 2^-$. 
In this case, the gap $(j_{\nu_\alpha, n+1} - j_{\nu_\alpha, n})_n$ is nonincreasing and converges to $\pi$, hence
$$ \forall n \geq 1, \quad \pi \leq j_{\nu_\alpha, n+1} - j_{\nu_\alpha, n} \leq j_{\nu_\alpha, 2} - j_{\nu_\alpha, 1} ,$$
hence
$$ \forall n \geq 1, \quad \sqrt{\lambda _{\alpha,n+1}} - \sqrt{\lambda _{\alpha,n}}
\leq \gamma _{max} \quad \text{ with } \gamma _{max } = \ell ^{\frac{\alpha}{2}-1} \kappa _\alpha (j_{\nu_\alpha, 2} - j_{\nu_\alpha, 1}) ;$$
but this time, we already noted that $j_{\nu_\alpha, 2} - j_{\nu_\alpha, 1}$ behaves as $\nu _\alpha ^{1/3}$ (see \eqref{*p16}), hence
$$ \gamma _{max } = c_\alpha \ell ^{\frac{\alpha}{2}-1} \kappa _\alpha ^{2/3} ,$$
with some uniformly bounded $c_\alpha$.

On the other hand, we proved in Lemma \ref{lem-Sturm} that
$$ \forall n \geq \nu_\alpha +1, \quad j_{\nu_\alpha, n+1} - j_{\nu_\alpha, n} \leq 2\pi ,$$
hence 
$$ \forall n \geq N_*, \quad \sqrt{\lambda _{\alpha,n+1}} - \sqrt{\lambda _{\alpha,n}}
\leq \gamma _{max}^*, \quad \text{ with } N_* = [\nu_\alpha ] +1, \text{ and } 
\gamma _{max}^* = 2\pi \ell ^{\frac{\alpha}{2}-1} \kappa _\alpha   .$$
Note that
$$ \frac{\gamma _{max}}{\gamma _{max} ^*} = \frac{c_\alpha}{2\pi \kappa _\alpha ^{1/3}} \to \infty \quad \text{ as } \alpha \to 2 ^- .$$
In that context, when there is a 'bad' global gap $\gamma _{max}$, and a 'good' (much smaller) asymptotic gap $\gamma _{max} ^*$, it is interesting to use the following extension of Theorem \ref{thm-guichal-gen}:

\begin{Theorem} (Theorem 2.2 in \cite{CMV-biortho-general})
\label{thm-guichal-gen*}
Assume that 
$$ \forall n\geq 1, \quad \lambda_n \geq 0, $$
and that there are $0<\gamma _{min} \leq \gamma _{\text{max}}^* \leq \gamma _{\text{max}}$ such that
\begin{equation}
\label{gap-max**}
\forall n \geq 1, \quad \gamma _{\text{min}} \leq \sqrt{\lambda _{n+1}} - \sqrt{\lambda _{n}}  \leq \gamma _{\text{max}} ,
\end{equation}
and
\begin{equation}
\label{gap-max*}
\forall n \geq N_*, \quad \sqrt{\lambda _{n+1}} - \sqrt{\lambda _{n}}  \leq \gamma _{\text{max}}^* .
\end{equation}
Then any family $(\sigma _{m} ^+)_{m\geq 1}$ which is biorthogonal to the family $(e^{\lambda _{n}t})_{n\geq 1}$ in $L^2(0,T)$
satisfies:
\begin{equation}
\label{csq-gap-max-bis}
\Vert \sigma _{m} ^+ \Vert _{L^2(0,T)} ^2
\geq e^{-2\lambda _{m} T} \, e^{\frac{2}{T (\gamma _{\text{max}}^*) ^2}} \, b^* (T,\gamma_{max},\gamma_{max} ^*, N_*, \lambda _1,m)^2 ,
\end{equation}
where $b^*$ is rational in $T$ (and explictly given in Lemma 4.4 of \cite{CMV-biortho-general}).
\end{Theorem}

Applying Theorem \ref{thm-guichal-gen*}, we obtain that
any family $(\sigma _{m} ^+)_{m\geq 1}$ which is biorthogonal to the family $(e^{\lambda _{\alpha, n}t})_{n\geq 1}$ in $L^2(0,T)$ satisfies
$$ \Vert \sigma _{m} ^+ \Vert _{L^2(0,T)} ^2
\geq e^{-2\lambda _{\alpha, m} T} \, e^{\frac{2}{4\pi ^2 \kappa _\alpha ^2 T \ell ^{\alpha -2}}} \, b^* (T,\gamma_{max},\gamma_{max} ^*, N_*, \lambda _{\alpha ,1},m)^2  $$
with an explicit value of $b^*$ (see Lemma 4.4 in \cite{CMV-biortho-general}): 
when $m \leq N_*$, we have
\begin{equation}
\label{min-distance}
b^*(T,\gamma_{max},\gamma_{max}^*,N_*, \lambda _1,m)
= C^* \frac{\sqrt{1+T\lambda _1}}{\sqrt{T}} \, \frac{(T \, (\gamma _{max} ^*) ^2)^{K_*+K' _*+2}}{(1+ (T \, (\gamma _{max} ^*) ^2))^{N_*+K_*+K' _*+3}} ,
\end{equation}
where 
$$ K_* = [\frac{2\sqrt{\lambda _1}+(N_*+m)\gamma _{max}}{\gamma_{max}^*}] - N_* + 2 ,$$
$$ K' _* =  [\frac{\gamma_{max}}{\gamma_{max}^*}(N_*-m)] -N_*+2 ,$$
$$ C^* = 
\frac{1}{(N_*+K_*+K'_*+3)!} 
\frac{c_u (\gamma _{max} ^*)^{2(N_*-1)} }{C^{(+)} C^{(-)}},$$
where
$$
C^{(+)} = (\frac{\gamma _{max}}{\gamma _{max} ^*}) ^{N_* -1} \frac{(N_*+m+ [\frac{2\sqrt{\lambda _1}}{\gamma_{max}}]+1)!}{(m+ [\frac{2\sqrt{\lambda _1}}{\gamma_{max}}]+1)! \, ([\frac{2\sqrt{\lambda _1}+(N_*+m)\gamma_{max}}{\gamma_{max}^*}]+1)! \, (2m+ [\frac{2\sqrt{\lambda _1}}{\gamma_{max}}]+1)}  ,
$$
and
$$ C^{(-)} =  (\frac{\gamma _{max}}{\gamma _{max} ^*}) ^{N_*-1} \, 
\frac{(m-1)! \, (N_*-m)! }{(1+ [\frac{\gamma_{max}}{\gamma_{max}^*}(N_*-m)])!} .
$$
These expressions seem be a little frightening, but we are looking for the behavior as $\alpha \to 2^-$, and this is not difficult to study:
one immediately sees that
$$ K_* + K' _* = c \nu _\alpha ^{4/3} + c' \ell ^{1-\alpha /2} \nu _\alpha + c'' \nu _\alpha ^{1/3} m ,$$
$$ (\gamma _{max} ^*)^{2(N_*-1)}
= e^{-c \nu _\alpha \ln \nu_\alpha - c' \nu _\alpha \ln \ell } ,$$
$$ \frac{1}{(N_*+K_*+K'_*+3)!}
\geq e^{-C (\nu_\alpha ^{4/3} + \ell ^{1-\alpha /2} \nu _\alpha + \nu _\alpha ^{1/3} m)(\ln \nu_\alpha + (1-\frac{\alpha}{2})\ln \ell + \ln m)},$$
and finally
$$ \frac{1}{C^{(+)} C^{(-)}} \geq e^{-c(\nu_\alpha +m)(\ln \nu _\alpha + \ln m)} \frac{1}{(m-1)!},$$
hence we obtain that
$$ C^* \geq e^{-C (\nu_\alpha ^{4/3} + \ell ^{1-\alpha /2} \nu _\alpha + \nu _\alpha ^{1/3} m)(\ln \nu_\alpha + (1-\frac{\alpha}{2})\ln \ell + \ln m)} \, \frac{1}{(m-1)!} .$$
This gives that
$$ b^*
\geq 
e^{-C (\nu_\alpha ^{4/3} + \ell ^{1-\alpha /2} \nu _\alpha + \nu _\alpha ^{1/3} m)(\ln \nu_\alpha + (1-\frac{\alpha}{2})\ln \ell + \ln m + \ln \frac{1}{T})} 
\, \frac{1}{(m-1)!} \, \frac{\sqrt{1+T}}{\sqrt{T}} ,
$$
hence
$$ \Vert \sigma _{m} ^+ \Vert _{L^2(0,T)} ^2 \geq \underline{b} (T,\alpha,m) ^2,
$$
with
\begin{multline}
\label{def-undb}
\underline{b} (T,\alpha,m)
:= e^{-\lambda _{\alpha, m} T}
 e^{C_u \frac{\ell^{2-\alpha}}{T \kappa _\alpha ^2}}
 \, \frac{1}{(m-1)!} \, \frac{\sqrt{1+T}}{\sqrt{T}}
\\
e^{-C (\nu_\alpha ^{4/3} + \ell ^{1-\alpha /2} \nu _\alpha + \nu _\alpha ^{1/3} m)(\ln \nu_\alpha + (1-\frac{\alpha}{2})\ln \ell + \ln m + \ln \frac{1}{T})} .
\end{multline}
This will give the expected blow-up of the cost, as $\alpha \to 2^-$ and/or as $T\to 0^+$.


\subsection{The exponential blow-up of the cost} \hfill
\label{subsec-expbl}

In the previous subsection, we obtained a bound from below for any biorthogonal sequence. But we already noted that if $u_0= \phi _{\alpha,m}$, and if $H_{\alpha, m}$ is any control that drives $u_0$ to rest in time $T$, then 
$(r_{\alpha ,m} H_{\alpha, m})_{m\geq 1}$ is biorthogonal to $(e^{\lambda _{\alpha, n} t})_{n\geq 1}$ in $L^2(0,T)$.  Hence
$$ \Vert r_{\alpha ,m} H_{\alpha, m} \Vert _{L^2(0,T)}
\geq \underline{b} (T,\alpha,m) ,$$ 
where $\underline{b} (T,\alpha,m)$ is given in \eqref{def-undb}.
By definition of the cost, we obtain that
$$ \forall m\geq 1, \quad C_{ctr-bd} \geq \frac{1}{\vert r_{\alpha ,m} \vert }\underline{b} (T,\alpha,m) ,$$
which gives the expected exponential blow-up of the cost:  choosing $m=1$, and using the fact that $\vert r_{\alpha,1} \vert = \sqrt{2 \kappa _\alpha \lambda _{\alpha ,1}} \ell ^{(\alpha -1)/2}$
and that $ \lambda _{\alpha ,1} = \kappa _\alpha j_{\nu _\alpha,1} \ell ^{(\alpha -2)/2}
\geq \frac{\pi}{4} \kappa _\alpha \nu _\alpha \ell ^{(\alpha -2)/2}$,
we obtain that there exists some $C_u$ independent of all the other parameters such that
\begin{multline*}
C_{ctr-bd} \geq C_u \frac{\ell ^{2-\alpha}}{\sqrt{\kappa _\alpha T \ell} }
e^{-\lambda _{\alpha ,1} T}
 e^{C_u \frac{\ell ^{2-\alpha}}{T \kappa _\alpha ^2}}
e^{-C (\nu_\alpha ^{4/3} + \ell ^{1-\alpha /2} \nu _\alpha )(\ln \nu_\alpha + \ln \ell ^{1-\alpha /2} + \ln \frac{1}{T})}  
\\
\geq 
C_u \frac{\ell ^{2-\alpha}}{\sqrt{\kappa _\alpha T \ell} }
e^{- \pi ^2 \frac{T}{\ell ^{2-\alpha}}}
 e^{C_u \frac{\ell ^{2-\alpha}}{T \kappa _\alpha ^2}}
e^{-C (\nu_\alpha ^{4/3} + \ell ^{1-\alpha /2} \nu _\alpha )(\ln \nu_\alpha + \ln \ell ^{1-\alpha /2} + \ln \frac{1}{T})}  ,
\end{multline*}
and this is \eqref{thm-cost-bd-low} and concludes the proof of Theorem \ref{thm-cost-bd-low}. \qed


\section{Proof of Theorem \ref{thm-cost-bd-up}}
\label{sec-Thm2}

We will use the following result (Theorem 2.4 in \cite{cost-weak}):

\begin{Theorem} (Existence of a suitable biorthogonal family and upper bounds)
\label{thm-biortho1-gen}
Assume that 
$$ \forall n\geq 0, \quad \lambda_n \geq 0, $$
and that there is some $\gamma _{\text{min}}>0$ such that
\begin{equation}
\label{gap-min}
\forall n \geq 0, \quad \sqrt{\lambda _{n+1}} - \sqrt{\lambda _{n}}  \geq \gamma _{\text{min}} .
\end{equation}
Then there exists a family $(\sigma _{m} ^+)_{m\geq 0}$ which is biorthogonal to the family $(e^{\lambda _{n}t})_{n\geq 0}$ in $L^2(0,T)$:
\begin{equation}
\label{*famillebi_qm-gen}
\forall m,n \geq 0, \quad \int _0 ^T \sigma _{m} ^+ (t)  e^{\lambda _{n}t} \, dt = \delta _{mn} .
\end{equation}
Moreover, it satisfies: there is some universal constant $C_u$ independent of $T$, $\gamma _{\text{min}}$ and $m$ such that, for all $m\geq 0$, we have
\begin{equation}
\label{*famillebi_qm-norme-gen}
\Vert \sigma _{m} ^+ \Vert _{L^2(0,T)} ^2  
\leq C_u e^{-2\lambda _{m} T} 
e^{C_u \frac{\sqrt{\lambda _{m}}}{ \gamma _{min}} }  
B(T,\gamma _{min}),
\end{equation}
with
\begin{equation}
\label{eq(B}
B(T,\gamma _{min}) =
\begin{cases} 
\Bigl( \frac{1}{T} + \frac{1}{T^2 \gamma _{min} ^2} \Bigr) 
\, e^{\frac{C_u}{\gamma _{min} ^2T}}  \quad & \text{ if } T \leq \frac{1}{\gamma _{min} ^2} , \\
C_u \gamma _{min} ^2 \quad & \text{ if } T \geq \frac{1}{\gamma _{min} ^2} .
\end{cases}
\end{equation}
\end{Theorem}

Note that \eqref{*famillebi_qm-norme-gen} and \eqref{eq(B} imply that there is some universal constant $C_u$ independent of $T$, $\gamma _{\text{min}}$ and $m$ such that, for all $m\geq 0$, we have
\begin{equation}
\label{*famillebi_qm-norme-gen**}
\Vert \sigma _{m} ^+ \Vert _{L^2(0,T)} ^2  
\leq C_u e^{-2\lambda _{m} T} 
e^{C_u \frac{\sqrt{\lambda _{m}}}{ \gamma _{min}} }  
e^{\frac{C_u}{\gamma _{min} ^2T}}
B^* (T,\gamma _{min}),
\end{equation}
with
\begin{equation}
\label{eq(B**}
B^* (T,\gamma _{min}) =
\frac{C_u}{T} \max \{T\gamma _{min} ^2 , \frac{1}{T\gamma _{min} ^2} \} .
\end{equation}

Now, as we have already noted, the eigenvalues of the problem satisfy
$$ \forall n\geq 1, \quad \sqrt{\lambda _{\alpha,n+1}} - \sqrt{\lambda _{\alpha,n}}
= \ell ^{\frac{\alpha}{2}-1} \kappa _\alpha (j_{\nu _\alpha,n+1} -  j_{\nu _\alpha,n})
\geq \begin{cases} 
 \ell ^{\frac{\alpha}{2}-1} \kappa _\alpha (j_{\nu _\alpha,2} -  j_{\nu _\alpha,1}) \text{ if } \nu _\alpha \in [0, \frac{1}{2}], \\
\ell ^{\frac{\alpha}{2}-1} \kappa _\alpha \pi \text{ if } \nu_\alpha \geq \frac{1}{2} 
\end{cases}
.$$
Define artificially 
$$ \lambda _{\alpha,0} :=0 .$$
Then 
$$ \sqrt{\lambda _{\alpha,1}} - \sqrt{\lambda _{\alpha,0}}
= \ell ^{\frac{\alpha}{2}-1} \kappa _\alpha j_{\nu _\alpha,1} .$$
Then consider
$$ c_\alpha := \begin{cases} 
 \min \{ j_{\nu _\alpha,2} -  j_{\nu _\alpha,1}, j_{\nu _\alpha,1}\} \text{ if } \nu _\alpha \in [0, \frac{1}{2}], \\
\min \{ \pi , j_{\nu _\alpha,1} \} \text{ if } \nu_\alpha \geq \frac{1}{2} 
\end{cases} ,
$$
and 
$$ \underline{c} := \inf _{\alpha \in [0,2) } c_\alpha .$$
It is clear from \eqref{*eq-Lorch}-\eqref{*eq-Lorch2} that $ \underline{c} >0$,
and by construction we have
$$  \forall n\geq 0, \quad \sqrt{\lambda _{\alpha,n+1}} - \sqrt{\lambda _{\alpha,n}}
\geq \gamma _{min} \quad \text{ with } \quad \gamma _{min} := \ell ^{\frac{\alpha}{2}-1} \kappa _\alpha \underline{c} .$$

Then, applying Theorem \ref{thm-biortho1-gen} with $\gamma _{min} = \ell ^{\frac{\alpha}{2}-1} \kappa _\alpha \underline{c}$, we obtain that there exists
a family $(\sigma _{\alpha,m} ^+)_{m\geq 0}$ biorthogonal to $(e^{\lambda _{\alpha,n}t})_{n\geq 0}$ in $L^2(0,T)$, and such that
\begin{multline*}
\Vert \sigma _{\alpha,m} ^+ \Vert _{L^2(0,T)} ^2  
\leq C_u e^{-2\lambda _{\alpha,m} T} 
e^{C_u \frac{\sqrt{\lambda _{\alpha,m}}}{ \gamma _{min}} } 
B(T,\gamma _{min}) 
\\
= C_u e^{-2\lambda _{\alpha,m} T} e^{C_u j _{\nu_\alpha,m}} B(T,\gamma _{min}).
\end{multline*}

 Then define
\begin{equation}
\label{*def-controle}
K (t) := - \sum _{m=1} ^\infty \frac{\lambda _{\alpha,m} \mu_{\alpha ,m} ^0}{r_{\alpha ,m}} \sigma_{\alpha, m}^+ (t), \quad \text{ and } \quad H (t) := \int _0 ^t K (\tau) \, d\tau ,
\end{equation}
and let us check that $H$ is an admissible control that drives the solution of \eqref{*pbm-controle1} to $0$ in time $T$:
\begin{itemize}
\item first we check that $K \in L^2 (0,T)$: using \eqref{*infl-fp2} and \eqref{*famillebi_qm-norme-gen}, we have
$$
\sum _{m=1} ^\infty \frac{\vert \lambda _{\alpha,m} \mu_{\alpha ,m} ^0\vert}{\vert r_{\alpha ,m}\vert} \Vert \sigma_{\alpha, m}^+ \Vert _{L^2(0,T)}
\\
\leq \Bigl( \sum _{m=1} ^\infty  \vert \mu_{\alpha ,m} ^0\vert ^2 \Bigr) ^{1/2} 
\Bigl( \sum _{m=1} ^\infty \frac{\vert \lambda _{\alpha,m} \vert ^2 }{\vert r_{\alpha ,m}\vert ^2 } \Vert \sigma_{\alpha, m}^+ \Vert _{L^2(0,T)} ^2 \Bigr) ^{1/2}
$$
which is finite (we will come back on this in the following);
this implies that $H \in H^1(0,T)$, and of course $H(0)=0$, and also $H(T)=0$ using \eqref{*famillebi_qm-gen} with $n=0$;
\item next, we check that $H$ satisfies the moment problem \eqref{*moment-bd-G'}:
\begin{multline*}
\forall n\geq 1, \quad -  \frac{r_{\alpha ,n}}{\lambda _{\alpha,n}} \int _0 ^T H'(t)  e^{\lambda _{\alpha,n}t} \, dt
+ \frac{r_{\alpha ,n}}{\lambda _{\alpha,n}}  \Bigl[ H(T)  e^{\lambda _{\alpha,n}T} - H(0) \Bigr]
\\
= -  \frac{r_{\alpha ,n}}{\lambda _{\alpha,n}} \int _0 ^T K(t)  e^{\lambda _{\alpha,n}t} \, dt
= \mu_{\alpha ,n} ^0 ;
\end{multline*}
\item finally we check that the solution of \eqref{*pbm-controle1} satisfies $u(T)=0$: multiplying the first equation of \eqref{*pbm-controle1} by $w_{\alpha ,n} (x,t):= \Phi_{\alpha ,n} (x) e^{\lambda_{\alpha ,n} (t-T)}$ and integrating by parts, we obtain that
$$ \forall n\geq 1, \quad \int _0 ^\ell u(x,T) \Phi_{\alpha ,n} (x) \, dx = 0 ,$$
hence $u(T)=0$.
\end{itemize}
Hence $H$ is an admissible control, and therefore
$$ C _{ctr-bd} \leq \frac{\Vert H \Vert _{H^1 (0,T)}}{\Vert u_0 \Vert _{L^2 (0,\ell)} } \leq C \frac{\Vert K \Vert _{L^2 (0,T)}}{\Vert u_0 \Vert _{L^2 (0,\ell)} } ,$$
hence
\begin{equation*}
C _{ctr-bd} 
\leq C \Bigl( \sum _{m=1} ^\infty \frac{\vert \lambda _{\alpha,m} \vert ^2 }{\vert r_{\alpha ,m}\vert ^2 } \Vert \sigma_{\alpha,m}^+ \Vert _{L^2(0,T)} ^2 \Bigr) ^{1/2}.
\end{equation*}
Since $ r_{\alpha ,m} ^2 = 2 \kappa _\alpha \ell ^{\alpha -1} \lambda _{\alpha,m}$,
we have 
$$ \frac{\vert \lambda _{\alpha,m} \vert ^2 }{\vert r_{\alpha ,m}\vert ^2 }
= \frac{\kappa _\alpha j_{\nu_\alpha ,m} ^2 }{2 \ell} ,$$
hence
\begin{multline*} 
C _{ctr-bd} 
\leq \frac{C}{\sqrt{\ell}} \sqrt{\kappa _\alpha} \Bigl( \sum _{m=1} ^\infty j_{\nu_\alpha ,m} ^2 \Vert \sigma_{\alpha,m}^+ \Vert _{L^2(0,T)} ^2 \Bigr) ^{1/2} 
\\
\leq \frac{C}{\sqrt{\ell}} \sqrt{\kappa _\alpha} \Bigl( \sum _{m=1} ^\infty j_{\nu_\alpha ,m} ^2 e^{-2\lambda _{\alpha,m} T} 
e^{C_u \frac{\sqrt{\lambda _{\alpha,m}}}{ \gamma _{min}} }  
B(T,\gamma _{min}) \Bigr) ^{1/2} 
\\
= \frac{C}{\sqrt{\ell}} \sqrt{\kappa _\alpha} \sqrt{B(T,\gamma _{min})} \Bigl( \sum _{m=1} ^\infty j_{\nu_\alpha ,m} ^2 e^{-2\lambda _{\alpha,m} T} 
e^{C_u \frac{\sqrt{\lambda _{\alpha,m}}}{ \gamma _{min}} } \Bigr) ^{1/2} 
\end{multline*}
But
$$ C_u \frac{\sqrt{\lambda _{\alpha,m}}}{ \gamma _{min}} \leq 
\lambda _{\alpha,m}T + \frac{C_u ^2}{T \gamma _{min} ^2} ,$$
hence
\begin{multline}
\label{*20oct7}
C _{ctr-bd} 
\leq \frac{C}{\sqrt{\ell}} \sqrt{\kappa _\alpha} \sqrt{B(T,\gamma _{min})} 
e^{\frac{C_u ^2}{2T \gamma _{min} ^2}}
\Bigl( \sum _{m=1} ^\infty j_{\nu_\alpha ,m} ^2 e^{-\lambda _{\alpha,m} T} 
\Bigr) ^{1/2}  
\\
= \frac{C}{\sqrt{\ell}} \sqrt{\kappa _\alpha} \sqrt{B(T,\gamma _{min})} 
e^{\frac{C_u ^2 }{2 \underline{c} ^2} \frac{\ell ^{2-\alpha}}{T \kappa _\alpha ^2}}
\Bigl( \sum _{m=1} ^\infty j_{\nu_\alpha ,m} ^2 e^{-j_{\nu_\alpha,m} ^2 \frac{\kappa _\alpha ^2 T }{\ell^{2-\alpha}}} 
\Bigr) ^{1/2}
\end{multline}
It remains to estimate the last sum. We distinguish the cases $\nu_\alpha \leq \frac{1}{2}$ and $\nu_\alpha \geq \frac{1}{2}$. 

Take $Y>0$. When $\nu_\alpha \leq \frac{1}{2}$, using \eqref{*eq-Lorch2} we see that
$$ j_{\nu_\alpha,m} ^2 e^{-j_{\nu_\alpha,m} ^2 Y}
\leq \pi ^2 (m+\frac{\nu_\alpha -\frac{1}{2}}{4} )^2 e^{- Y \pi ^2 (m+\frac{\nu_\alpha -\frac{1}{2}}{2} )^2} .
$$
When $m\geq 1$, we have
$$ (m+\frac{\nu_\alpha -\frac{1}{2}}{2} )^2 \geq (m-\frac{1}{4})^2 \geq \frac{1}{2} m^2 + \frac{1}{16},$$
hence
$$ j_{\nu_\alpha,m} ^2 e^{-j_{\nu_\alpha,m} ^2 Y}
\leq \pi ^2 m^2 e^{- Y \pi ^2 m^2 /2}  e^{- Y \pi ^2 /16}.
$$
The function $V: x\mapsto \pi ^2 x^2 e^{- Y \pi ^2 x^2 /2}$ attains its maximum at $x_Y := \sqrt{\frac{2}{\pi^2Y}}$,
is increasing on $[0, x_Y]$,
decreasing on $[x_Y, +\infty)$, and its maximum is $\frac{2}{eY}$. If $x_Y \leq 1$, then 
\begin{multline*}
\sum _{m=1} ^\infty \pi ^2 m^2 e^{- Y \pi ^2 m^2 /2}
= \sum _{m=1} ^\infty V(m)  = V(1) + \sum _{m=1} ^\infty V(m+1) 
\\
\leq V(1) + \sum _{m=1} ^\infty \int _m ^{m+1} V(x) \, dx = V(1) + \int _1 ^\infty V(x) \, dx .
\end{multline*}
If $x_Y \geq 1$, then 
\begin{multline*}
\sum _{m=1} ^\infty \pi ^2 m^2 e^{- Y \pi ^2 m^2 /2}
= \sum _{m=0} ^\infty V(m)  
\\
= \sum _{m=0} ^{[x_Y]-1}  V(m) + V([x_Y]) + V([x_Y]+1) + \sum _{[x_Y]+2} ^\infty  V(m)
\\
\leq \sum _{m=0} ^{[x_Y]-1}   \int _m ^{m+1} V(x) \, dx + V([x_Y]) + V([x_Y]+1) + \sum _{[x_Y]+2} ^\infty  \int _{m-1} ^{m} V(x) \, dx
\\
\leq V([x_Y]) + V([x_Y]+1) + \int _0 ^\infty V(x) \, dx
\leq \frac{4}{eY} + \int _0 ^\infty V(x) \, dx .
\end{multline*}
Hence in any case, 
\begin{multline}
\label{*20oct-4}
\sum _{m=1} ^\infty \pi ^2 m^2 e^{- Y \pi ^2 m^2 /2}
\leq \frac{4}{eY} + \int _0 ^\infty \pi ^2 x^2 e^{- Y \pi ^2 x^2 /2} \, dx
\\
= \frac{4}{eY} + \frac{2^{3/2}}{\pi Y^{3/2}} \int _0 ^\infty s^2 e^{-s^2} \, ds ,
\end{multline}
and 
$$ \sum _{m=1} ^\infty j_{\nu_\alpha,m} ^2
e^{-j_{\nu_\alpha,m} ^2 Y} 
\leq \Bigl( \frac{4}{eY} + \frac{2^{3/2}}{\pi Y^{3/2}} \int _0 ^\infty s^2 e^{-s^2} \, ds \Bigr) e^{- Y \pi ^2 /16} .$$
hence there exists some $C_u$ such that, when $\nu_\alpha \leq \frac{1}{2}$,
\begin{equation}
\label{*20oct-3}
\forall Y>0, \quad  \sum _{m=1} ^\infty j_{\nu_\alpha,m} ^2
e^{-j_{\nu_\alpha,m} ^2 Y} 
\leq \frac{1}{C_u  Y ^{3/2}} e^{-C_u Y } .
\end{equation}

When $\nu_\alpha \geq \frac{1}{2}$, we proceed in the same way, using \eqref{*eq-Lorch}: we see that
$$ j_{\nu_\alpha,m} ^2 e^{-j_{\nu_\alpha,m} ^2 Y}
\leq \pi ^2 (m+\frac{\nu_\alpha -\frac{1}{2}}{2} )^2 e^{- Y \pi ^2 (m+\frac{\nu_\alpha -\frac{1}{2}}{4} )^2} .$$
But
$$  (m+\frac{\nu_\alpha -\frac{1}{2}}{4} )^2 \geq \frac{1}{2}m^2 + (\frac{\nu_\alpha -\frac{1}{2}}{4} )^2 ,$$
hence
\begin{multline*}
j_{\nu_\alpha,m} ^2 e^{-j_{\nu_\alpha,m} ^2 Y} \leq \pi ^2 (m+\frac{\nu_\alpha -\frac{1}{2}}{2} )^2  e^{- Y \pi ^2 m^2/2}  e^{- Y \pi ^2 (\frac{\nu_\alpha -\frac{1}{2}}{4} )^2 }
\\
\leq 2 \pi ^2 \Bigl( m^2 + (\frac{\nu_\alpha -\frac{1}{2}}{2})^2 \Bigr) e^{- Y \pi ^2 m^2/2}  e^{- Y \pi ^2 (\frac{\nu_\alpha -\frac{1}{2}}{4} )^2 } 
\\ 
\leq 2 \Bigl( 1 + (\frac{\nu_\alpha -\frac{1}{2}}{2})^2 \Bigr) e^{- Y \pi ^2 (\frac{\nu_\alpha -\frac{1}{2}}{4} )^2 } 
\pi ^2 m^2 e^{- Y \pi ^2 m^2/2}  .
\end{multline*}
Hence, using \eqref{*20oct-4}, we obtain
$$  \sum _{m=1} ^\infty j_{\nu_\alpha,m} ^2 e^{-j_{\nu_\alpha,m} ^2 Y} 
\leq 2 \Bigl( 1 + (\frac{\nu_\alpha -\frac{1}{2}}{2})^2 \Bigr) e^{- Y \pi ^2 (\frac{\nu_\alpha -\frac{1}{2}}{4} )^2 } 
\Bigl( \frac{4}{eY} + \frac{2^{3/2}}{\pi Y^{3/2}} \int _0 ^\infty s^2 e^{-s^2} \, ds \Bigr) , $$
hence there exists some $C_u$ such that, when $\nu_\alpha \geq \frac{1}{2}$,
\begin{equation}
\label{*20oct-5}
\forall Y>0, \quad   \sum _{m=1} ^\infty j_{\nu_\alpha,m} ^2
e^{-j_{\nu_\alpha,m} ^2 Y} 
\leq \frac{\nu_\alpha ^2}{C_u  Y^{3/2}} e^{-C_u \nu_\alpha ^2 Y} .
\end{equation}
And then, there is some $C_u$ independent of $\alpha \in [1,2)$ of $T>0$ and of $\ell >0$ such that
\begin{equation}
\label{*20oct-6}
 \sum _{m=1} ^\infty j_{\nu_\alpha,m} ^2
e^{-j_{\nu_\alpha,m} ^2 \Bigl(\frac{\kappa _\alpha ^2 T }{\ell^{2-\alpha}}\Bigr)} 
\leq \frac{1}{C_u  \kappa_\alpha ^2 \Bigl(\frac{\kappa _\alpha ^2 T }{\ell^{2-\alpha}}\Bigr)^{3/2}} e^{-\frac{C_u}{\kappa_\alpha ^2} \Bigl(\frac{\kappa _\alpha ^2 T }{\ell^{2-\alpha}}\Bigr)} .
\end{equation}
That allows us to complete the estimate from above of the null controllability cost: we deduce from \eqref{*20oct7}, \eqref{*20oct-6}, \eqref{eq(B} that

$$
C _{ctr-bd} 
\leq C'_u 
\frac{1}{\sqrt{\kappa_\alpha T \ell}} 
e^{C'_u \frac{\ell^{2-\alpha}}{\kappa _\alpha ^2 T }} 
e^{-\frac{C_u}{2} \frac{T }{\ell^{2-\alpha}} }  ,
$$
which is \eqref{*borne-bd-up}. This completes the proof of Theorem \ref{thm-cost-bd-up}. \qed



\section{Proof of Theorem \ref{thm-cost-ld-up}}
\label{sec-Thm4}

In section \ref{sec-mmloc}, we constructed, at least formally, a control that drives the initial condition $u_0$ to $0$ in time $T$. This control is given by \eqref{*contr}, and depends of a suitable biorthogonal family $ \sigma_{\alpha ,m}^+$ satisfying \eqref{biortho+1}, and of the norm of the eigenfunctions in the control region.
Theorem \ref{thm-biortho1-gen} (in fact \eqref{*famillebi_qm-norme-gen**} and \eqref{eq(B**}) gives the existence and bounds for a biorthogonal family $ (\sigma_{\alpha ,m}^+)_{m\geq 1}$ satisfying \eqref{biortho+1}. Proposition \ref{prop-fctpr} gives an estimate of the norm of the eigenfunctions in the control region (and will be proved in section \ref{sec-complfp}). Here we use these results to prove Theorem \ref{thm-cost-ld-up}: using Theorem \ref{thm-biortho1-gen} and Proposition \ref{prop-fctpr}, we have
\begin{multline*}
\Vert \sigma_{\alpha ,m}^+ \Vert ^2 _{L^2(0,T)} \frac{1}{\Bigl( \int_a^b \Phi_{\alpha ,m}^2 \Bigr) ^2}
\leq C_u B^* (T,\gamma_{min}) e^{C_u \frac{\ell^{2-\alpha}}{\kappa _\alpha ^2 T }} e^{- \lambda _{\alpha,m} T}
\frac {1}{\Bigl( \gamma _0 ^* (2-\alpha) \Bigr)^2}
\\
\leq \frac{C_u B^* (T,\gamma_{min}) }{\vert \gamma _0 ^* \vert ^2 (2-\alpha)^2} e^{C_u \frac{\ell^{2-\alpha}}{\kappa _\alpha ^2 T }} e^{- \kappa _\alpha ^2 j_{\nu_\alpha,1} ^2 \frac{T}{\ell ^{2-\alpha}}} .
\end{multline*}
Hence, there is some $C_u$ independent of $T>0$, $\ell >0$, $\alpha \in [1,2)$, $m\geq 1$  such that
$$ \Vert \sigma_{\alpha ,m}^+ \Vert ^2 _{L^2(0,T)} \frac{1}{\Bigl( \int_a^b \Phi_{\alpha ,m}^2 \Bigr) ^2}
\leq \frac{C_u B^* (T,\gamma_{min})}{\vert \gamma _0 ^* \vert ^2  (2-\alpha)^2} e^{C_u \frac{\ell^{2-\alpha}}{\kappa _\alpha ^2 T }} e^{- \frac{1}{C_u} \frac{T}{\ell ^{2-\alpha}}} .$$
Of course, if $(\mu_{\alpha ,m} ^0)_m \in \ell ^2(\Bbb N)$, then the series
$$ \sum_{m\geq 1} \vert \mu_{\alpha ,m} ^0 \vert ^2 \Vert \sigma_{\alpha ,m}^+ \Vert ^2 _{L^2(0,T)} \frac{1}{\Bigl( \int_a^b \Phi_{\alpha ,m}^2 \Bigr) ^2} $$
is convergent. Hence the control given by the formula \eqref{*contr} is in $L^2((0,\ell) \times (0,T))$, and 
\begin{multline*}
\Vert h \Vert _{L^2((0,\ell)\times (0,T))} ^2 = \sum_{m\geq 1} \vert \mu_{\alpha ,m} ^0 \vert ^2 \Vert \sigma_{\alpha ,m}^+ \Vert ^2 _{L^2(0,T)} \frac{1}{\Bigl( \int_a^b \Phi_{\alpha ,m}^2 \Bigr) ^2}  
\\
\leq \frac{C_u B^* (T,\gamma_{min})}{\vert \gamma _0 ^* \vert ^2  (2-\alpha)^2} e^{C_u \frac{\ell^{2-\alpha}}{\kappa _\alpha ^2 T }} e^{- \frac{1}{C_u} \frac{T}{\ell ^{2-\alpha}}} \sum_{m\geq 1} \vert \mu_{\alpha ,m} ^0 \vert ^2 .
\end{multline*}
Hence
$$ C_{ctr-loc} ^2 \leq \frac{C_u B^* (T,\gamma_{min})}{ \vert \gamma _0 ^* \vert ^2 (2-\alpha)^2} e^{C_u \frac{\ell^{2-\alpha}}{\kappa _\alpha ^2 T }} e^{- \frac{1}{C_u} \frac{T}{\ell ^{2-\alpha}}} .$$
This gives \eqref{*borne-loc-up} (with another constant $C_u$). In particular, note that the dependence in the control region appears only in $\gamma  ^* _0 $. \qed


\section{Proof of Theorem \ref{thm-cost-ld-low}}
\label{sec-Thm3}

Given $u_0 \in L^2(0,\ell)$, assume that $h\in L^2((a,b)\times (0,T))$ is a control that drives the solution of \eqref{*pbm-controle2} to $0$ in time $T$. Denote
$$ H (t) := u(a,t) .$$
Then the function $u$ satisfies
\begin{equation}
\label{*pbm-controle2-bd}
\begin{cases}
u_t - (x^\alpha u_x)_x =0  & \qquad x\in(0,a),\ t>0,\\
(x^\alpha u_x) (0,t)=0, & \qquad t>0, \\
u(a,t)= H (t)  & \qquad t>0, \\
u(x,0)=u_0(x), &\qquad x\in(0,a) 
\end{cases}
\end{equation}
and 
$$ u(x,T)=0, \quad x\in(0,a) ,$$
hence $H$ is a boundary control that drives the solution of \eqref{*pbm-controle2-bd} to $0$ in time $T$. Let us choose $m\geq 1$ and 
$$ u_0 (x) := 
\begin{cases}
 \frac{\sqrt{2 \kappa _\alpha }}{a^{\kappa _\alpha} \vert J'_{\nu_\alpha} (j_{\nu_\alpha,m} ) \vert} 
x^{(1-\alpha)/2} J_{\nu _\alpha} (j_{\nu_\alpha,m} (\frac{x}{a}) ^{\kappa_\alpha}) , & \qquad x\in(0,a),\\
0, & \qquad x\in(a,\ell) ,
\end{cases}$$
in such a way that the initial condition of \eqref{*pbm-controle2-bd} is exactly an eigenfunction of the associated Sturm-Liouville problem. Then we know from subsection \ref{subsec-expbl} that

\begin{equation}
\label{*vendredi16oct2}
\vert r_{\alpha ,m} \vert  \Vert u(a, \cdot )  \Vert _{L^2(0,T)} 
= \Vert r_{\alpha ,m} H_m \Vert _{L^2(0,T)}
\geq \underline{b} (T,\alpha,m) ,
\end{equation}
where $\underline{b}$ is defined in \eqref{def-undb}, but where $a$ replaces $\ell$ in the expressions of $r_{\alpha ,m}$ and $\underline{b} (T,\alpha,m)$.

On the other hand, energy methods tell us that the control $h$ and the initial condition dominate the solution of \eqref{*pbm-controle2}: indeed, first we have
$$ \forall y \geq a, \quad -u(y,t) = \int _y ^\ell u_x (x,t) \, dx ,$$
hence
$$ u(y,t) ^2 = \Bigl( \int _y ^\ell u_x (x,t) \, dx \Bigr) ^2 
\leq \Bigl( \int _y ^\ell x^\alpha u_x ^2 (x,t) \, dx \Bigr) \Bigl( \int _y ^\ell x^{-\alpha} \, dx \Bigr)
,$$
hence
$$ \forall y \in [a, \ell), \quad u(y,t) ^2 \leq C(\alpha, a, \ell) \int _0 ^\ell x^\alpha u_x ^2 (x,t) \, dx $$
with 
$$ C(\alpha, a, \ell) = 
\begin{cases}
\frac{1}{(\alpha -1) a^{\alpha -1}} & \text{ if } \alpha \in (1,2), \\
\ln \frac{\ell}{a} & \text{ if } \alpha =1 
\end{cases} .
$$
Then, multiplying the first equation of \eqref{*pbm-controle2} by $u$, we have
$$ 
\int _0 ^T \int _0 ^\ell u h \chi _{(a,b)}
= \int _0 ^T \int _0 ^\ell u (u_t - (x^\alpha u_x)_x )
= -\frac{1}{2} \int _0 ^\ell u_0 ^2 + \int _0 ^T \int _0 ^\ell x^\alpha u_x ^2 ,
$$
hence
\begin{multline*}
 \int _0 ^T \int _0 ^\ell x^\alpha u_x ^2
= \frac{1}{2} \int _0 ^\ell u_0 ^2 + \int _0 ^T \int _a ^b u h 
\\
\leq \frac{1}{2} \int _0 ^\ell u_0 ^2 +   \int _0 ^T \int _a ^b  \Bigl( C(\alpha, a, \ell) \int _0 ^\ell x^\alpha u_x ^2 (x,t) \, dx \Bigr) ^{1/2} \vert h \vert
\\
\leq \frac{1}{2} \int _0 ^\ell u_0 ^2 + \frac{1}{2} \int _0 ^T  \int _0 ^\ell x^\alpha u_x ^2 (x,t) \, dx \, dt 
+ 
\frac{(b-a) C(\alpha, a, \ell)}{2} \int _0 ^T  \int _a ^b h(x,t)^2 \, dx \, dt .
\end{multline*}
We obtain that
$$
 \int _0 ^T \int _0 ^\ell x^\alpha u_x ^2 
\leq  \int _0 ^\ell u_0 ^2 
+ (b-a) C(\alpha, a, \ell) \int _0 ^T  \int _a ^b h(x,t)^2 \, dx \, dt ,$$
hence
\begin{multline*}
 \int _0 ^T u(a,t)^2 \, dt 
\leq C(\alpha,a,\ell) \int _0 ^T \int _0 ^\ell x^\alpha u_x ^2
\\
\leq C(\alpha,a,\ell) \int _0 ^\ell u_0 ^2 
+ (b-a) C(\alpha, a, \ell)^2  \int _0 ^T  \int _a ^b h(x,t)^2 \, dx \, dt .
\end{multline*}
The initial condition $u_0$ that we have chosen has an $L^2$-norm equal to $1$, hence
$$
\int _0 ^T  \int _a ^b h(x,t)^2 \, dx \, dt
\geq \frac{1}{(b-a) C(\alpha, a, \ell)^2} \int _0 ^T u(a,t)^2 \, dt - \frac{1}{(b-a) C(\alpha, a, \ell)}  ,$$
and the lower bound \eqref{*vendredi16oct2} of $\Vert u(a, \cdot ) \Vert _{L^2(0,T)}$ implies that 
$$ 
\int _0 ^T  \int _a ^b h(x,t)^2 \, dx \, dt
\geq \frac{1}{(b-a) C(\alpha, a, \ell)^2} 
\frac{\underline{b} (T,\alpha,m)^2}{ r_{\alpha ,m} ^2 } 
- \frac{1}{(b-a) C(\alpha, a, \ell)} .$$

As we did in subsection \ref{subsec-expbl}, choosing $m=1$, this implies that
\begin{multline*}
\int _0 ^T  \int _a ^b h(x,t)^2 \, dx \, dt
\\
\geq 
\frac{1}{(b-a) C(\alpha, a, \ell)^2}
\Bigl[ 
C_u \frac{a ^{2(2-\alpha)}}{\kappa _\alpha T a }
e^{- 2\pi ^2 \frac{T}{a ^{2-\alpha}}}
 e^{2C_u \frac{a ^{2-\alpha}}{T (2-\alpha) ^2}}
\\
e^{-\frac{2}{C_u} (\frac{1}{(2-\alpha)^{4/3}} +  \frac{a ^{1-\alpha /2}}{2-\alpha} )(\ln \frac{a ^{1-\alpha /2}}{2-\alpha} + \ln \frac{1}{T})}\Bigr] - \frac{1}{(b-a) C(\alpha, a, \ell)} .
\end{multline*}
Then the null controllability
cost for \eqref{*pbm-controle2} blows up at least exponentially fast when $\alpha \to 2^-$, as stated in Theorem \ref{thm-cost-ld-low}. 
One can note that the bound from below is very poor when $\ell$ is large. But this is due to the method: indeed, we concentrate 
the initial condition on the zone at the left of the control region.
\qed 


\section{The eigenfunctions in the control region (Proposition \ref{prop-fctpr})}
\label{sec-complfp}

The goal of this section is to prove Proposition \ref{prop-fctpr}, that was be useful to prove Theorem \ref{thm-cost-ld-up}.


\subsection{The reduction to an ordinary differential equation question} \hfill

Using Proposition \ref{*prop-vp}, we note that 
\begin{multline*}
\int_a^b \Phi_{\alpha ,m}(x) ^2 \, dx 
= \int_a^b  \frac{2 \kappa _\alpha }{\ell^{2\kappa _\alpha} \vert J'_{\nu_\alpha} (j_{\nu_\alpha,m} ) \vert ^2} 
x^{1-\alpha} J_{\nu _\alpha} (j_{\nu_\alpha,m} (\frac{x}{\ell}) ^{\kappa_\alpha}) ^2 \, dx 
\\
=\frac{2 \kappa _\alpha }{\ell^{2\kappa _\alpha} \vert J'_{\nu_\alpha} (j_{\nu_\alpha,m} ) \vert ^2} 
\int _{j_{\nu_\alpha,m} (\frac{a}{\ell}) ^{\kappa_\alpha}} ^{j_{\nu_\alpha,m} (\frac{b}{\ell}) ^{\kappa_\alpha}}
\Bigl( \ell (\frac{y}{j_{\nu_\alpha,m}}) ^{1/\kappa_\alpha} \Bigr) ^{1-\alpha} J_{\nu _\alpha} (y) ^2
\frac{\ell}{j_{\nu_\alpha,m}^{1/\kappa_\alpha}} \frac{1}{\kappa_\alpha} y^{\frac{1}{\kappa_\alpha} -1} \, dy
\\
= \frac{2 \kappa _\alpha \ell ^{2-\alpha} }{\kappa_\alpha \ell^{2\kappa _\alpha} 
j_{\nu_\alpha,m}^{\frac{1-\alpha}{\kappa_\alpha} + \frac{1}{\kappa_\alpha}}
\vert J'_{\nu_\alpha} (j_{\nu_\alpha,m} ) \vert ^2} 
\int _{j_{\nu_\alpha,m} (\frac{a}{\ell}) ^{\kappa_\alpha}} ^{j_{\nu_\alpha,m} (\frac{b}{\ell}) ^{\kappa_\alpha}}
y^{\frac{1-\alpha}{\kappa_\alpha} + \frac{1}{\kappa_\alpha} -1} J_{\nu _\alpha} (y) ^2 \, dy
\\
=  \frac{2}{  j_{\nu_\alpha,m}^2 \vert J'_{\nu_\alpha} (j_{\nu_\alpha,m} ) \vert ^2} 
\int _{j_{\nu_\alpha,m} (\frac{a}{\ell}) ^{\kappa_\alpha}} ^{j_{\nu_\alpha,m} (\frac{b}{\ell}) ^{\kappa_\alpha}}
y J_{\nu _\alpha} (y) ^2 \, dy ,
\end{multline*}
where we used the change of variables
$$ y= j_{\nu_\alpha,m} (\frac{x}{\ell}) ^{\kappa_\alpha}, \quad x = \ell (\frac{y}{j_{\nu_\alpha,m}}) ^{1/\kappa_\alpha},
\quad dx = \frac{\ell}{j_{\nu_\alpha,m}^{1/\kappa_\alpha}} \frac{1}{\kappa_\alpha} y^{\frac{1}{\kappa_\alpha} -1} \, dy .$$
Now introduce the function
$$ K_{\nu_\alpha,m} (y) := \sqrt{y} \frac{J_{\nu _\alpha} (y) }{J'_{\nu_\alpha} (j_{\nu_\alpha,m} )} .$$
With the help of $K_{\nu_\alpha,m}$, we have
\begin{equation}
\label{norm-K}
\int_a^b \Phi_{\alpha ,m}(x) ^2 \, dx 
= \frac{2}{  j_{\nu_\alpha,m}^2 } 
\int _{j_{\nu_\alpha,m} (\frac{a}{\ell}) ^{\kappa_\alpha}} ^{j_{\nu_\alpha,m} (\frac{b}{\ell}) ^{\kappa_\alpha}}
K_{\nu_\alpha,m} (y) ^2 \, dy .
\end{equation}
Moreover, it is well known that $K_{\nu_\alpha,m} $ is solution of the second order ordinary differential equation
$$ K_{\nu_\alpha,m} '' + h_{\nu_\alpha}(y) K_{\nu_\alpha,m} = 0 ,$$
where 
$$ h_{\nu_\alpha}(y) = 1 - \frac{\nu _\alpha ^2 - \frac{1}{4}}{y^2} .$$
(We already recalled this in the proof of Lemma \ref{lem-Sturm}.)
Hence in fact $K_{\nu_\alpha,m} $ solves the Cauchy problem
\begin{equation}
\label{cauchy-K}
\begin{cases}
K_{\nu_\alpha,m} '' + h_{\nu_\alpha}(y) K_{\nu_\alpha,m} = 0 , \\
K_{\nu_\alpha,m} (j_{\nu_\alpha,m}) = 0 , \\
K'_{\nu_\alpha,m} (j_{\nu_\alpha,m}) = \sqrt{j_{\nu_\alpha,m}} ,
\end{cases}
\end{equation}
and we want to estimate it on the zone $[j_{\nu_\alpha,m} (\frac{a}{\ell}) ^{\kappa_\alpha} , j_{\nu_\alpha,m} (\frac{b}{\ell}) ^{\kappa_\alpha}]$. Let us normalize the localization of the Cauchy conditions and the localization of the integration interval, using a suitable change of variables: consider
\begin{equation}
L_{\nu_\alpha,m} (z) = \frac{-1}{\sqrt{j_{\nu_\alpha,m}}} K_{\nu_\alpha,m} (j_{\nu_\alpha,m} - z j_{\nu_\alpha,m}) .
\end{equation}
Then $L_{\nu_\alpha,m}$ solves the Cauchy problem
\begin{equation}
\label{cauchy-L}
\begin{cases}
L_{\nu_\alpha,m} '' + k_{\nu_\alpha,m}(z) L_{\nu_\alpha,m} = 0 , \quad \text{ with } \quad k_{\nu_\alpha,m}(z) =j_{\nu_\alpha,m} ^2 - \frac{\nu _\alpha ^2 - \frac{1}{4}}{(1-z)^2} , \\
L_{\nu_\alpha,m} (0) = 0 , \\
L'_{\nu_\alpha,m} (0) = j_{\nu_\alpha,m} ,
\end{cases}
\end{equation}
and
\begin{equation}
\label{norm-L}
\int_a^b \Phi_{\alpha ,m}(x) ^2 \, dx  = 2
\int _{1-(\frac{b}{\ell}) ^{\kappa_\alpha}} ^{1-(\frac{a}{\ell}) ^{\kappa_\alpha}}
L_{\nu_\alpha,m} (z) ^2 \, dz .
\end{equation}
Once again, the term we are interested in is the norm of the solution of a Cauchy problem, but now with Cauchy conditions at the point $0$, and we have to estimate its norm on some fixed interval (that does not contain $0$).
To do this, we are going to study the Cauchy problem \eqref{cauchy-L}.


\subsection{The study of the Cauchy problem: a uniform bound on $L_{\nu_\alpha,m}$} \hfill

We begin by the following observation: when $\alpha \to 2^-$, then 
$$ 1-(\frac{b}{\ell}) ^{\kappa_\alpha} 
= 1- e^{(\frac{1}{2} \ln \frac{b}{\ell} ) (2-\alpha) } 
= -(\frac{1}{2} \ln \frac{b}{\ell} ) (2-\alpha) + O((2-\alpha)^2) ,$$
hence
$$ 1-(\frac{b}{\ell}) ^{\kappa_\alpha}
\sim \kappa_\alpha \ln \frac{\ell}{b} \quad \text{ as } \alpha \to 2^-,$$
and in the same way
$$
1-(\frac{a}{\ell}) ^{\kappa_\alpha} \sim \kappa_\alpha \ln \frac{\ell}{a} \quad \text{ as } \alpha \to 2^- ;$$
hence the integration interval shrinks to $0$, and its length satisfies
$$ \Bigl( 1-(\frac{a}{\ell}) ^{\kappa_\alpha} \Bigr) - \Bigl( 1-(\frac{b}{\ell}) ^{\kappa_\alpha} \Bigr)
\sim \kappa_\alpha \ln \frac{b}{a} \quad \text{ as } \alpha \to 2^- .$$
In particular, there exists some $0 <\gamma _* (a,b,\ell) \leq \gamma ^* (a,b,\ell)$ and $\underline \gamma (a,b,\ell) >0$ such that, for all $\alpha \in [1,2)$, 
\begin{gather}
\label{bounds-bornes}
\gamma _* (a,b,\ell) \kappa_\alpha \leq 1-(\frac{b}{\ell}) ^{\kappa_\alpha} < 1-(\frac{a}{\ell}) ^{\kappa_\alpha} \leq \gamma ^* (a,b,\ell) \kappa_\alpha , \\
\label{bounds-bornes2}
(1-(\frac{a}{\ell}) ^{\kappa_\alpha} ) - (1-(\frac{b}{\ell}) ^{\kappa_\alpha}) \geq \underline \gamma (a,b,\ell) \kappa_\alpha .
\end{gather}

Let us prove the following uniform bound:

\begin{Lemma}
\label{lem-uniformbound}
There exists $C_u$ independent of $\alpha \in [1,2)$ and of $m\geq 1$ such that
\begin{equation}
\label{bound-L}
\forall \alpha \in [1,2[, \forall m\geq 1, \forall z \in [0, 1-(\frac{a}{\ell}) ^{\kappa_\alpha}], \quad \vert  L_{\nu_\alpha,m} (z) \vert  \leq C_u .
\end{equation}
\end{Lemma}

\noindent {\it Proof of Lemma \ref{lem-uniformbound}.} 
To obtain an integral equation satisfied by $L_{\nu_\alpha,m}$, we write the Cauchy problem \eqref{cauchy-L} under the form
\begin{equation}
\label{cauchy-L2}
\begin{cases}
L_{\nu_\alpha,m} '' +   j_{\nu_\alpha,m} ^2  L_{\nu_\alpha,m} = 
\frac{\nu _\alpha ^2 - \frac{1}{4}}{(1-z)^2} L_{\nu_\alpha,m} , \\
L_{\nu_\alpha,m} (0) = 0 , \\
L'_{\nu_\alpha,m} (0) = j_{\nu_\alpha,m} .
\end{cases}
\end{equation}
Since the solution of the Cauchy problem
\begin{equation*}
\begin{cases}
Y '' +   \omega ^2  Y  = g(z), \\
Y (0) = 0 , \\
Y' (0) = \omega 
\end{cases}
\end{equation*}
is
$$ Y(z) = \sin (\omega z) + \frac{1}{\omega} \int _0 ^z g(s) \sin (\omega (z-s)) \, ds ,$$
we deduce that $L_{\nu_\alpha,m}$ satisfies
\begin{equation}
\label{ei1}
L_{\nu_\alpha,m} (z) = \sin ( j_{\nu_\alpha,m} z )
+ \frac{1}{j_{\nu_\alpha,m}} \int _0 ^z \frac{\nu _\alpha ^2 - \frac{1}{4}}{(1-s)^2} L_{\nu_\alpha,m} (s) \sin  (j_{\nu_\alpha,m} (z-s) ) \, ds .
\end{equation}
Hence
$$ \vert L_{\nu_\alpha,m} (z) \vert \leq 1 + \frac{1}{j_{\nu_\alpha,m}} \int _0 ^z \frac{\vert \nu _\alpha ^2 - \frac{1}{4}\vert }{(1-s)^2} \vert L_{\nu_\alpha,m} (s) \vert  \, ds ,$$
and the classical Gronwall inequality gives that
$$ \vert L_{\nu_\alpha,m} (z) \vert  \leq e^{\frac{\vert \nu _\alpha ^2 - \frac{1}{4}\vert }{j_{\nu_\alpha,m}} \int _0 ^z \frac{1}{(1-s)^2} \, ds}
= e^{\frac{\vert \nu _\alpha ^2 - \frac{1}{4}\vert }{j_{\nu_\alpha,m}} \frac{z}{1-z} } . $$
But then, since we know from \eqref{eq-Qu-Wong} that $j_{\nu_\alpha,m} \geq \nu_\alpha$, we have
$$ \forall \alpha \in [1,2[, \forall m\geq 1, \forall z \in [0, 1-(\frac{a}{\ell}) ^{\kappa_\alpha}], \quad 
\vert L_{\nu_\alpha,m} (z) \vert  \leq e^{\frac{\vert \nu _\alpha ^2 - \frac{1}{4}\vert }{\nu_\alpha} \frac{1-(\frac{a}{\ell}) ^{\kappa_\alpha}}{(\frac{a}{\ell}) ^{\kappa_\alpha}} } .$$
Using \eqref{bounds-bornes}, we see that this is uniformly bounded with respect to $\alpha$, hence we obtain \eqref{bound-L}. \qed


\subsection{The $L^2$ norm of $L_{\nu_\alpha,m}$ for fixed values of $\nu _\alpha$} \hfill

The integral expression \eqref{ei1} and the uniform bound \eqref{bound-L} allow us to prove the following 

\begin{Lemma}
\label{lem-minLL2-0}
There exists $\overline \gamma = \overline \gamma (a,b,\ell)$ such that
\begin{equation}
\label{minLL2-0}
\forall \alpha \in [1,2), \forall m\geq 1, \quad \int _{1-(\frac{b}{\ell}) ^{\kappa_\alpha}} ^{1-(\frac{a}{\ell}) ^{\kappa_\alpha}}
L_{\nu_\alpha,m} (z) ^2 \, dz  \geq \frac{1}{2} \underline \gamma \kappa _\alpha - \frac{\overline \gamma}{j_{\nu_\alpha,m} } .
\end{equation}
\end{Lemma}

\noindent {\it Proof of Lemma \ref{lem-minLL2-0}.}
From \eqref{ei1} we have
\begin{multline*}
\int _{1-(\frac{b}{\ell}) ^{\kappa_\alpha}} ^{1-(\frac{a}{\ell}) ^{\kappa_\alpha}}
L_{\nu_\alpha,m} (z) ^2 \, dz 
= \int _{1-(\frac{b}{\ell}) ^{\kappa_\alpha}} ^{1-(\frac{a}{\ell}) ^{\kappa_\alpha}} \sin ^2 ( j_{\nu_\alpha,m} z) \, dz
\\
+ \int _{1-(\frac{b}{\ell}) ^{\kappa_\alpha}} ^{1-(\frac{a}{\ell}) ^{\kappa_\alpha}} \Bigl( \frac{1}{j_{\nu_\alpha,m}} \int _0 ^z \frac{\nu _\alpha ^2 - \frac{1}{4}}{(1-s)^2} L_{\nu_\alpha,m} (s) \sin ( j_{\nu_\alpha,m} (z-s)) \, ds  \Bigr) ^2 \, dz 
\\
+ \int _{1-(\frac{b}{\ell}) ^{\kappa_\alpha}} ^{1-(\frac{a}{\ell}) ^{\kappa_\alpha}} 2 \sin  ( j_{\nu_\alpha,m} z ) \Bigl( \frac{1}{j_{\nu_\alpha,m}} \int _0 ^z \frac{\nu _\alpha ^2 - \frac{1}{4}}{(1-s)^2} L_{\nu_\alpha,m} (s) \sin  (j_{\nu_\alpha,m} (z-s)) \, ds  \Bigr)  \, dz.
\end{multline*}
Of course
$$ \int _{1-(\frac{b}{\ell}) ^{\kappa_\alpha}} ^{1-(\frac{a}{\ell}) ^{\kappa_\alpha}} \Bigl( \frac{1}{j_{\nu_\alpha,m}} \int _0 ^z \frac{\nu _\alpha ^2 - \frac{1}{4}}{(1-s)^2} L_{\nu_\alpha,m} (s) \sin ( j_{\nu_\alpha,m} (z-s)) \, ds  \Bigr) ^2 \, dz \geq 0 ,$$
and 
\begin{multline*}
\int _{1-(\frac{b}{\ell}) ^{\kappa_\alpha}} ^{1-(\frac{a}{\ell}) ^{\kappa_\alpha}} \sin ^2  (j_{\nu_\alpha,m} z) \, dz
= \int _{1-(\frac{b}{\ell}) ^{\kappa_\alpha}} ^{1-(\frac{a}{\ell}) ^{\kappa_\alpha}} \frac{1- \cos (2 j_{\nu_\alpha,m} z)}{2} \, dz
\\
= \frac{1}{2}\Bigl( (1-(\frac{a}{\ell}) ^{\kappa_\alpha}) - (1-(\frac{b}{\ell}) ^{\kappa_\alpha}) \Bigr)
-\frac{1}{2} [\frac{\sin 2 j_{\nu_\alpha,m} z}{2 j_{\nu_\alpha,m} }] _{1-(\frac{b}{\ell}) ^{\kappa_\alpha}} ^{1-(\frac{a}{\ell}) ^{\kappa_\alpha}} 
\\
\geq \frac{1}{2}\underline \gamma  \kappa _\alpha - \frac{1}{2 j_{\nu_\alpha,m} } .
\end{multline*}
And using \eqref{bound-L} and \eqref{bounds-bornes}, we have
\begin{multline*}
\Bigl\vert \int _{1-(\frac{b}{\ell}) ^{\kappa_\alpha}} ^{1-(\frac{a}{\ell}) ^{\kappa_\alpha}} 2 \sin  ( j_{\nu_\alpha,m} z) \Bigl( \frac{1}{j_{\nu_\alpha,m}} \int _0 ^z \frac{\nu _\alpha ^2 - \frac{1}{4}}{(1-s)^2} L_{\nu_\alpha,m} (s) \sin  (j_{\nu_\alpha,m} (z-s)) \, ds  \Bigr)  \, dz \Bigr\vert 
\\
\leq \frac{2}{j_{\nu_\alpha,m}} \int _{1-(\frac{b}{\ell}) ^{\kappa_\alpha}} ^{1-(\frac{a}{\ell}) ^{\kappa_\alpha}} \Bigl(  \int _0 ^z \frac{\vert \nu _\alpha ^2 - \frac{1}{4}\vert }{(1-s)^2} C_u \, ds  \Bigr)  \, dz
= 2 C_u \frac{\vert \nu _\alpha ^2 - \frac{1}{4}\vert}{j_{\nu_\alpha,m}} 
\int _{1-(\frac{b}{\ell}) ^{\kappa_\alpha}} ^{1-(\frac{a}{\ell}) ^{\kappa_\alpha}} \frac{z}{1-z} \, dz 
\\
\leq \frac{2 C_u }{(\frac{a}{\ell}) ^{\kappa_\alpha}}\frac{\vert \nu _\alpha ^2 - \frac{1}{4}\vert}{j_{\nu_\alpha,m}} 
[\frac{z^2}{2}]_{1-(\frac{b}{\ell}) ^{\kappa_\alpha}} ^{1-(\frac{a}{\ell}) ^{\kappa_\alpha}}
\\
= \frac{ C_u }{(\frac{a}{\ell}) ^{\kappa_\alpha}}\frac{\vert \nu _\alpha ^2 - \frac{1}{4}\vert}{j_{\nu_\alpha,m}} 
\Bigl( (1-(\frac{a}{\ell}) ^{\kappa_\alpha}) - (1-(\frac{b}{\ell}) ^{\kappa_\alpha}) \Bigr)
\Bigl( (1-(\frac{a}{\ell}) ^{\kappa_\alpha}) + (1-(\frac{b}{\ell}) ^{\kappa_\alpha}) \Bigr)
\\
\leq \frac{ C_u }{(\frac{a}{\ell}) ^{\kappa_\alpha}}\frac{\vert \nu _\alpha ^2 - \frac{1}{4}\vert}{j_{\nu_\alpha,m}} 
2 \gamma ^* (\gamma ^* - \gamma _*) \kappa _\alpha ^2 .
\end{multline*}
Hence
$$
\int _{1-(\frac{b}{\ell}) ^{\kappa_\alpha}} ^{1-(\frac{a}{\ell}) ^{\kappa_\alpha}}
L_{\nu_\alpha,m} (z) ^2 \, dz 
\geq \frac{1}{2} \underline \gamma \kappa _\alpha - \frac{1}{2j_{\nu_\alpha,m} }
- \frac{ C_u }{(\frac{a}{\ell}) ^{\kappa_\alpha}}\frac{\vert \nu _\alpha ^2 - \frac{1}{4}\vert}{j_{\nu_\alpha,m}} 
2 \gamma ^*  (\gamma ^* - \gamma _*) \kappa _\alpha ^2 .$$
Hence there exists $\overline \gamma = \overline \gamma (a,b,\ell)$ such that
$$ \int _{1-(\frac{b}{\ell}) ^{\kappa_\alpha}} ^{1-(\frac{a}{\ell}) ^{\kappa_\alpha}}
L_{\nu_\alpha,m} (z) ^2 \, dz  \geq \frac{1}{2} \underline \gamma \kappa _\alpha - \frac{\overline \gamma}{j_{\nu_\alpha,m} } .$$
This implies \eqref{minLL2-0} . \qed


\subsection{A first consequence for the eigenfunctions on the control region } \hfill

The consequence of \eqref{minLL2-0} is immediate: combining \eqref{norm-L} and  \eqref{minLL2-0}, we obtain

\begin{equation}
\label{norm-L-phi}
\int_a^b \Phi_{\alpha ,m}(x) ^2 \, dx  \geq 
\underline \gamma \kappa _\alpha - \frac{2\overline \gamma}{j_{\nu_\alpha,m} },
\end{equation}
hence we find again what is well-known at least in the nondegenerate case (\cite{Lagnese}):
$$ \lim _{m\to \infty} \int_a^b \Phi_{\alpha ,m}(x) ^2 \, dx   >0,$$
hence the sequence of positive terms $(\int_a^b \Phi_{\alpha ,m} ^2)_{m\geq 1}$ is bounded from below by a positive constant. But this constant may depend on $\alpha$.
At least, we obtain the following uniform result:
given $\overline \alpha ^* \in [1,2)$, there exists $\overline{m}^* \geq 1$ such that
$$ \forall \alpha \in [1,\overline \alpha ^*], \forall m\geq \overline{m}^*, \quad  \int_a^b \Phi_{\alpha ,m}(x) ^2 \, dx
\geq \frac{1}{2} \underline \gamma \kappa _{\overline \alpha ^*} .$$
Hence, since $\Phi_{\alpha ,m}$ depends continuously on the parameter $\alpha$, we obtain that, given $\overline \alpha ^* \in [1,2)$, the sequences
$(\int_a^b \Phi_{\alpha ,m} ^2)_{m\geq 1}$ are bounded from below by a positive constant, uniformly with respect to $\alpha \in [1, \overline \alpha ^*]$: 
\begin{equation}
\label{pdtVic}
\forall \overline \alpha ^* \in [1,2), \exists \underline \gamma ^* >0, \forall \alpha \in [1, \overline \alpha ^*], \forall m\geq 1, \quad
\int_a^b \Phi_{\alpha ,m}(x) ^2 \, dx  \geq \underline \gamma ^* \kappa _{\overline \alpha ^*} .
\end{equation}
This is not sufficient to conclude, since we want a lower estimate valid for all $\alpha \in [1,2)$, but this is a first step, and we will use this partial result later.


\subsection{Another integral equation for $L_{\nu_\alpha,m}$ when $\nu_\alpha$ is large } \hfill

Now, we would like to obtain bounds from below when $\nu_\alpha$ is large. In this case, we have to  integrate $L _{\nu_\alpha,m}^2$ in an interval close to $0$.
So, since we are interested in looking what happens near $0$, it is more interesting to write
\begin{multline*}
k_{\nu_\alpha,m}(z) =j_{\nu_\alpha,m} ^2 - \frac{\nu _\alpha ^2 - \frac{1}{4}}{(1-z)^2}
= j_{\nu_\alpha,m} ^2 - (\nu _\alpha ^2 - \frac{1}{4}) \Bigl( 1 + \frac{1}{(1-z)^2} -1 \Bigr) 
\\
= \Bigl( j_{\nu_\alpha,m} ^2 - \nu _\alpha ^2 + \frac{1}{4} \Bigr) - (\nu _\alpha ^2 - \frac{1}{4}) \frac{2z-z^2}{(1-z)^2} .
\end{multline*}
Then we can write the Cauchy problem \eqref{cauchy-L} under the form
\begin{equation}
\label{cauchy-L3}
\begin{cases}
L_{\nu_\alpha,m} '' +  \Bigl( j_{\nu_\alpha,m} ^2 - \nu _\alpha ^2 + \frac{1}{4} \Bigr) L_{\nu_\alpha,m} = 
(\nu _\alpha ^2 - \frac{1}{4}) \frac{2z-z^2}{(1-z)^2} L_{\nu_\alpha,m} , \\
L_{\nu_\alpha,m} (0) = 0 , \\
L'_{\nu_\alpha,m} (0) = j_{\nu_\alpha,m} .
\end{cases}
\end{equation}
Since the solution of the Cauchy problem
\begin{equation*}
\begin{cases}
Y '' +   \omega ^2  Y  = g(z), \\
Y (0) = 0 , \\
Y' (0) = \rho 
\end{cases}
\end{equation*}
is
$$ Y(z) = \frac{\rho}{\omega} \sin \omega z + \frac{1}{\omega} \int _0 ^z g(s) \sin \omega (z-s) \, ds ,$$
we obtain a new integral equation satisfied by $L_{\nu_\alpha,m}$: denote
$$ \omega _{\nu_\alpha,m} := \sqrt{j_{\nu_\alpha,m} ^2 - \nu _\alpha ^2 + \frac{1}{4}};$$
then we have 
\begin{multline*}
L_{\nu_\alpha,m} (z) = \frac{j_{\nu_\alpha,m}}{\omega _{\nu_\alpha,m}} \sin (\omega _{\nu_\alpha,m} z )
\\
+ \frac{1}{\omega _{\nu_\alpha,m}} \int _0 ^z (\nu _\alpha ^2 - \frac{1}{4}) \frac{2s-s^2}{(1-s)^2} L_{\nu_\alpha,m}(s) \sin (\omega _{\nu_\alpha,m} (z-s)) \, ds .
\end{multline*}
Hence
\begin{multline}
\label{LL2}
\int _{1-(\frac{b}{\ell}) ^{\kappa_\alpha}} ^{1-(\frac{a}{\ell}) ^{\kappa_\alpha}}
L_{\nu_\alpha,m} (z) ^2 \, dz 
= \int _{1-(\frac{b}{\ell}) ^{\kappa_\alpha}} ^{1-(\frac{a}{\ell}) ^{\kappa_\alpha}}
\Bigl( \frac{j_{\nu_\alpha,m}}{\omega _{\nu_\alpha,m}} \sin (\omega _{\nu_\alpha,m} z )
\\
+ \frac{1}{\omega _{\nu_\alpha,m}} \int _0 ^z (\nu _\alpha ^2 - \frac{1}{4}) \frac{2s-s^2}{(1-s)^2} L_{\nu_\alpha,m}(s) \sin (\omega _{\nu_\alpha,m} (z-s) )\, ds \Bigr) ^2 \, dz
\\
= \frac{j_{\nu_\alpha,m}^2}{\omega _{\nu_\alpha,m}^2} \int _{1-(\frac{b}{\ell}) ^{\kappa_\alpha}} ^{1-(\frac{a}{\ell}) ^{\kappa_\alpha}} \sin ^2 (\omega _{\nu_\alpha,m} z) \, dz
\\
+ \frac{1}{\omega _{\nu_\alpha,m} ^2} \int _{1-(\frac{b}{\ell}) ^{\kappa_\alpha}} ^{1-(\frac{a}{\ell}) ^{\kappa_\alpha}}
\Bigl( \int _0 ^z (\nu _\alpha ^2 - \frac{1}{4}) \frac{2s-s^2}{(1-s)^2} L_{\nu_\alpha,m}(s) \sin (\omega _{\nu_\alpha,m} (z-s)) \, ds \Bigr) ^2 \, dz
\\
+ 2 \frac{j_{\nu_\alpha,m}}{\omega _{\nu_\alpha,m}^2}\int _{1-(\frac{b}{\ell}) ^{\kappa_\alpha}} ^{1-(\frac{a}{\ell}) ^{\kappa_\alpha}}  \sin (\omega _{\nu_\alpha,m} z) \Bigl( \int _0 ^z (\nu _\alpha ^2 - \frac{1}{4}) \frac{2s-s^2}{(1-s)^2} L_{\nu_\alpha,m}(s) \sin ( \omega _{\nu_\alpha,m} (z-s)) \, ds \Bigr) \, dz .
\end{multline}
We are going to study the behavior of the first and third term of the right hand side of \eqref{LL2}, the second one being nonnegative.


\subsection{The $L^2$ norm of $L_{\nu_\alpha,m}$ for large values of $\nu _\alpha$} \hfill


\subsubsection{The first term of \eqref{LL2}} \hfill

We study 
\begin{equation}
\label{20nov-0}
\int _{1-(\frac{b}{\ell}) ^{\kappa_\alpha}} ^{1-(\frac{a}{\ell}) ^{\kappa_\alpha}} \sin ^2 \omega _{\nu_\alpha,m} z \, dz .
\end{equation}
It appears that we need to distinguish the cases $\omega _{\nu_\alpha,m} \kappa _\alpha$ small and
$\omega _{\nu_\alpha,m} \kappa _\alpha$ not small:
indeed, 
$$ \int _{1-(\frac{b}{\ell}) ^{\kappa_\alpha}} ^{1-(\frac{a}{\ell}) ^{\kappa_\alpha}} \sin ^2 (\omega _{\nu_\alpha,m} z) \, dz
= \frac{1}{\omega _{\nu_\alpha,m}} \int _{\omega _{\nu_\alpha,m}[1-(\frac{b}{\ell}) ^{\kappa_\alpha}]} ^{\omega _{\nu_\alpha,m}[1-(\frac{a}{\ell}) ^{\kappa_\alpha}]} \sin ^2 x \, dx ,
$$
and we derive from the elementary convexity inequalities:
\begin{equation}
\label{20nov-1}
\forall \mu \geq 0, \forall u \in [0,1], \quad (1 -e^{-\mu}) u \leq 1- e^{-\mu u} \leq \mu u ,
\end{equation}
that
\begin{itemize}
\item first, using \eqref{20nov-1} with $u= \kappa _\alpha$ and $\mu = \ln \frac{\ell}{a}$, we have
\begin{equation}
\label{20nov-2}
[1-(\frac{a}{\ell}) ^{\kappa_\alpha}]
=  [1- e^{- \kappa _\alpha \ln \frac{\ell}{a}}]
\leq \ln \frac{\ell}{a} \kappa _\alpha  , 
\end{equation}
\item next, using \eqref{20nov-1} with $u= \kappa _\alpha$ and $\mu = \ln \frac{\ell}{b}$, we have
\begin{equation}
\label{20nov-3}
[1-(\frac{b}{\ell}) ^{\kappa_\alpha}]
=  [1- e^{- \kappa _\alpha \ln \frac{\ell}{b}}]
\geq  (1- e^{-\ln \frac{\ell}{b}}) \kappa _\alpha
= (1- \frac{b}{\ell}) \kappa _\alpha , 
\end{equation}
\end{itemize}
and \eqref{20nov-2}-\eqref{20nov-3} give that
$$ (1- \frac{b}{\ell})\,  \omega _{\nu_\alpha,m} \kappa _\alpha
\leq \omega _{\nu_\alpha,m}[1-(\frac{b}{\ell}) ^{\kappa_\alpha}] 
\leq \omega _{\nu_\alpha,m}[1-(\frac{a}{\ell}) ^{\kappa_\alpha}] 
\leq (\ln \frac{\ell}{a}) \,  \omega _{\nu_\alpha,m} \kappa _\alpha .$$
Hence the bounds of the integral appearing in \eqref{20nov-0} are both small or both non small,
depending on the value of $\omega _{\nu_\alpha,m} \kappa _\alpha$.

We prove the following
\begin{Lemma}
\label{lem-LL2-1}
There exists $\eta _0 = \eta _0 (a,b,\ell)>0$, and $\gamma _0 = \gamma _0 (a,b,\ell)>0$ both independent of $\alpha \in [1,2)$ and of $m\geq 1$ such that
\begin{itemize}
\item if $\omega _{\nu_\alpha,m} \kappa _\alpha \leq \eta _0$, then
\begin{equation}
\label{LL2-1-small}
\int _{1-(\frac{b}{\ell}) ^{\kappa_\alpha}} ^{1-(\frac{a}{\ell}) ^{\kappa_\alpha}} \sin ^2 (\omega _{\nu_\alpha,m} z) \, dz
\geq \gamma _0  \, \omega _{\nu_\alpha,m} ^2 \, \kappa_\alpha ^3 .
\end{equation}

\item  if $\omega _{\nu_\alpha,m} \kappa _\alpha \geq \eta _0$, then
\begin{equation}
\label{LL2-1-large}
\int _{1-(\frac{b}{\ell}) ^{\kappa_\alpha}} ^{1-(\frac{a}{\ell}) ^{\kappa_\alpha}} \sin ^2 (\omega _{\nu_\alpha,m} z) \, dz
\geq \gamma _0 \kappa _\alpha .
\end{equation}
\end{itemize}
\end{Lemma}

\begin{Remark} \em
A similar property of the function sinus appears in Haraux \cite{Haraux} (Lemma 1.3.2) and \cite{Privat} (Theorem 1).
In our case, we have to bound from below the integrals of $z\mapsto \sin ^2 (\omega _{\nu_\alpha,m} z)$ with respect to the size of the integration zone (and this size is small, of the order $\kappa _\alpha$), and the coefficients $\omega _{\nu_\alpha,m}$ that appear are non integer and possibly small.

\end{Remark}

\noindent {\it Proof of Lemma \ref{lem-LL2-1}.} It comes from the following observations: 
\begin{itemize}
\item First,
\begin{equation}
\label{obsABsmall}
0<A<B\leq \frac{\pi}{2} \quad \implies \quad \int _A ^B \sin ^2 x \, dx \geq \frac{4}{\pi^2} A^2 (B-A) .
\end{equation}
Indeed, if $0<A<B\leq \frac{\pi}{2}$, then
$$ \forall x \in [A,B], \quad \sin  x  \geq \frac{2}{\pi} x , $$
hence
$$ \int _A ^B \sin ^2 x \, dx \geq \frac{4}{\pi^2} A^2 (B-A) .$$
\item On the other hand, 
consider $\eta _1$, $A$ and $B$ such that
$$ 0 < \eta _1 \leq A < B, \quad \text{ and } \quad B-A \geq 2 \eta _1 .$$
Then, first the function 
$$ s \mapsto \int _s ^{s+\eta _1} \sin ^2 x \, dx $$
is continuous, $2\pi$-periodic, and positive, hence is bounded from below by a positive constant, denoted $\gamma _1$.
Next, there exists one and only one integer $k$ such that
$$ k\eta _1 \leq B-A < (k+1) \eta _1 ,$$
then
\begin{multline*}
\int _A ^B \sin ^2 x \, dx \geq \int _A ^{A+k\eta_1} \sin ^2 x \, dx
= \sum _{j=0} ^{k-1} \int _{A+j\eta _1} ^{A+(j+1)\eta _1} \sin ^2 x \, dx
\\
\geq \sum _{j=0} ^{k-1} \gamma _1 
= k \gamma _1 \geq \gamma _1 (\frac{B-A}{\eta _1} - 1)
= \gamma _1 (\frac{B-A}{2\eta _1} + \frac{B-A}{2\eta _1}- 1) \geq \frac{\gamma _1 }{2\eta _1} (B-A) .
\end{multline*}
Hence
\begin{multline}
\label{obsABlarge}
\Bigl( 0 < \eta _1 \leq A < B, 
\quad 
\text{ and } 
\quad
 B-A \geq 2 \eta _1 \Bigr) 
\\ 
\implies 
\quad 
\int _A ^B \sin ^2 x \, dx \geq \frac{\gamma _1 }{2\eta _1} (B-A) 
\quad 
\text{ with } \gamma _1=\gamma _1(\eta _1)  .
\end{multline}

\end{itemize}
Now we are in position to conclude the proof of Lemma \ref{lem-LL2-1}. This is based on the observation that we are in 
one of the two situations studied previously: using \eqref{bounds-bornes}, we see that there is some $\eta _0$ such that
$$ \omega _{\nu_\alpha,m} \kappa _\alpha \leq \eta _0 \quad \implies \quad 
0 < \omega _{\nu_\alpha,m}[1-(\frac{b}{\ell}) ^{\kappa_\alpha}] < \omega _{\nu_\alpha,m}[1-(\frac{a}{\ell}) ^{\kappa_\alpha}] \leq \frac{\pi}{2} ;$$
in this case, \eqref{obsABsmall} gives that
\begin{multline*}
\int _{\omega _{\nu_\alpha,m}[1-(\frac{b}{\ell}) ^{\kappa_\alpha}]} ^{\omega _{\nu_\alpha,m}[1-(\frac{a}{\ell}) ^{\kappa_\alpha}]} \sin ^2 x \, dx
\\
\geq \frac{4}{\pi^2} \Bigl( \omega _{\nu_\alpha,m}[1-(\frac{b}{\ell}) ^{\kappa_\alpha}] \Bigr) ^2 
\Bigl(\omega _{\nu_\alpha,m}[1-(\frac{a}{\ell}) ^{\kappa_\alpha}]- \omega _{\nu_\alpha,m}[1-(\frac{b}{\ell}) ^{\kappa_\alpha}] \Bigr) 
\\
= \frac{4}{\pi^2} \omega _{\nu_\alpha,m} ^3 [1-(\frac{b}{\ell}) ^{\kappa_\alpha}]  ^2 \Bigl([1-(\frac{a}{\ell}) ^{\kappa_\alpha}]- [1-(\frac{b}{\ell}) ^{\kappa_\alpha}] \Bigr) ,
\end{multline*}
hence, thanks to \eqref{bounds-bornes} and \eqref{bounds-bornes2}, there exists some $\gamma (a,b,\ell) >0$ such that
$$
\omega _{\nu_\alpha,m} \kappa _\alpha \leq \eta _0 \quad \implies \quad 
\int _{\omega _{\nu_\alpha,m}[1-(\frac{b}{\ell}) ^{\kappa_\alpha}]} ^{\omega _{\nu_\alpha,m}[1-(\frac{a}{\ell}) ^{\kappa_\alpha}]} \sin ^2 x \, dx
\geq \gamma (a,b,\ell) \, \omega _{\nu_\alpha,m} ^3 \, \kappa_\alpha ^3 ,$$
which gives \eqref{LL2-1-small}.

Now if $\omega _{\nu_\alpha,m} \kappa _\alpha \geq \eta _0$: then, thanks to \eqref{bounds-bornes} and \eqref{bounds-bornes2}, there is some $\eta _1 = \eta _1 (a,b,\ell)$
such that
$$ \eta _1 < \omega _{\nu_\alpha,m}[1-(\frac{b}{\ell}) ^{\kappa_\alpha}] < \omega _{\nu_\alpha,m}[1-(\frac{a}{\ell}) ^{\kappa_\alpha}] ,$$
and
$$ \omega _{\nu_\alpha,m}[1-(\frac{a}{\ell}) ^{\kappa_\alpha}] - \omega _{\nu_\alpha,m}[1-(\frac{b}{\ell}) ^{\kappa_\alpha}] 
\geq 2 \eta _1 .$$
Then \eqref{obsABlarge} gives that there exists $\gamma _1 = \gamma _1 (a,b,\ell)$ such that
$$ \int _{\omega _{\nu_\alpha,m}[1-(\frac{b}{\ell}) ^{\kappa_\alpha}]} ^{\omega _{\nu_\alpha,m}[1-(\frac{a}{\ell}) ^{\kappa_\alpha}]} \sin ^2 x \, dx \geq \frac{\gamma _1 }{2\eta _1} \Bigl(\omega _{\nu_\alpha,m} [1-(\frac{a}{\ell}) ^{\kappa_\alpha}]- \omega _{\nu_\alpha,m} [1-(\frac{b}{\ell}) ^{\kappa_\alpha}] \Bigr) ,$$
hence, once again with \eqref{bounds-bornes2}, we obtain that
$$ 
\omega _{\nu_\alpha,m} \kappa _\alpha \geq \eta _0 \quad \implies \quad 
\int _{\omega _{\nu_\alpha,m}[1-(\frac{b}{\ell}) ^{\kappa_\alpha}]} ^{\omega _{\nu_\alpha,m}[1-(\frac{a}{\ell}) ^{\kappa_\alpha}]} \sin ^2 x \, dx \geq  \frac{\gamma _1 }{2\eta _1} \omega _{\nu_\alpha,m} \underline \gamma \kappa _\alpha ,$$
which gives \eqref{LL2-1-large}. \qed


\subsubsection{The third term of \eqref{LL2}} \hfill

We study
$$ M_3:= \int _{1-(\frac{b}{\ell}) ^{\kappa_\alpha}} ^{1-(\frac{a}{\ell}) ^{\kappa_\alpha}}  \sin (\omega _{\nu_\alpha,m} z ) \Bigl( \int _0 ^z (\nu _\alpha ^2 - \frac{1}{4}) \frac{2s-s^2}{(1-s)^2} L_{\nu_\alpha,m}(s) \sin (\omega _{\nu_\alpha,m} (z-s)) \, ds \Bigr) \, dz , $$
and we prove the following

\begin{Lemma}
\label{lem-LL2-3}
Choose $\eta _0 = \eta _0 (a,b,\ell)>0$ given in Lemma \ref{lem-LL2-1}. Then there exists
$\gamma ' _0 = \gamma ' _0 (a,b,\ell)>0$ independent of $\alpha \in [1,2)$ and of $m\geq 1$ such that
\begin{itemize}
\item if $\omega _{\nu_\alpha,m} \kappa _\alpha \leq \eta _0$, then
\begin{equation}
\label{LL2-3-small}
\vert M_3 \vert 
\leq \gamma _0 ' \omega _{\nu_\alpha,m} \kappa_\alpha ^2 ,
\end{equation}

\item  if $\omega _{\nu_\alpha,m} \kappa _\alpha \geq \eta _0$, then
\begin{equation}
\label{LL2-3-large}
\vert M_3 \vert 
\leq \gamma _0 ' \kappa_\alpha ,
\end{equation}
\end{itemize}
hence in any case
\begin{equation}
\label{LL2-3-unif}
\vert M_3 \vert 
\leq \gamma _0 ' \kappa_\alpha .
\end{equation}

\end{Lemma}

\noindent {\it Proof of Lemma \ref{lem-LL2-3}.} First we prove \eqref{LL2-3-small}. Assume that $\omega _{\nu_\alpha,m} \kappa _\alpha \leq \eta _0$.
Then, since we have proved in Lemma \ref{lem-uniformbound} that $L_{\nu_\alpha,m}$ is uniformly bounded in $[0, 1-(\frac{a}{\ell}) ^{\kappa_\alpha}]$, we have
\begin{multline*}
\vert M_3 \vert 
\leq \int _{1-(\frac{b}{\ell}) ^{\kappa_\alpha}} ^{1-(\frac{a}{\ell}) ^{\kappa_\alpha}}  \vert \sin \omega _{\nu_\alpha,m} z \vert \Bigl( \int _0 ^z \vert \nu _\alpha ^2 - \frac{1}{4} \vert \frac{\vert2s-s^2\vert}{(1-s)^2} \vert L_{\nu_\alpha,m}(s) \vert \vert \sin \omega _{\nu_\alpha,m} (z-s) \vert \, ds \Bigr) \, dz 
\\
\leq \int _{1-(\frac{b}{\ell}) ^{\kappa_\alpha}} ^{1-(\frac{a}{\ell}) ^{\kappa_\alpha}}  \omega _{\nu_\alpha,m} z \Bigl( \int _0 ^z \vert \nu _\alpha ^2 - \frac{1}{4}  \vert  \frac{2s}{(1-s)^2}C_u  \, ds \Bigr) \, dz  
\\ \leq  \omega _{\nu_\alpha,m} \vert \nu _\alpha ^2 - \frac{1}{4}  \vert C_u ' \int _{1-(\frac{b}{\ell}) ^{\kappa_\alpha}} ^{1-(\frac{a}{\ell}) ^{\kappa_\alpha}} z^3 \, dz 
\\
\leq  \omega _{\nu_\alpha,m} \vert \nu _\alpha ^2 - \frac{1}{4}  \vert C_u ' (1-(\frac{a}{\ell}) ^{\kappa_\alpha}) ^3
\Bigl( (1-(\frac{a}{\ell}) ^{\kappa_\alpha}) - 1-(\frac{b}{\ell}) ^{\kappa_\alpha}) \Bigr)
\\
\leq C_u '' \omega _{\nu_\alpha,m} \vert \nu _\alpha ^2 - \frac{1}{4}  \vert \kappa_\alpha ^4 
\leq C_u '' \omega _{\nu_\alpha,m} \kappa_\alpha ^{-2}  \kappa_\alpha ^4 = 
C_u '' \omega _{\nu_\alpha,m} \kappa_\alpha ^2 ,
\end{multline*}
which is \eqref{LL2-3-small} (which implies \eqref{LL2-3-unif}).

 Now we prove \eqref{LL2-3-large}. Assume that $\omega _{\nu_\alpha,m} \kappa _\alpha \geq \eta _0$. Then, once again using the fact that $L_{\nu_\alpha,m}$ is uniformly bounded in $[0, 1-(\frac{a}{\ell}) ^{\kappa_\alpha}]$ (Lemma \ref{lem-uniformbound}), we have
\begin{multline*}
\vert M_3 \vert 
\leq \int _{1-(\frac{b}{\ell}) ^{\kappa_\alpha}} ^{1-(\frac{a}{\ell}) ^{\kappa_\alpha}}  \vert \sin \omega _{\nu_\alpha,m} z \vert \Bigl( \int _0 ^z \vert \nu _\alpha ^2 - \frac{1}{4} \vert \frac{\vert2s-s^2\vert}{(1-s)^2} \vert L_{\nu_\alpha,m}(s) \vert \vert \sin \omega _{\nu_\alpha,m} (z-s) \vert \, ds \Bigr) \, dz 
\\
\leq \int _{1-(\frac{b}{\ell}) ^{\kappa_\alpha}} ^{1-(\frac{a}{\ell}) ^{\kappa_\alpha}}   \Bigl( \int _0 ^z \vert \nu _\alpha ^2 - \frac{1}{4} \vert C_u ' s  \, ds \Bigr) \, dz 
\leq C_u ' \vert \nu _\alpha ^2 - \frac{1}{4} \vert [\frac{z^3}{6}] _{1-(\frac{b}{\ell}) ^{\kappa_\alpha}} ^{1-(\frac{a}{\ell}) ^{\kappa_\alpha}} 
\\
\leq C_u '' \vert \nu _\alpha ^2 - \frac{1}{4} \vert (1-(\frac{a}{\ell}) ^{\kappa_\alpha}) ^2
\Bigl( (1-(\frac{a}{\ell}) ^{\kappa_\alpha}) - 1-(\frac{b}{\ell}) ^{\kappa_\alpha}) \Bigr)
\leq C_u ''' \kappa _\alpha ^{-2} \kappa _\alpha ^3 = C_u ''' \kappa _\alpha 
\end{multline*}
which is \eqref{LL2-3-large}. \qed


\subsubsection{The $L^2$ norm of $L_{\nu_\alpha,m}$ for large values of $\nu _\alpha$} \hfill

We prove the following
\begin{Lemma}
\label{lem-LL2-entier}
Choose $\eta _0 = \eta _0 (a,b,\ell)>0$ and $\gamma _0 = \gamma _0 (a,b,\ell)>0$ 
given in Lemma \ref{lem-LL2-1}. Then there exists
$\overline \alpha \in [1,2)$ such that
\begin{equation}
\label{LL2-entier}
\forall \alpha \in [\overline \alpha, 2), \forall m\geq 1, \quad \int _{1-(\frac{b}{\ell}) ^{\kappa_\alpha}} ^{1-(\frac{a}{\ell}) ^{\kappa_\alpha}} L_{\nu_\alpha,m} (z) ^2 \, dz 
\geq \frac{\gamma _0}{2}  \frac{ j_{\nu_\alpha,m}^2 \kappa_\alpha ^3}{1+\omega _{\nu_\alpha,m}^2  \kappa_\alpha ^2}.
\end{equation}
\end{Lemma}

\noindent {\it Proof of Lemma \ref{lem-LL2-entier}.}
We start from \eqref{LL2}, that gives
\begin{multline}
\label{LL2-2}
\int _{1-(\frac{b}{\ell}) ^{\kappa_\alpha}} ^{1-(\frac{a}{\ell}) ^{\kappa_\alpha}}
L_{\nu_\alpha,m} (z) ^2 \, dz 
\geq  \frac{j_{\nu_\alpha,m}^2}{\omega _{\nu_\alpha,m}^2} \int _{1-(\frac{b}{\ell}) ^{\kappa_\alpha}} ^{1-(\frac{a}{\ell}) ^{\kappa_\alpha}} \sin ^2 (\omega _{\nu_\alpha,m} z) \, dz
\\
+ 2 \frac{j_{\nu_\alpha,m}}{\omega _{\nu_\alpha,m}^2}\int _{1-(\frac{b}{\ell}) ^{\kappa_\alpha}} ^{1-(\frac{a}{\ell}) ^{\kappa_\alpha}}  \sin (\omega _{\nu_\alpha,m} z) \Bigl( \int _0 ^z (\nu _\alpha ^2 - \frac{1}{4}) \frac{2s-s^2}{(1-s)^2} L_{\nu_\alpha,m}(s) \sin (\omega _{\nu_\alpha,m} (z-s)) \, ds \Bigr) \, dz .
\end{multline}
First we prove \eqref{LL2-entier} when $\omega _{\nu_\alpha,m} \kappa _\alpha \leq \eta _0$.
Using \eqref{LL2-1-small} and \eqref{LL2-3-small} in \eqref{LL2-2}, we have
\begin{multline*}
\int _{1-(\frac{b}{\ell}) ^{\kappa_\alpha}} ^{1-(\frac{a}{\ell}) ^{\kappa_\alpha}}
L_{\nu_\alpha,m} (z) ^2 \, dz 
\geq  \frac{j_{\nu_\alpha,m}^2}{\omega _{\nu_\alpha,m}^2} \gamma _0  \, \omega _{\nu_\alpha,m} ^2 \, \kappa_\alpha ^3
- 2 \frac{j_{\nu_\alpha,m}}{\omega _{\nu_\alpha,m}^2} \gamma _0 ' \omega _{\nu_\alpha,m} \kappa_\alpha ^2
\\
= \gamma _0 j_{\nu_\alpha,m}^2 \kappa_\alpha ^3 -2 \gamma _0 ' \frac{j_{\nu_\alpha,m}}{\omega _{\nu_\alpha,m}} \kappa_\alpha ^2 
=  j_{\nu_\alpha,m}^2 \kappa_\alpha ^3 \Bigl( \gamma _0 - \frac{2 \gamma _0 '}{j_{\nu_\alpha,m} \omega _{\nu_\alpha,m} \kappa_\alpha} \Bigr) .
\end{multline*}
Now remember that \eqref{eq-Qu-Wong} (\cite{QuWong}) says that
$$ j_{\nu,k}   \geq  \nu - \frac{a_1}{2^{1/3}} \nu ^{1/3} ,$$
where $a_1 <0$. Hence
$$ j_{\nu_\alpha,m} - \nu _\alpha \geq - \frac{a_1}{2^{1/3}} \nu _\alpha ^{1/3} ,$$
and
$$  j_{\nu_\alpha,m} + \nu _\alpha \geq 2 \nu _\alpha , $$
therefore
$$
\omega _{\nu_\alpha,m} = \sqrt{j_{\nu_\alpha,m}^2 - \nu_\alpha ^2 + \frac{1}{4}}
\geq \sqrt{- \frac{2a_1}{2^{1/3}} \nu _\alpha ^{4/3}} .
$$
Since $j_{\nu_\alpha,m} \geq \nu_\alpha$, this implies that
$$ j_{\nu_\alpha,m} \omega _{\nu_\alpha,m} \kappa_\alpha \geq \Bigl( - \frac{2a_1}{2^{1/3}}\Bigr) ^{1/2} \nu _\alpha ^{5/3} \kappa_\alpha .$$
That last quantity goes to infinity when $\alpha \to 2^-$, hence there exists $\alpha _0 \in [1,2)$ such that
$$\forall \alpha \in [\alpha _0,2), \forall m\geq 1, \quad \int _{1-(\frac{b}{\ell}) ^{\kappa_\alpha}} ^{1-(\frac{a}{\ell}) ^{\kappa_\alpha}}
L_{\nu_\alpha,m} (z) ^2 \, dz \geq \frac{\gamma _0}{2} j_{\nu_\alpha,m}^2 \kappa_\alpha ^3 ,$$
which implies \eqref{LL2-entier} when $\omega _{\nu_\alpha,m} \kappa _\alpha \leq \eta _0$.

Next we prove \eqref{LL2-entier} when $\omega _{\nu_\alpha,m} \kappa _\alpha \geq \eta _0$.
Using \eqref{LL2-1-large} and \eqref{LL2-3-large} in \eqref{LL2-2}, we have

$$
\int _{1-(\frac{b}{\ell}) ^{\kappa_\alpha}} ^{1-(\frac{a}{\ell}) ^{\kappa_\alpha}}
L_{\nu_\alpha,m} (z) ^2 \, dz 
\geq  \frac{j_{\nu_\alpha,m}^2}{\omega _{\nu_\alpha,m}^2} \gamma _0 \kappa _\alpha
- 2 \frac{j_{\nu_\alpha,m}}{\omega _{\nu_\alpha,m}^2}\gamma _0 ' \kappa_\alpha  
=\frac{j_{\nu_\alpha,m}^2 \kappa _\alpha }{\omega _{\nu_\alpha,m}^2} \Bigl( \gamma _0 -2 \frac{\gamma ' _0}{j_{\nu_\alpha,m}} \Bigr) .
$$
Hence, once again, there exists $\alpha _1 \in [1,2)$ such that
$$\forall \alpha \in [\alpha _1,2), \forall m\geq 1, \quad \int _{1-(\frac{b}{\ell}) ^{\kappa_\alpha}} ^{1-(\frac{a}{\ell}) ^{\kappa_\alpha}}
L_{\nu_\alpha,m} (z) ^2 \, dz \geq \frac{\gamma _0 }{2} \frac{j_{\nu_\alpha,m}^2 \kappa _\alpha }{\omega _{\nu_\alpha,m}^2} .$$
Finally, since
$$ \inf \{a, \frac{1}{b} \} \geq \frac{a}{1+ab},$$
we obtain that
\begin{multline*}
\int _{1-(\frac{b}{\ell}) ^{\kappa_\alpha}} ^{1-(\frac{a}{\ell}) ^{\kappa_\alpha}}
L_{\nu_\alpha,m} (z) ^2 \, dz \geq \inf \{ \frac{\gamma _0}{2} j_{\nu_\alpha,m}^2 \kappa_\alpha ^3, \frac{\gamma _0 }{2} \frac{j_{\nu_\alpha,m}^2 \kappa _\alpha }{\omega _{\nu_\alpha,m}^2}\}
\\
= \frac{\gamma _0}{2} j_{\nu_\alpha,m}^2 \kappa_\alpha ^3 \inf \{ 1, \frac{1}{\omega _{\nu_\alpha,m}^2  \kappa_\alpha ^2}\} 
\geq \frac{\gamma _0}{2} j_{\nu_\alpha,m}^2 \kappa_\alpha ^3  \frac{1}{1+\omega _{\nu_\alpha,m}^2  \kappa_\alpha ^2}, 
\end{multline*}
which is \eqref{LL2-entier}. \qed


\subsection{Proof of Proposition \ref{prop-fctpr}} \hfill

We deduce from \eqref{norm-L} and \eqref{LL2-entier} that, for $\alpha \geq \overline \alpha$,
\begin{multline*}
\forall m \geq 1, \quad  \int _a ^b \Phi _{\alpha ,m} (x) ^2 \, dx 
\geq \gamma _0 \frac{ j_{\nu_\alpha,m}^2 \kappa_\alpha ^3}{1+\omega _{\nu_\alpha,m}^2  \kappa_\alpha ^2} 
= \gamma _0 \kappa_\alpha \frac{ j_{\nu_\alpha,m}^2 \kappa_\alpha ^2}{1+\omega _{\nu_\alpha,m}^2  \kappa_\alpha ^2} 
\\
= \gamma _0 \kappa_\alpha \frac{ \omega_{\nu_\alpha,m}^2 \kappa_\alpha ^2 + (\nu_\alpha ^2 - \frac{1}{4}) \kappa _\alpha ^2}{1+\omega _{\nu_\alpha,m}^2  \kappa_\alpha ^2} 
.
\end{multline*}
Since $\nu _\alpha \kappa _\alpha \to \frac{1}{2}$ as $\alpha \to 2^-$, there exists
$\overline \alpha ^* \in [\overline \alpha, 2)$ such that $(\nu_\alpha ^2 - \frac{1}{4}) \kappa _\alpha ^2  \geq \frac{1}{8}$ for all $\alpha \in [\overline \alpha ^*,2)$.
Then, for all $\alpha \in [\overline \alpha ^*,2)$, we have
$$ \frac{ \omega_{\nu_\alpha,m}^2 \kappa_\alpha ^2 + (\nu_\alpha ^2 - \frac{1}{4}) \kappa _\alpha ^2}{1+\omega _{\nu_\alpha,m}^2  \kappa_\alpha ^2} 
\geq \frac{ \omega_{\nu_\alpha,m}^2 \kappa_\alpha ^2 + \frac{1}{8}}{1+\omega _{\nu_\alpha,m}^2  \kappa_\alpha ^2} \geq \frac{1}{8}. $$ 
Hence 
\begin{equation}
\label{pdtVic2}
\forall \alpha \in [\overline \alpha ^*,2), \forall m \geq 1, \quad  \int _a ^b \Phi _{\alpha ,m} (x) ^2 \, dx 
\geq \frac{\gamma _0}{8} \kappa_\alpha .
\end{equation}
But we already proved in \eqref{pdtVic} that $(\int _a ^b \Phi _{\alpha ,m} ^2 )_{m\geq 1}$ is uniformly bounded from below when $\alpha \in [1,\overline \alpha ^*]$. Hence \eqref{pdtVic}  and \eqref{pdtVic2} give \eqref{mmmp} and the proof of Proposition \ref{prop-fctpr} is completed. \qed

\bigskip

{\bf Acknowledgements.} The authors wish to thank the referees for their constructive comments.


\end{document}